\documentclass[reqno]{amsart}
\usepackage{amsmath, amsthm, amssymb, color,amsfonts,mathrsfs}
\usepackage{graphicx}
\usepackage{latexsym}
\usepackage{xcolor}
\usepackage{enumerate}
\usepackage{setspace}

\usepackage{framed}
\usepackage{cases}
\allowdisplaybreaks[4]

\newcommand{\R}{\mathbb{R}}
\newcommand{\T}{\mathbb{T}}
\theoremstyle{plain}

\newtheorem{theorem}{Theorem}[section]
\newtheorem{lemma}[theorem]{Lemma}
\newtheorem{proposition}[theorem]{Proposition}

\newtheorem{definition}[theorem]{Definition}
\newtheorem{remark}[theorem]{Remark}
\theoremstyle{remark}
\numberwithin{equation}{section}

\def\v{\varepsilon}

\def\d{\delta}

\def\r{\rho}

\def\i{\infty}
\def\f{\frac}

\usepackage{tikz}

\newcommand{\mI}{\mathcal{I}}
\newcommand{\mT}{\mathcal{T}}

\author[S. Jiang]{Song Jiang}
\address[S. Jiang]{Institute of Applied Physics and Computational Mathematics, Beijing 100088, China}
\email{jiang@iapcm.ac.cn}

\author[Q.C. Ju]{Qiangchang Ju}
\address[Q.C. Ju]{Institute of Applied Physics and Computational Mathematics, Beijing 100088, China}
\email{ju\_qiangchang@iapcm.ac.cn}

\author[J. Ouyang]{Jing Ouyang}
\address[J. Ouyang]{Institute of Applied Physics and Computational Mathematics, Beijing 100088, China}
\email{ouyangjing@amss.ac.cn}

\date{\today}
\begin{document}
\title[Diffusive and Vanishing-Scattering Limits]{Diffusive and Vanishing-Scattering Limits for a Radiative Heat Transfer System with Boundary Layers}
	
   \begin{abstract}
   	We study the diffusive and vanishing-scattering limits of a
   	steady radiative heat transfer system with isotropic scattering
   	in a periodic slab. Here $\v>0$ denotes the dimensionless mean
   	free path, and $\kappa\in[0,1)$ is the normalized scattering
   	coefficient, with $1-\kappa$ representing the absorption
   	coefficient. Under positivity, finite-order regularity, and
   	spectral assumptions, we establish high-order tangential
   	derivative estimates for the boundary-layer profiles and
   	their far-field temperatures, together with corresponding
   	spatial derivative estimates for the interior profiles.
   	The dependence of these estimates on $1-\kappa$ is tracked
   	explicitly. These estimates close the recursive construction
   	of the high-order interior profiles and both lower and upper boundary layers
   	and ensure exact boundary matching. For each fixed $\kappa<1$,
   	we then justify the diffusive limit as $\v\to0$ and establish
   	an $O(\v)$ approximation of the full solution in $L^\infty$
   	by the leading interior and boundary-layer profiles,
   	quantifying the dependence of the error bound and the
   	admissible range of $\v$ on $1-\kappa$.
   	Under uniform hypotheses and an additional spectral condition,
   	we further prove that the leading profiles differ from
   	their non-scattering counterparts by $O(\kappa)$.
   	This yields an $O(\v+\kappa)$ approximation of the full
   	solution by the non-scattering leading profiles in the
   	joint limit $(\v,\kappa)\to(0,0)$.
   \end{abstract}

\maketitle

\setcounter{tocdepth}{2}
\tableofcontents
	\section{Introduction}\label{int}
	\subsection{Problem statement}
		We study the following radiative heat transfer system:
	\begin{align}
		\begin{cases}
			\v\beta\cdot\nabla I^\v+\kappa(I^\v-\bar{I}^\v)=(1-\kappa)((T^\v)^{4}-I^\v),\\
			-\v^{2}\Delta T^\v=(1-\kappa)4\pi(\bar{I}^\v-(T^\v)^{4}),
		\end{cases}
		\label{1.1}
	\end{align}
	where $I^\v(x,\beta)$ denotes the radiation density and $T^\v(x)$ denotes the temperature. Here $x\in\Omega=\{(x_{1},x_{2},x_{3}):(x_{1},x_{2})\in\T^{2},\,x_{3}\in[0,1]\}$ and $\beta\in\mathbb{S}^{2}$. The term $\kappa(I-\bar I)$ shows isotropic scattering effect
	while $(1-\kappa)((T^\v)^4-I^\v)$ represents nonlinear
	emission and absorption. The parameters satisfy $0<\v\ll1$ and $\kappa\in[0,1)$. The angular average of $I^\v$ is defined by
	\begin{align}
		\bar{I}^\v=\frac{1}{4\pi}\int_{\mathbb{S}^{2}}I^\v(x,\beta)\,\mathrm{d}\beta.
		\label{1.2}
	\end{align}
Let $\gamma=\partial\Omega\times\mathbb S^2$ denote the phase
	boundary, and let $\vec n(x)$ be the outward unit normal to
	$\partial\Omega$. We define the outgoing, incoming, and grazing
	sets by
	\begin{align*}
		\begin{aligned}
			\gamma_+=\{(x,\beta):x\in\partial\Omega, \vec{n}(x)\cdot\beta>0\},\\
			\gamma_-=\{(x,\beta):x\in\partial\Omega, \vec{n}(x)\cdot\beta<0\},\\
			\gamma_0=\{(x,\beta):x\in\partial\Omega, \vec{n}(x)\cdot\beta=0\}.
		\end{aligned}
	\end{align*} 
	 We impose an inflow condition for the radiation and a Dirichlet condition for the temperature:
	\begin{align}
		I^\v(x,\beta)|_{\gamma_{-}}=I_{b}(x,\beta),\quad T^\v(x)|_{\partial\Omega}=T_{b}(x).
		\label{1.3}
	\end{align}
	The incoming data satisfy
	\begin{equation}\label{ue}
		\begin{aligned}
			I_b&=T_b^4
			&&\text{on }\{x_3=1,\ \beta_3<0\},\\
			I_b&\not\equiv T_b^4
			&&\text{on }\{x_3=0,\ \beta_3>0\}.
		\end{aligned}
	\end{equation}
	Throughout the paper, we assume that there are constants $0<a_1\leq a_2$ such that
	\begin{align}\label{bb}
		a_1\leq T_b\leq a_2\quad\text{on }\partial\Omega,\qquad
		0\leq I_b\leq a_2^4\quad\text{on }\gamma_-.
	\end{align}
	
	Under the spectral assumptions stated below, we study the
	diffusive limit $\v\to0$ and the joint diffusive and
	vanishing-scattering limit $(\v,\kappa)\to(0,0)$.
	\subsection{Related works}
	Radiative heat transfer systems typically couple a radiative
	transfer equation with a heat conduction equation for the
	material temperature. For the physical background, we refer
	to the monographs by Modest and Pomraning
	\cite{Modest,Pomraning}. Larsen, Pomraning and Badham \cite{Larsen-1983} studied a
	time-dependent, frequency-dependent radiative transfer model. Under well-prepared initial and boundary conditions, they used a formal asymptotic expansion to derive the Rosseland diffusion approximation and the associated boundary layer
	equations.
	
To the best of our knowledge, the rigorous mathematical study of radiative transfer equations dates back to \cite{BGP-1987}, where Bardos, Golse and Perthame used nonlinear
	semigroup theory to establish well-posedness under suitable
	assumptions on the coefficients and justified the Rosseland
	approximation in $L^1$.
	For a grey radiative transfer model, Bardos, Golse, Perthame,
	and Sentis \cite{BGP-1988} subsequently proved existence of
	solutions and convergence to a degenerate nonlinear diffusion
	equation, allowing the opacity coefficient to be singular
	without requiring it to be monotone.
	
	The glass cooling model couples radiative transfer with heat
	conduction and typically neglects scattering; see Thommes
	et al.\ \cite{gl-1} for a detailed derivation. Klar and Siedow \cite{KS-1998} developed domain decomposition methods coupling radiative transfer and diffusion models, with applications to glass manufacturing. In the frequency-independent setting, Klar and Schmeiser
	\cite{Klar-2000} justified the equilibrium diffusion limit
	and studied the associated nonlinear Milne problem.  In 2022, Ghattassi, Huo, and Masmoudi \cite{GHM-2022} proved global
	existence of weak solutions and justified the diffusive limit
	under well-prepared initial and boundary conditions.
	Under additional regularity assumptions on the limiting
	solution, they also obtained a convergence rate by the
	relative entropy method.
	
When the boundary data do not satisfy $I_b=T_b^4$, the
diffusion approximation requires a boundary layer correction.
Ghattassi, Huo, and Masmoudi \cite{GHM-2023-JMPA} constructed
a Hilbert expansion with boundary layers in a periodic slab
and proved a quantitative $L^\infty$ diffusion approximation
under a spectral assumption.
The corresponding stability theory for the nonlinear Milne
problem and its linearization was developed in \cite{GHM-2023-ARMA}. Curved boundaries introduce additional effects in the boundary layer equations. Following the geometric correction approach developed by Wu and Guo~\cite{Guo-2015} for neutron transport, Ghattassi, Huo, and Masmoudi~\cite{GHM-2025} constructed geometrically corrected boundary layers in a two-dimensional disk and obtained an $L^\infty$ diffusion approximation. Related regularity estimates and diffusive limits in convex domains were established in~\cite{Guo-2017,Wu-2020,Wu-2021}. 

	
The expansion in~\cite{GHM-2023-JMPA} requires tangential derivatives of lower-order boundary layers in the Milne sources and boundary traces of normal derivatives of the interior profiles in the Taylor coefficients.
Building on this construction and the Milne stability theory in~\cite{GHM-2023-ARMA}, we establish the derivative estimates needed to close the coupled construction for the system with scattering. Starting from finite-order boundary regularity and the stated spectral assumptions, we control the boundary layers, their far-field temperatures, and the interior profiles in the derivative ranges required by the expansion. The resulting estimates yield the $L^\infty$ diffusion approximation in Theorem~\ref{thm1}, with explicit dependence of both the error bound and the admissible range of $\v$ on $1-\kappa$. Under uniform data and spectral bounds, the approximation is uniform for $0\leq\kappa\leq\kappa_0<1$.

We further compare the leading profiles of the scattering and non-scattering systems. Under the stated spectral assumptions, their difference is $O(\kappa)$, which yields the $O(\v+\kappa)$ approximation in Theorem~\ref{thm1.5} by the non-scattering leading interior profile and its boundary layer correction.

\subsection{Hilbert expansion}
\subsubsection{Interior expansion} Away from the boundary, we seek the formal expansions
\begin{align}\label{3.1}
	I^\v(x,\beta)\sim\sum_{k=0}^\infty \v^kI_k(x,\beta), \quad T^\v(x)\sim \sum_{k=0}^\infty \v^kT_k(x).
\end{align}
Substituting \eqref{3.1} into \eqref{1.1}, one has 
\begin{align}
	\begin{cases}
		\v\beta\cdot\nabla\left(\sum\limits_{k=0}^\infty \v^kI_k\right)+\kappa\left(\sum\limits_{k=0}^\infty \v^kI_k-\sum\limits_{k=0}^\infty \v^k\bar{I}_k\right)=(1-\kappa)\left((\sum\limits_{k=0}^\infty \v^kT_k)^4-\sum\limits_{k=0}^\infty \v^kI_k\right),\\
		-\v^2\Delta\left(\sum\limits_{k=0}^\infty \v^kT_k\right)=(1-\kappa)4\pi\left[\sum\limits_{k=0}^\infty \v^k\bar{I}_k-\left((\sum\limits_{k=0}^\infty \v^kT_k)^4\right)\right].
	\end{cases}
\end{align}
Comparing equal powers of $\v$, we obtain
\begin{align}\label{3.3}
	\begin{split}
	\v^0&:\begin{cases}\kappa(I_0-\bar{I}_0)=(1-\kappa)(T_0^4-I_0),\\
0=(1-\kappa)4\pi(\bar{I}_0-T_0^4),\end{cases}\\
	\v^1&:\begin{cases}\beta\cdot\nabla I_0+\kappa(I_1-\bar{I}_1)=(1-\kappa)(4T_0^3T_1-I_1),\\
 0=(1-\kappa)4\pi(\bar{I}_1-4T_0^3T_1),\end{cases}\\
	\v^2&:
	\begin{cases}\beta\cdot\nabla I_1+\kappa(I_2-\bar{I}_2)=(1-\kappa)(\mathcal{C}(T,2)-I_2),\\
	-\Delta T_0=(1-\kappa)4\pi(\bar{I}_2-\mathcal{C}(T,2)),\end{cases}\\
	\vdots\\
	\v^k&:
	\begin{cases}\beta\cdot\nabla I_{k-1}+\kappa(I_k-\bar{I}_k)=(1-\kappa)[\mathcal{C}(T,k)-I_k],\\
	-\Delta T_{k-2}=(1-\kappa)4\pi[\bar{I}_k-\mathcal{C}(T,k)],\end{cases}
	\end{split}
\end{align}
where
\begin{align*}
\mathcal C(T,k)
:=\sum_{\substack{i+j+l+q=k\\i,j,l,q\geq0}}T_iT_jT_lT_q.
\end{align*}

The equations at orders $\v^0$ and $\v^1$ give
\begin{align}\label{3.4}
	I_0=T_0^4,\quad \text{and}\quad  I_1=4T_0^3T_1-\beta\cdot\nabla I_0.
\end{align}
Integrating the radiation equation in~\eqref{3.3} with $k=2$
over $\mathbb S^2$, and using the temperature equation
at the same order together with~\eqref{3.4}, we obtain
\[
\Delta T_0
=\int_{\mathbb S^2}\beta\cdot\nabla I_1\,\mathrm d\beta
=-\int_{\mathbb S^2}(\beta\cdot\nabla)^2T_0^4\,\mathrm d\beta
=-\frac{4\pi}{3}\Delta T_0^4.
\]
Hence
\begin{align}\label{3.5}
	\Delta\left(T_0+\frac{4\pi}{3}T_0^4\right)=0.
\end{align}

Set $T_{-1}:=0$. Eliminating $\bar I_k$ from~\eqref{3.3}
and using~\eqref{3.4} for $k=1$, we obtain, for $k\geq1$,
\begin{align}\label{3.9}
	I_k=-\beta\cdot\nabla I_{k-1}
	-\frac{\kappa}{4\pi(1-\kappa)}\Delta T_{k-2}
	+\mathcal C(T,k).
\end{align}
Integrating the radiation equation at order $\v^{k+2}$
in~\eqref{3.3} over $\mathbb S^2$ and using the temperature
equation at the same order, we obtain, for $k\geq0$,
\begin{align}\label{3.6}
	\Delta T_k
	=\int_{\mathbb S^2}\beta\cdot\nabla I_{k+1}\,\mathrm d\beta.
\end{align}
Substituting~\eqref{3.9} with $k$ replaced by $k+1$
into~\eqref{3.6}, the two terms independent of $\beta$
vanish after angular integration because
$\int_{\mathbb S^2}\beta\,\mathrm d\beta=0$. Thus
\begin{align}\label{3.8}
	\Delta T_k
	=-\int_{\mathbb S^2}(\beta\cdot\nabla)^2I_k\,\mathrm d\beta.
\end{align}
Substituting~\eqref{3.9} into~\eqref{3.8} then gives,
for $k\geq1$,
\begin{align}\label{3.10}
	\Delta T_k+\frac{4\pi}{3}\Delta\mathcal C(T,k)
	=\int_{\mathbb S^2}(\beta\cdot\nabla)^3I_{k-1}\,\mathrm d\beta
	+\frac{\kappa}{3(1-\kappa)}\Delta^2T_{k-2}.
\end{align}
Since $\mathcal C(T,k)-4T_0^3T_k$ depends only on
$T_0,\ldots,T_{k-1}$, equation~\eqref{3.10} is linear
in $T_k$ once the lower-order profiles are known.
The boundary conditions for~\eqref{3.5} and~\eqref{3.10}
will be determined by matching with the lower and upper
boundary layers.

\subsubsection{Lower boundary layer expansion}
Denote the boundary layer as $(\mathcal{I}^\v,\mathcal{T}^\v)$ and 
\begin{align}
	\mathcal{I}^\v\sim\sum\limits_{k=0}^\infty\v^k\mathcal{I}_k,\quad\mathcal{T}^\v\sim\sum\limits_{k=0}^\infty\v^k\mathcal{T}_k.
\end{align}
Substituting $(I^\v+\mI^\v,T^\v+\mT^\v)$ into \eqref{1.1} formally gives
\begin{align}\label{3.12}
	\begin{cases}
	\v\beta\cdot\nabla \mI^\v+\kappa(\mI^\v-\bar{\mI}^\v)=(1-\kappa)[(T^\v+\mT^\v)^4-(T^\v)^4-\mI^\v],\\
	-\v^2\Delta\mT^\v=(1-\kappa)4\pi\left\{\bar{\mI}^\v-[(T^\v+\mT^\v)^4-(T^\v)^4]\right\}.	
	\end{cases}
\end{align}
Set $x_3=\v\eta$, so that $\mathrm dx_3=\v\,\mathrm d\eta$,
$\partial_{x_3}=\v^{-1}\partial_\eta$, and
$\partial_{x_3}^2=\v^{-2}\partial_\eta^2$.
For the traces at $x_3=0$, write
\[
 \partial_{x_3}^\ell T_j^0
 :=\partial_{x_3}^\ell T_j(x_{\shortparallel},0),
 \qquad
 \partial_{x_3}^\ell I_j^0
 :=\partial_{x_3}^\ell I_j(x_{\shortparallel},0,\beta).
\]
Accordingly, define the coefficient of total order $i$ by
\begin{align}\label{3.12-1}
 Q_i(\eta,x_{\shortparallel})
 :=\sum_{j+\ell=i}\frac{\eta^\ell}{\ell!}
 \partial_{x_3}^\ell T_j^0.
\end{align}
Fix an integer $0\leq\mathfrak b\leq N$ and assume that
$T_j$ has continuous normal derivatives through order
$\mathfrak b-j+1$ on $\overline\Omega$ for
$0\leq j\leq\mathfrak b$.
Applying Taylor's formula to each $T_j$ at $x_3=0$
and collecting terms of the same total order using
\eqref{3.12-1}, we obtain, for $0\leq\v\eta\leq1$,
\begin{align}\label{tr}
	\sum_{j=0}^{\mathfrak b}\v^jT_j(x_\shortparallel,\v\eta)
	=
	\sum_{i=0}^{\mathfrak b}\v^iQ_i(\eta,x_\shortparallel)
	+\v^{\mathfrak b+1}
	\sum_{j=0}^{\mathfrak b}
	\frac{\eta^{\mathfrak b-j+1}}{(\mathfrak b-j+1)!}
	\partial_{x_3}^{\mathfrak b-j+1}T_j
	(x_\shortparallel,\theta_j\v\eta),
\end{align}
where $\theta_j\in(0,1)$ may depend on
$\mathfrak b,j,x_\shortparallel$, and $\v\eta$.
In particular, the remainder satisfies
\[
 \left|\sum_{j=0}^{\mathfrak b}\v^jT_j(x_\shortparallel,\v\eta)
       -\sum_{i=0}^{\mathfrak b}\v^iQ_i(\eta,x_\shortparallel)\right|
 \leq\v^{\mathfrak b+1}\sum_{j=0}^{\mathfrak b}
       \frac{\eta^{\mathfrak b-j+1}}{(\mathfrak b-j+1)!}
       \|\partial_{x_3}^{\mathfrak b-j+1}T_j\|_{L^\infty(\Omega)}.
\]
Below, $Q=(Q_i)_{i\geq0}$ denotes the sequence of Taylor coefficients.

Using
\[
\v\beta\cdot\nabla
=\beta_3\partial_\eta
+\v\beta_\shortparallel\cdot\nabla_\shortparallel,
\qquad
\v^2\Delta
=\partial_\eta^2+\v^2\Delta_\shortparallel,
\]
we compare powers of $\v$ in~\eqref{3.12}.
The coefficient of the quartic difference at order $\v^k$ is
\begin{align}\label{3.15-0}
	\mathcal P(\mathcal T,k)
	:=\mathcal C(Q+\mathcal T,k)-\mathcal C(Q,k),
	\qquad k\geq0.
\end{align}
At order $\v^0$, we have
$\mathcal P(\mathcal T,0)
=(T_0^0+\mathcal T_0)^4-(T_0^0)^4$.
Thus, setting
\begin{align}\label{3.15-1}
	\mathbf I=(T_0^0)^4+\mathcal I_0,\qquad
	\mathbf T=T_0^0+\mathcal T_0,
\end{align}
the zeroth-order equations become
\begin{align}\label{3.16}
	\begin{cases}
		\beta_3\partial_{\eta}\mathbf{I}+\kappa(\mathbf{I}-\bar{\mathbf{I}})=(1-\kappa)(\mathbf{T}^4-\mathbf{I}),\\
		-\partial_\eta^2\mathbf{T}=(1-\kappa)4\pi(\bar{\mathbf{I}}-\mathbf{T}^4).
	\end{cases}
\end{align}
It follows from \eqref{1.3} that the boundary conditions for \eqref{3.16} are 
\begin{align}\label{3.17}
	\mathbf{I}(x_\shortparallel,0,\beta)|_{\gamma_-}=I_b(x_{\shortparallel},0,\beta),\quad \mathbf{T}(x_{\shortparallel},0)=T_b(x_{\shortparallel},0).
\end{align}
In the following theorem, $x_\shortparallel$ is fixed
and omitted from the notation. The proof is given in Section~\ref{nlm}.
\begin{theorem}[Existence and stability for the nonlinear problem]\label{thm2}
	Under the boundary bounds \eqref{bb}, there exists a bounded non-negative weak solution $(\mathbf{I},\mathbf{T})\in L^2_{\text{loc}}(\mathbb{R}_+\times\mathbb{S}^2)\times L^2_{\text{loc}}(\mathbb{R}_+)$ to \eqref{3.16}--\eqref{3.17}. We have
	\begin{align}\label{4.00}
		0\leq 	\mathbf{I}(\eta,\beta)\leq a_2^4,\quad a_*:=\frac{a_1}{1+4\pi a_2^3}\leq \mathbf{T}(\eta)\leq a_2.
	\end{align}
	Moreover, there exists a constant $T_{0,\infty}\in[a_*,a_2]$ such that for any $0<\alpha<1-\kappa$,
	\begin{align}\label{4.0}
		\begin{aligned}
			e^{\alpha\eta}|\mathbf{I}(\eta,\beta)-T_{0,\infty}^4|\leq  C\left(\frac{\kappa}{1-\kappa-\alpha}+\frac{1-\kappa}{\alpha^{\f12}(1-\kappa-\alpha)^{\f32}}\right)\|I_b-T_b^4\|_{L^\infty(\{\beta_3>0\})},	\quad a.e.,\\
			e^{\alpha\eta}|\mathbf{T}(\eta)-T_{0,\infty}|\leq \frac{C}{\alpha^{\f12}(1-\kappa-\alpha)^{\f12}}\|I_b-T_b^4\|_{L^\infty(\{\beta_3>0\})}, \quad \text{for any }\, \eta\in\R_+.
		\end{aligned}
	\end{align} 
	 Here $C>0$ is independent of $\kappa$, $\alpha$.

In addition, suppose that two sets of boundary data
$(I_{b,i},T_{b,i})$, $i=1,2$, satisfy
\[
 0\leq I_{b,i}\leq a_2^4,
 \qquad
 0<a_1\leq T_{b,i}\leq a_2.
\]
The corresponding half-space solutions can be chosen so that
\begin{align}\label{hshift}
 &\|\mathbf I_1-\mathbf I_2\|_
       {L^\infty(\mathbb R_+\times\mathbb S^2)}
  +\|\mathbf T_1-\mathbf T_2\|_{L^\infty(\mathbb R_+)}
       \nonumber\\
 &\qquad\leq C\Bigl(
       \|I_{b,1}-I_{b,2}\|_{L^\infty(\{\beta_3>0\})}
       +|T_{b,1}-T_{b,2}|
       \Bigr),
\end{align}
where $C$ depends only on $a_2$ and is independent of $\kappa$.
Denoting their far-field temperatures by
\[
 T_{0,\infty,i}:=\lim_{\eta\to\infty}\mathbf T_i(\eta),
 \qquad i=1,2,
\]
we also have
\begin{align}\label{hfarstab}
 |T_{0,\infty,1}-T_{0,\infty,2}|
 \leq C\Bigl(
       \|I_{b,1}-I_{b,2}\|_{L^\infty(\{\beta_3>0\})}
       +|T_{b,1}-T_{b,2}|
       \Bigr).
\end{align}
\end{theorem}

We introduce a non-increasing smooth cutoff function $\chi(z)$ by
\begin{align}\label{cut}
	\chi(z)=\begin{cases}
		1,\, z\in[0,\f14\rho],\\
		0,\, z\geq\f38\rho,\\
		\in(0,1)\, \text{otherwise},
	\end{cases}
\end{align}
where $0<\rho\ll1$ and the cutoff is chosen so that
$|\chi'|\leq C\rho^{-1}$ and $|\chi''|\leq C\rho^{-2}$. We define the truncated leading boundary layer by
\begin{align}\label{Bd}
	\mathscr{I}_0(\eta,x_{\shortparallel},\beta):=\chi(\v\eta)\left(\mathbf{I}(\eta,x_{\shortparallel},\beta)-T_{0,\infty}^4\right), \quad \mathscr{T}_0(\eta,x_{\shortparallel}):=\chi(\v\eta)\left(\mathbf{T}(\eta,x_{\shortparallel})-T_{0,\infty}\right),\end{align} 
	and thus we derive the boundary conditions for \eqref{3.5} as 
\begin{align}\label{Tb}
	T_0(0,x_{\shortparallel})=T_{0,\infty}(x_{\shortparallel}),\quad T_0(1,x_{\shortparallel})=T_b(x_\shortparallel,1).
\end{align}

Set $\mathcal T_{-1}:=0$. For $k\geq1$, define
\[
\phi_k=\mathcal I_k+4(T_0^0)^3T_k^0,
\qquad
g_k=\mathcal T_k+T_k^0.
\]
Comparing the coefficients of $\v^k$ in~\eqref{3.12} and substituting these definitions, we obtain
\begin{align}\label{3.18} 
	\begin{cases}
		\beta_3\partial_\eta\phi_k+\kappa(\phi_k-\bar{\phi}_k)=(1-\kappa)(4\mathbf{T}^3g_k-\phi_k)+\mathcal{S}_{k,1},\\
		-\partial_\eta^2g_k=(1-\kappa)4\pi(\bar{\phi}_k-4\mathbf{T}^3g_k)+\mathcal{S}_{k,2},
	\end{cases}
\end{align}
where 
\begin{align}\label{3.19}
	\begin{split}
		\mathcal{S}_{k,1}=(1-\kappa)(\mathcal{P}(\mT,k)-4\mathbf{T}^3(\mT_k+T_k^0)+4(T_0^0)^3T_k^0)-\beta_{\shortparallel}\cdot \nabla_{\shortparallel}\mathcal{I}_{k-1},\\
		\mathcal{S}_{k,2}=-(1-\kappa)4\pi(\mathcal{P}(\mT,k)-4\mathbf{T}^3(\mT_k+T_k^0)+4(T_0^0)^3T_k^0)+\Delta_{\shortparallel}\mathcal{T}_{k-2}.
	\end{split}
\end{align}
By~\eqref{3.12-1} and~\eqref{3.15-0}, the terms in
$\mathcal P(\mathcal T,k)$ containing the order-$k$
profiles are
\[
4\mathbf T^3\mathcal T_k
+4\bigl[\mathbf T^3-(T_0^0)^3\bigr]T_k^0.
\]
These terms cancel in~\eqref{3.19}, so both source terms depend only on profiles of orders less than $k$.

Matching the coefficients of $\v^k$ in the boundary
conditions~\eqref{1.3} at $x_3=0$, and using the
definitions of $\phi_k$ and $g_k$, we obtain, for $k\geq1$,
\begin{align}\label{lkb}
	\phi_k(0,\beta)|_{\beta_3>0}=4(T_0^0)^3T_k^0-I_k(0,\beta),\quad g_k(0)=0.
\end{align}
By \eqref{3.9}, we get 
\begin{align}
	\phi_k(0,\beta)|_{\beta_3>0}=\beta\cdot\nabla I_{k-1}(0)+\f{\kappa}{4\pi(1-\kappa)}\Delta T_{k-2}(0)-\left(\mathcal{C}(T,k)(0)-4(T_0^0)^3T_k^0\right),
\end{align}
where $\mathcal C(T,k)-4T_0^3T_k$ contains only $T_j$ with $j<k$.

The system \eqref{3.18} has the following standard form, with $S_1=\mathcal S_{k,1}$ and $S_2=-\mathcal S_{k,2}$:
\begin{align}\label{5.1}
	\begin{cases}
		\beta_3\partial_\eta\phi+\kappa(\phi-\bar{\phi})=(1-\kappa)(4\mathbf{T}^3g-\phi)+S_1,\\
		\partial_\eta^2g+(1-\kappa)4\pi(\bar{\phi}-4\mathbf{T}^3g)=S_2,
	\end{cases}
\end{align}
with boundary conditions 
\begin{align}
	\phi(0,\beta)|_{\beta_3>0}=\phi_b(\beta),\quad g(0)=g_b.\label{5.1-2}
\end{align}
Before we state the result for \eqref{5.1}--\eqref{5.1-2}, we propose the following spectral assumption:
\begin{framed}
	\underline{\textbf{Spectral Assumption:}}
	\begin{enumerate}
		\item[(A)] We say that a positive function $T\in C^1(\R_+)$ satisfies the spectral assumption if there exist constants $\vartheta>0$ and $0<M<1$ such that
		\begin{align}\label{asA}
			M\int_0^\infty e^{2\vartheta\eta}(2T^{\f32})^2|\partial_\eta f|^2\,\mathrm{d}\eta\geq 4\int_0^\infty e^{2\vartheta\eta}|\partial_\eta(2T^{\f32})|^2f^2\,\mathrm{d}\eta
		\end{align}
		for every $f\in H^1_{\mathrm{loc}}([0,\infty))$ with $f(0)=0$
        and finite weighted derivative energy
        $\int_0^\infty e^{2\vartheta\eta}4T^3|\partial_\eta f|^2
        \,\mathrm d\eta<\infty$.
    \end{enumerate}
\end{framed}

The spectral assumption gives existence, uniqueness, and exponential decay for the linear problem. 
For $p=2,\infty$, we use the convention
\[
 \|e^{\vartheta\eta}(S_1,S_2)\|_{L^p}
 :=
 \|e^{\vartheta\eta}S_1\|_{L^p(\mathbb R_+\times\mathbb S^2)}
 +\|e^{\vartheta\eta}S_2\|_{L^p(\mathbb R_+)}.
\]
Define
\begin{align}\label{5.01-1}
	\mathcal N_\vartheta:=|g_b|+\frac1{\sqrt{\vartheta}}\Biggl[
	&\frac{|\phi_b|_{L^2(\gamma_-)}}{\sqrt{1-\vartheta}}
	+\frac{\|e^{\vartheta\eta}(S_1,S_2)\|_{L^2}}
	{\sqrt{(1-\vartheta)(1-\kappa-\vartheta)}}\nonumber\\
	&+\frac{\|e^{\vartheta\eta}(\langle S_1\rangle-S_2)\|_{L^2}}
	{\vartheta(1-\vartheta)}
	+\frac{(1-\kappa)|g_b|}
	{(1-\vartheta)(1-\kappa-\vartheta)}\Biggr],
\end{align}
and
\begin{align}\label{5.01-2}
	\mathcal M_\vartheta:=\frac1{1-\kappa-\vartheta}\Bigl[
	&(1-\vartheta)\bigl(|\phi_b|_{L^\infty(\gamma_-)}
	+4a_2^3\mathcal N_\vartheta\bigr)
	+\|e^{\vartheta\eta}S_1\|_{L^\infty}\nonumber\\
	&+4(1-\kappa)\left[a_2^3+
	\frac{3a_2^2}{\sqrt{2\vartheta}}
	\left(\frac{2\pi}{3(1-\kappa-\vartheta)}
	\int_{\beta_3>0}\beta_3(I_b-T_b^4)^2\,\mathrm d\beta\right)^{1/2}
	\right]\mathcal N_\vartheta\Bigr].
\end{align}
\begin{theorem}\label{co}
 Let $0\leq\kappa<1$, and let $(\mathbf I,\mathbf T)$ be the nonlinear Milne profile from Theorem~\ref{thm2}. Suppose $\mathbf T$ satisfies \eqref{asA} with
 $0<\vartheta\leq\tilde b(1-\kappa)$ and $0<M<1$, where $0<\tilde b<1$ is fixed.
 Assume $\phi_b\in L^\infty(\gamma_-)$ and
 $e^{\vartheta\eta}(S_1,S_2)\in L^2\cap L^\infty$.
 Then \eqref{5.1}--\eqref{5.1-2} has a unique bounded solution
 $(\phi,g)\in L^\infty(\mathbb R_+\times\mathbb S^2)
 \times W^{2,\infty}_{\mathrm{loc}}(\mathbb R_+)$.
 It satisfies
 \[
 \|g\|_{L^\infty}\leq C\mathcal N_\vartheta,\qquad
 \|\phi\|_{L^\infty}\leq C\mathcal N_\vartheta
 +|\phi_b|_{L^\infty(\gamma_-)}
 +\frac{\|S_1\|_{L^\infty}}{1-\kappa}.
 \]
 There exists $g_\infty\in\mathbb R$, with
 $|g_\infty|\leq C\mathcal N_\vartheta$, such that
 \begin{align}\label{5.01}
 |g(\eta)-g_\infty|&\leq C\mathcal N_\vartheta e^{-\vartheta\eta},\\
 |\phi(\eta,\beta)-4T_{0,\infty}^3g_\infty|
 &\leq C\mathcal M_\vartheta e^{-\vartheta\eta}
 \quad\text{for }(\eta,\beta)\in \R_+\times\mathbb{S}^2\backslash \gamma_-.\nonumber
 \end{align}
 In addition,
 \begin{align}\label{gg}
  \|e^{\vartheta\eta}\partial_\eta g\|_{L^2(\mathbb R_+)}
  \leq C\sqrt{\vartheta}\,\mathcal N_\vartheta.
 \end{align}

 For $i=1,2$, let $(\phi_i,g_i)$ be the solution corresponding to
 data $(\phi_{b,i},g_{b,i},S_{1,i},S_{2,i})$ satisfying the above
 assumptions, with the same background $\mathbf T$, scattering
 parameter $\kappa$, and weight $\vartheta$. Define
 \[
  g_{\infty,i}:=\lim_{\eta\to\infty}g_i(\eta).
 \]
 Let $\mathcal N_{\vartheta,12}$ and $\mathcal M_{\vartheta,12}$
 be the quantities in~\eqref{5.01-1} and~\eqref{5.01-2},
 respectively, with the data replaced by
 \[
 (\phi_b,g_b,S_1,S_2)
 =
 (\phi_{b,1}-\phi_{b,2},\,
  g_{b,1}-g_{b,2},\,
  S_{1,1}-S_{1,2},\,
  S_{2,1}-S_{2,2}).
 \]
 Then
 \begin{align}\label{ldt}
  &|g_{\infty,1}-g_{\infty,2}|
   +\|g_1-g_2\|_{L^\infty}
   +\|e^{\vartheta\eta}
       [(g_1-g_{\infty,1})-(g_2-g_{\infty,2})]\|_{L^\infty}
   \nonumber\\
  &\qquad
   +\vartheta^{-1/2}
       \|e^{\vartheta\eta}\partial_\eta(g_1-g_2)\|_{L^2}
   \leq C\mathcal N_{\vartheta,12},
 \end{align}
 and
 \begin{align}
  &\|\phi_1-\phi_2\|_{L^\infty}
   \leq C\mathcal N_{\vartheta,12}
   +|\phi_{b,1}-\phi_{b,2}|_{L^\infty(\gamma_-)}
   +\frac{\|S_{1,1}-S_{1,2}\|_{L^\infty}}{1-\kappa},
   \label{ldr}\\
  &\|e^{\vartheta\eta}
    [\phi_1-\phi_2-4T_{0,\infty}^3(g_{\infty,1}-g_{\infty,2})]
    \|_{L^\infty}
   \leq C\mathcal M_{\vartheta,12}.
   \label{ldd}
 \end{align}
 The constant $C$ in all these estimates depends on $a_2$ and
 $\tilde b$ and is independent of $\kappa$ and $\vartheta$.
\end{theorem}

\begin{remark}
	Every positive constant function satisfies~\eqref{asA}.
	For the nonlinear profile, Lemma~\ref{lem-sp} verifies
	this assumption when~\eqref{4.62} holds with $C_b$
	sufficiently small.
	Theorem~\ref{thm3} gives local uniqueness under the
	spectral assumption.
\end{remark}

Once the source and boundary estimates required by Theorem \ref{co} are verified, the theorem gives a unique solution $(\phi_k,g_k)$ to \eqref{3.18} and a constant $\mathfrak{g}_{k,\infty}\in \R$ such that
\begin{align*}
	(g_k,\phi_k)(\eta)
	\longrightarrow
	(\mathfrak g_{k,\infty},4T_{0,\infty}^3\mathfrak g_{k,\infty})
	\qquad (\eta\to\infty).
\end{align*}
The decaying profiles are defined by
\begin{align}\label{cen-1}
	\mathcal T_k=g_k-\mathfrak g_{k,\infty},
	\qquad
	\mathcal I_k=\phi_k-4T_{0,\infty}^3\mathfrak g_{k,\infty}.
\end{align}
We define the approximate $k$-th boundary layer as 
\begin{align}\label{Bdk}
	\mathscr{I}_k(\eta,x_{\shortparallel},\beta):=\chi(\v\eta)\left(\phi_k(\eta,x_{\shortparallel},\beta)-4T_{0,\infty}^3\mathfrak{g}_{k,\infty}\right),\quad \mathscr{T}_k(\eta,x_{\shortparallel}):=\chi(\v\eta)\left(g_k(\eta,x_{\shortparallel})-\mathfrak{g}_{k,\infty}\right).
\end{align}

\subsubsection{Upper boundary layer expansion}\label{ubl}
We write $f|_j=f(x_\shortparallel,j)$, retaining the angular variable
for radiation traces. In upper boundary-layer equations and estimates,
$T_b$ and $I_b$ are evaluated at $x_3=1$.
Since $T_0|_1=T_b$ and $I_0=T_0^4$, the upper-boundary condition
in~\eqref{ue} implies that the leading upper layer vanishes.
All negative-index profiles are set to zero. Higher upper layers will
determine the upper Dirichlet values of the same interior hierarchy.
Set $\xi=(1-x_3)/\v$.
The upper coefficients are
$\mathcal T_k^+(\xi,x_\shortparallel)$ and
$\mathcal I_k^+(\xi,x_\shortparallel,\beta)$, with
$\mathcal T_0^+=\mathcal I_0^+=0$.

The normal scaling gives
\[
 |dx_3|=\v\,d\xi,\qquad \partial_{x_3}=-\v^{-1}\partial_\xi,
 \qquad
 \v\beta\cdot\nabla=-\beta_3\partial_\xi+\v\beta_\shortparallel\cdot\nabla_\shortparallel.
\]
The Taylor coefficients of the interior temperature at the upper boundary are
\begin{align}\label{uq}
 Q_k^+(\xi,x_\shortparallel)
 =\sum_{j+\ell=k}\frac{(-\xi)^\ell}{\ell!}
             \partial_{x_3}^\ell T_j(x_\shortparallel,1).
\end{align}
Under the regularity assumptions used to derive~\eqref{tr}, Taylor's formula at $x_3=1$ and regrouping by total order using~\eqref{uq} give, for $0\leq\v\xi\leq1$,
\begin{align}\label{utr}
 &\sum_{j=0}^{\mathfrak b}\v^jT_j(x_\shortparallel,1-\v\xi)
 =\sum_{k=0}^{\mathfrak b}\v^kQ_k^+(\xi,x_\shortparallel)\nonumber\\
 &\quad+\v^{\mathfrak b+1}\sum_{j=0}^{\mathfrak b}
   \frac{(-\xi)^{\mathfrak b-j+1}}{(\mathfrak b-j+1)!}
   \partial_{x_3}^{\mathfrak b-j+1}T_j
       (x_\shortparallel,1-\theta_j^+\v\xi).
\end{align}
Here each $\theta_j^+\in(0,1)$ may depend on
$\mathfrak b,j,x_\shortparallel$, and $\v\xi$.
In particular, the remainder satisfies
\[
 \left|\sum_{j=0}^{\mathfrak b}\v^jT_j(x_\shortparallel,1-\v\xi)
       -\sum_{k=0}^{\mathfrak b}\v^kQ_k^+(\xi,x_\shortparallel)\right|
 \leq\v^{\mathfrak b+1}\sum_{j=0}^{\mathfrak b}
       \frac{\xi^{\mathfrak b-j+1}}{(\mathfrak b-j+1)!}
       \|\partial_{x_3}^{\mathfrak b-j+1}T_j\|_{L^\infty(\Omega)}.
\]
Thus $Q_0^+=T_b$ and $Q_1^+=T_1|_1-\xi\partial_{x_3}T_0|_1$.

In this notation, $T=(T_j)$, $Q^+=(Q_j^+)$, and
$\mathcal T^+=(\mathcal T_j^+)$ denote coefficient sequences, and
$Q^++\mathcal T^+$ is their termwise sum.
For $1\leq k\leq N$, the upper quartic coefficients and sources are
\begin{align}\label{up}
 \mathcal P^+(\mathcal T^+,k)
 &=\mathcal C(Q^++\mathcal T^+,k)-\mathcal C(Q^+,k),\nonumber\\
 \mathcal S_{k,1}^+
 &=(1-\kappa)\bigl[\mathcal P^+(\mathcal T^+,k)
                         -4T_b^3\mathcal T_k^+\bigr]
        -\beta_\shortparallel\cdot\nabla_\shortparallel\mathcal I_{k-1}^+,\nonumber\\
 \mathcal S_{k,2}^+
 &=-4\pi(1-\kappa)\bigl[\mathcal P^+(\mathcal T^+,k)
                         -4T_b^3\mathcal T_k^+\bigr]
        +\Delta_\shortparallel\mathcal T_{k-2}^+.
\end{align}
Substituting the interior and upper boundary layer expansions
into~\eqref{1.1}, subtracting the interior equations, and
collecting terms of order $\v^k$, we obtain
\begin{align}\label{uus}
 -\beta_3\partial_\xi\mathcal I_k^+
       +\mathcal I_k^+-\kappa\overline{\mathcal I_k^+}
 &=(1-\kappa)\mathcal P^+(\mathcal T^+,k)
       -\beta_\shortparallel\cdot\nabla_\shortparallel\mathcal I_{k-1}^+,
       \nonumber\\
 -\partial_\xi^2\mathcal T_k^+
 &=(1-\kappa)\langle\mathcal I_k^+
                         -\mathcal P^+(\mathcal T^+,k)\rangle
       +\Delta_\shortparallel\mathcal T_{k-2}^+.
\end{align}
Since $Q_0^+=T_b$ and $\mathcal T_0^+=0$, separating the terms with
one index equal to $k$ yields
\begin{align}\label{ust}
 \mathcal P^+(\mathcal T^+,k)-4T_b^3\mathcal T_k^+=\sum_{\substack{i_1+\cdots+i_4=k\\0\leq i_\nu<k}}
 \left[\prod_{\nu=1}^4(Q_{i_\nu}^++\mathcal T_{i_\nu}^+)-\prod_{\nu=1}^4Q_{i_\nu}^+\right].
\end{align}
In particular, $\mathcal S_{k,1}^+$ and $\mathcal S_{k,2}^+$ depend
only on profiles of orders less than $k$.

For each $1\leq k\leq N$, let
$g_k^+(\xi,x_\shortparallel)$ and
$\phi_k^+(\xi,x_\shortparallel,\beta)$ solve the linear Milne problem
\begin{align}\label{uks}
 \begin{cases}
 -\beta_3\partial_\xi\phi_k^++\phi_k^+-\kappa\bar\phi_k^+
       =(1-\kappa)4T_b^3g_k^++\mathcal S_{k,1}^+,\\
 -\partial_\xi^2g_k^+
       =(1-\kappa)\langle\phi_k^+-4T_b^3g_k^+\rangle
    +\mathcal S_{k,2}^+,\\
 g_k^+(0)=0,
 \end{cases}
\end{align}
with incoming boundary data
\begin{align}\label{uki}
 \phi_k^+(0,\beta)
 =\left[\beta\cdot\nabla I_{k-1}
       +\frac{\kappa}{4\pi(1-\kappa)}\Delta T_{k-2}
       -\bigl(\mathcal C(T,k)-4T_0^3T_k\bigr)\right]\Big|_1,
 \quad\beta_3<0.
\end{align}
The interior relation~\eqref{3.9} identifies
\begin{align}\label{uim}
 \phi_k^+(0,\beta)=(4T_b^3T_k-I_k)|_1,
 \qquad \beta_3<0.
\end{align}

Under the reflection $\beta_3\mapsto-\beta_3$, this is the linear
Milne problem~\eqref{5.1} with coefficient $4T_b^3$ and sources
$(S_1,S_2)=(\mathcal S_{k,1}^+,-\mathcal S_{k,2}^+)$. Whenever the incoming datum is bounded and
$e^{\vartheta\xi}(\mathcal S_{k,1}^+,\mathcal S_{k,2}^+)
\in L^2\cap L^\infty$ for some $0<\vartheta<1-\kappa$,
Theorem~\ref{co} gives a unique bounded solution for \eqref{uks}. 
There is a far-field temperature
$\mathfrak g_{k,\infty}^+(x_\shortparallel)$ such that
\begin{align}\label{ukf}
 (g_k^+,\phi_k^+)(\xi)
 \longrightarrow
 (\mathfrak g_{k,\infty}^+,4T_b^3\mathfrak g_{k,\infty}^+)
 \qquad (\xi\to\infty).
\end{align}
The decaying profiles are defined by
\begin{align}\label{cen}
 \mathcal T_k^+=g_k^+-\mathfrak g_{k,\infty}^+,
 \qquad
 \mathcal I_k^+=\phi_k^+-4T_b^3\mathfrak g_{k,\infty}^+.
\end{align}
Using the cutoff $\chi$ from~\eqref{cut}, with the same
width $0<\rho<1$, we define the truncated upper boundary
layers by
\begin{align}\label{ukc}
	\mathscr T_k^+(x)
	=\chi(1-x_3)\mathcal T_k^+((1-x_3)/\v,x_\shortparallel),\quad 
	\mathscr I_k^+(x,\beta)
	=\chi(1-x_3)\mathcal I_k^+((1-x_3)/\v,x_\shortparallel,\beta).
\end{align}

The lower cutoff profiles are $\mathscr T_k=\chi(\v\eta)\mathcal T_k$
and $\mathscr I_k=\chi(\v\eta)\mathcal I_k$, where $\eta=x_3/\v$.
Their supports lie in $x_3\leq3\rho/8$, while the upper supports lie
in $x_3\geq1-3\rho/8$; these are disjoint since $\rho<1$.

For $1\leq k\leq N$, separating $4T_0^3T_k$ from
$\mathcal C(T,k)$ in~\eqref{3.10} and prescribing the
two Milne far-field temperatures as boundary values,
we obtain
\begin{align}\label{ikm}
	\Delta\left(T_k+\frac{16\pi}{3}T_0^3T_k\right)
	&=-\frac{4\pi}{3}\Delta\bigl(\mathcal C(T,k)-4T_0^3T_k\bigr)\nonumber\\
	&\quad+\left\langle(\beta\cdot\nabla)^3I_{k-1}\right\rangle
	+\frac{\kappa}{3(1-\kappa)}\Delta^2T_{k-2},\nonumber\\
	T_k|_{x_3=0}&=\mathfrak g_{k,\infty},\qquad
	T_k|_{x_3=1}=\mathfrak g_{k,\infty}^+.
\end{align}
Since $1+\frac{16\pi}{3}T_0^3>0$, this is a Dirichlet
Poisson problem for
$(1+\frac{16\pi}{3}T_0^3)T_k$ and determines $T_k$
uniquely. The radiation profile $I_k$ is then determined by~\eqref{3.9}.
\subsubsection{Composite approximation}
Define the composite approximation by

\begin{align}\label{sc}
 T^a&=\sum_{k=0}^N\v^k(T_k+\mathscr T_k)
                   +\sum_{k=1}^N\v^k\mathscr T_k^+,\nonumber\\
 I^a&=\sum_{k=0}^N\v^k(I_k+\mathscr I_k)
                   +\sum_{k=1}^N\v^k\mathscr I_k^+.
\end{align}

For the profiles constructed in Section~\ref{prf},
the composite approximation satisfies the exact boundary
matching~\eqref{eb}, as proved in Proposition~\ref{mat}.

\subsubsection{Remainder equations}
For the composite approximation~\eqref{sc}, define the residuals
\begin{align}\label{rd}
 R_1&=\v\beta\cdot\nabla I^a+I^a-\kappa\bar I^a
                                          -(1-\kappa)(T^a)^4,\nonumber\\
 R_2&=-\v^2\Delta T^a-(1-\kappa)\langle I^a-(T^a)^4\rangle.
\end{align}
For a solution $(I^\v,T^\v)$ of~\eqref{1.1} with the boundary conditions~\eqref{1.3}, define
\[
 I_R^\v=I^\v-I^a,\qquad T_R^\v=T^\v-T^a.
\]
Subtracting~\eqref{rd} from~\eqref{1.1} and using
the boundary matching~\eqref{eb}, we obtain
\begin{align}\label{re}
 \begin{cases}
 \v\beta\cdot\nabla I_R^\v+I_R^\v-\kappa\bar I_R^\v
   =(1-\kappa)\bigl[(T^a+T_R^\v)^4-(T^a)^4\bigr]-R_1,\\
 -\v^2\Delta T_R^\v
   =(1-\kappa)\bigl\langle I_R^\v
                         -[(T^a+T_R^\v)^4-(T^a)^4]\bigr\rangle-R_2,\\
 T_R^\v|_{\partial\Omega}=0,\qquad I_R^\v|_{\gamma_-}=0.
 \end{cases}
\end{align}

The exact decomposition of these residuals and their estimates are
given in Section~\ref{res}. 
\subsection{Main results}

For the main convergence result, fix $N=6$ and let $s_0\geq23$. We assume that the boundary data satisfy
\begin{align}\label{pr}
 \sum_{|w|\leq s_0}\left(
 |D^w_{x_{\shortparallel}}I_b|_{L^\infty(\gamma_-)}
 +|D^w_{x_{\shortparallel}}T_b|_{L^\infty(\partial\Omega)}
 \right)\leq C_{\mathrm{data}}<
\infty.
\end{align}
\begin{theorem}\label{thm1}
Assume the hypotheses of Proposition~\ref{prop} with $N=6$.
Let $(I_k,T_k)$, $(\mathscr I_k,\mathscr T_k)$, and
$(\mathscr I_k^+,\mathscr T_k^+)$ be the interior profiles and the
truncated lower and upper boundary layers constructed there.
There are exponents $0<K_6\leq\widehat K_6$ and constants
$C,\rho_*,\widehat\varepsilon_*>0$ such that, with
\[
 \rho=\rho_*(1-\kappa)^5,\qquad
 0<\v\leq\widehat\varepsilon_*(1-\kappa)^{\widehat K_6},
\]
there exists a unique solution
$(I^\v,T^\v)\in L^\infty(\Omega\times\mathbb S^2)\times L^\infty(\Omega)$
to~\eqref{1.1} with boundary conditions~\eqref{1.3}, satisfying
\begin{align}\label{1.49}
 &\|I^\v-T_0^4-\mathscr I_0\|_{L^\infty(\Omega\times\mathbb S^2)}
 +\|T^\v-T_0-\mathscr T_0\|_{L^\infty(\Omega)}
 \leq \frac{C\v}{(1-\kappa)^{K_6}}.
\end{align}
The exponents are independent of $\v,\kappa$. 
\end{theorem}
\begin{remark}
	The exact definitions of $K_6$ and $\hat{K}_6$ are given in\eqref{KN} and \eqref{Khat}.
\end{remark}

We next consider the joint limit $(\v,\kappa)\to(0,0)$.
Let $(\mathbf{I}_{\mathrm{ns}},\mathbf{T}_{\mathrm{ns}})$ be the selected solution to the non-scattering boundary layer problem
\begin{align}\label{8.1}
	\begin{cases}
		\beta_3\partial_{\eta}\mathbf{I}_{\mathrm{ns}}=(\mathbf{T}_{\mathrm{ns}})^4-\mathbf{I}_{\mathrm{ns}},\\
		-\partial_\eta^2\mathbf{T}_{\mathrm{ns}}=4\pi(\bar{\mathbf{I}}_{\mathrm{ns}}-(\mathbf{T}_{\mathrm{ns}})^4),
	\end{cases}
\end{align}
with boundary condition
\begin{align}\label{8.2}
	\mathbf{I}_{\mathrm{ns}}(x_\shortparallel,0,\beta)|_{\gamma_-}=I_b(x_{\shortparallel},0,\beta),\quad \mathbf{T}_{\mathrm{ns}}(x_{\shortparallel},0)=T_b(x_{\shortparallel},0).
\end{align}
By Theorem~\ref{thm2} at $\kappa=0$, there is a $T_{0,\infty,\mathrm{ns}}$ such that $\lim_{\eta\to\infty} \mathbf{T}_{\mathrm{ns}}(\eta)=T_{0,\infty,\mathrm{ns}}$ and $\lim_{\eta\to\infty} \mathbf{I}_{\mathrm{ns}}(\eta,\beta)=(T_{0,\infty,\mathrm{ns}})^4$. Define
\begin{align}\label{8.2-1}
	\mathscr{I}_{0,\mathrm{ns}}=\chi(\v\eta)(\mathbf{I}_{\mathrm{ns}}-(T_{0,\infty,\mathrm{ns}})^4),\quad \mathscr{T}_{0,\mathrm{ns}}=\chi(\v\eta)(\mathbf{T}_{\mathrm{ns}}-T_{0,\infty,\mathrm{ns}}).
\end{align}
Let $T_{0,\mathrm{ns}}$ be the smooth solution to \eqref{3.5} with the non-scattering boundary condition
\begin{align}\label{8.2-2}
	T_{0,\mathrm{ns}}(0,x_{\shortparallel})=T_{0,\infty,\mathrm{ns}}(x_{\shortparallel}),\quad T_{0,\mathrm{ns}}(1,x_{\shortparallel})=T_b(x_{\shortparallel},1),
\end{align}
and set $I_{0,\mathrm{ns}}=(T_{0,\mathrm{ns}})^4$.
For the comparison of the two nonlinear Milne profiles, define the
secant coefficient
\begin{align}\label{sec}
 a_\kappa:=\mathbf T^3+\mathbf T^2\mathbf T_{\mathrm{ns}}
 +\mathbf T(\mathbf T_{\mathrm{ns}})^2+(\mathbf T_{\mathrm{ns}})^3.
\end{align}
\begin{theorem}\label{thm1.5}
Let $0<\kappa_0<1$. Suppose the hypotheses of Theorem~\ref{thm1}
hold for every $0\leq\kappa\leq\kappa_0$, with the same boundary data,
fixed decay parameters, and constants $a_1,a_2,C_{\mathrm{data}}$
independent of $\kappa$. Assume that the spectral bounds in
Proposition~\ref{prop} hold with a common $M_*<1$.
Assume that $a_\kappa^{1/3}$ satisfies the spectral
assumption~\eqref{asA} with a common weight
$0<\vartheta_*<1-\kappa_0$ and $M\leq M_*<1$, uniformly in
$\kappa$ and $x_\shortparallel$.
    Then, for $\v>0$ sufficiently small, the solution furnished by Theorem \ref{thm1} satisfies
	\begin{align}
		\|I^\varepsilon-(T_{0,\mathrm{ns}})^4-\mathscr I_{0,\mathrm{ns}}\|_{L^\infty(\Omega\times\mathbb S^2)}+\|T^\varepsilon-T_{0,\mathrm{ns}}-\mathscr T_{0,\mathrm{ns}}\|_{L^\infty(\Omega)}
		\le C(\varepsilon+\kappa).
	\end{align}
\end{theorem}
\subsection{Key points} 

$\bullet$ {\it Derivative estimates close the high-order construction.}

We first construct nonlinear solutions corresponding to perturbed
boundary data and compare them to obtain difference estimates
for the leading nonlinear Milne profile
(Lemma~\ref{tstab} and Theorem~\ref{thm2}).
Applying these estimates to tangential difference quotients
establishes the existence of the first tangential weak derivatives
(Lemma~\ref{ffr}). For higher derivatives, we subtract the far-field states and introduce the corrected unknowns $\widehat U_w,\widehat V_w$,
so that the sources contain only lower-order derivatives and
decay exponentially; see~\eqref{6.g} and~\eqref{cl}--\eqref{ct}
in the proof of Lemma~\ref{ffr}.
The linear Milne estimates of Theorem~\ref{co} then yield
the high-order weighted derivative bounds in Lemma~\ref{ffr}.
These bounds also control the derivatives of the far-field states
 and close the coupled construction of the higher-order boundary layers and interior profiles
(Proposition~\ref{prop}).

For $0\leq k\leq N$, Proposition~\ref{prop} gives weighted
tangential estimates for both boundary layers up to order
$s_0-3k$, together with $W^{r+2,q}(\Omega)$ bounds for $T_k$
and $W^{r+1,q}(\Omega)$ bounds for $I_k$, for $0\leq r\leq s_0-3k-2$
and any fixed $3<q<\infty$.
The condition $s_0\geq3N+5$ covers the derivatives required
by the Milne sources, Taylor coefficients, and remainder
estimates. The explicit powers of $(1-\kappa)^{-1}$ in
these bounds determine, through the remainder analysis,
the error bound and the smallness condition on $\v$
in Theorem~\ref{thm1}.

$\bullet$ {\it The upper boundary correction gives exact matching at every retained order.}

The upper equilibrium condition~\eqref{ue} makes the leading
upper layer vanish. At each positive order, the upper Milne far-field temperature determines the interior boundary value $T_k|_{x_3=1}=\mathfrak g_{k,\infty}^+$. Subtracting the upper far field then gives
\[
 (T_k+\mathscr T_k^+)|_{x_3=1}=0,\qquad
 (I_k+\mathscr I_k^+)|_{x_3=1}=0\quad(\beta_3<0),
 \qquad 1\leq k\leq N.
\]
Together with the lower boundary corrections, these identities
give exact boundary matching for the composite approximation, as proved in Proposition~\ref{mat}.  The upper coefficient $4T_b^3\geq4a_1^3$ is independent of $\xi$, so the spectral assumption~\eqref{asA} holds automatically.

%

$\bullet$ {\it Stability of the leading profiles as
	scattering vanishes.}

Under the hypotheses of Theorem~\ref{thm1.5}, subtracting
the non-scattering Milne equations from the scattering
equations gives~\eqref{8.3}, with secant coefficient
$a_\kappa$ and a source of size $O(\kappa)$ in weighted
$L^2$ and $L^\infty$ norms.
The spectral assumption on $a_\kappa^{1/3}$ yields
the $O(\kappa)$ stability estimate~\eqref{8.5} for the nonlinear profiles and their far-field temperatures.
By~\eqref{8.6}--\eqref{8.7}, the same bound holds for the leading interior
profiles. Combining these bounds with Theorem~\ref{thm1}
gives an $O(\v+\kappa)$ approximation by the non-scattering
leading interior and boundary-layer profiles, uniformly
for $0\leq\kappa\leq\kappa_0<1$.
\subsection{Notations}
For $w=(w_1,w_2)\in\mathbb N_0^2$, we write
$|w|=w_1+w_2$ and
$D^w=\partial_{x_1}^{w_1}\partial_{x_2}^{w_2}$.
Angular integrals and averages are defined by$
	\langle f\rangle:=\int_{\mathbb S^2}f(\beta)\,\mathrm d\beta,$ and  $
	\bar f:=\frac1{4\pi}\langle f\rangle.
$
Radiation boundary norms use the flux measure $
|f|_{L^2(\gamma_\pm)}^2
:=\int_{\gamma_\pm}|f|^2|\beta\cdot\vec n|\,
\mathrm dS_x\,\mathrm d\beta.
$
For Milne problems, norms are taken at fixed $x_\shortparallel$, and the boundary measure is $|\beta_3|\,\mathrm d\beta$.
The constant $C$ may change from line to line;
its parameter dependence is specified where used.
\subsection{Organization of the paper}
 Sections~\ref{nlm} and~\ref{lin} analyze the nonlinear and
 linear Milne problems. Section~\ref{prf} constructs the
 matched interior and boundary layer profiles and establishes
 their derivative estimates. Sections~\ref{res} and~\ref{vsl}
 prove the remainder estimate and the vanishing-scattering
 limit, respectively. Appendix~\ref{app:lem5.1} contains
 the proof of Lemma~\ref{lem5.1}.
\section{Nonlinear Milne problem}\label{nlm}
Initially, we fix $x_{\shortparallel}\in\mathbb T^2$ and suppress
this dependence in the notation. We first solve the problem on a bounded
interval. For $B>0$, consider
\begin{align}\label{4.1}
	\begin{cases}
	\beta_3\partial_{\eta}\mathbf{I}^B+\kappa(\mathbf{I}^B-\bar{\mathbf{I}}^B)=(1-\kappa)((\mathbf{T}^B)^4-\mathbf{I}^B),\\
	-\partial_\eta^2\mathbf{T}^B=(1-\kappa)4\pi(\bar{\mathbf{I}}^B-(\mathbf{T}^B)^4),
	\end{cases}
	\qquad (\eta,\beta)\in[0,B]\times\mathbb S^2,
\end{align}
with boundary conditions
\begin{align}\label{4.2}
	\begin{split}
	&\mathbf{I}^B(0,\beta)=I_b(\beta),\quad \beta_3>0,\\
	&\mathbf{I}^B(B,\beta_{\shortparallel},\beta_3)
	 =\mathbf{I}^B(B,\beta_{\shortparallel},-\beta_3),\quad \beta_3<0,\\
	&\mathbf{T}^B(0)=T_b,\qquad \partial_\eta\mathbf{T}^B(B)=0.
	\end{split}
\end{align}
We write $\Omega_B=[0,B]\times\mathbb S^2$. The incoming, outgoing, and grazing parts of the phase boundary are denoted by $\gamma_-$, $\gamma_+$, and $\gamma_0$, respectively. Thus $\gamma_-$ consists of $\{\eta=0,\ \beta_3>0\}$ and $\{\eta=B,\ \beta_3<0\}$, while $\gamma_+$ is obtained by reversing these inequalities.

\subsection{The truncated problem}
\subsubsection{Monotonicity}
We first recall the comparison principle for the second-order equation; see \cite[Lemma 1]{GHM-2023-ARMA}.
\begin{lemma}[\cite{GHM-2023-ARMA}]\label{lem4.1}
	Let $0\leq h_b\leq a_2$. Suppose that $l\in C(\mathbb R)$ is strictly increasing, $l(0)=0$, and $f\in C([0,B])$ satisfies
	\[
	0\leq f(\eta)\leq l(a_2),\qquad \eta\in[0,B].
	\]
	Then the problem
	\begin{align}\label{4.11}
	-\partial_\eta^2h+l(h)=f,\qquad h(0)=h_b,\qquad \partial_\eta h(B)=0,
	\end{align}
	has a unique solution $h\in C^2([0,B])$, and $0\leq h\leq a_2$ on $[0,B]$. Moreover, let $h_i$ solve \eqref{4.11} with source $f_i$ and boundary value $h_{b,i}$, $i=1,2$. If
	\[
	0\leq f_1\leq f_2\leq l(a_2),
	\qquad 0\leq h_{b,1}\leq h_{b,2}\leq a_2,
	\]
	then $0\leq h_1\leq h_2\leq a_2$ on $[0,B]$.
\end{lemma}

The transport equation has an analogous comparison principle.
\begin{lemma}\label{lem4.2}
	Let $\kappa\in[0,1)$, $0\leq h\leq a_2$ on $[0,B]$, and $0\leq\psi_b\leq a_2^4$ on $\{\beta_3>0\}$. Then there exists a unique bounded solution of
	\begin{align}\label{4.3}
	\beta_3\partial_\eta\psi+\psi-\kappa\bar\psi=(1-\kappa)h^4
	\end{align}
	with
	\[
	\psi(0,\beta)=\psi_b(\beta),\quad \beta_3>0,
	\qquad
	\psi(B,\beta_{\shortparallel},\beta_3)
	=\psi(B,\beta_{\shortparallel},-\beta_3),\quad \beta_3<0.
	\]
	Moreover, $0\leq\psi\leq a_2^4$ almost everywhere on $\Omega_B$. If $h\in C([0,B])$ and $\psi_b$ is continuous on $\{\beta\in\mathbb S^2:\beta_3>0\}$, then $\psi$ is continuous on $\Omega_B\setminus\gamma_0$.

	If $h_1\leq h_2$ and $\psi_{b,1}\leq\psi_{b,2}$, then the corresponding solutions satisfy $\psi_1\leq\psi_2$ almost everywhere.
\end{lemma}
\begin{proof}
	We construct a sequence by setting $\psi^0=0$ and, for each $i\geq0$, solving
	\begin{align}\label{ti}
	\begin{cases}
	\beta_3\partial_\eta\psi^{i+1}+\psi^{i+1}
	=\kappa\bar\psi^i+(1-\kappa)h^4,
	&0<\eta<B,\quad\beta\in\mathbb S^2,\\
	\psi^{i+1}(0,\beta)=\psi_b(\beta),&\beta_3>0,\\
	\psi^{i+1}(B,\beta_{\shortparallel},\beta_3)
	=\psi^{i+1}(B,\beta_{\shortparallel},-\beta_3),&\beta_3<0.
	\end{cases}
	\end{align}
	For $\beta_3\neq0$, let $\tau_-(\eta,\beta)$ be the travel time from  $(\eta,\beta)$ to the lower incoming boundary  along the specular characteristic:
	\[
	\tau_-(\eta,\beta)=
	\begin{cases}
	\eta/\beta_3,&\beta_3>0,\\
	(2B-\eta)/|\beta_3|,&\beta_3<0.
	\end{cases}
	\]
	Denote this characteristic by $(X_{\eta,\beta}(s),V_{\eta,\beta}(s))$, $0\leq s\leq\tau_-(\eta,\beta)$; it has one reflection at $\eta=B$ when $\beta_3<0$. Integrating \eqref{ti} along this characteristic, we obtain
	\begin{align}\label{4.6}
	\psi^{i+1}(\eta,\beta)
	&=\psi_b(V_{\eta,\beta}(0))e^{-\tau_-(\eta,\beta)}\nonumber\\
	&\quad+\int_0^{\tau_-(\eta,\beta)}
	\left[\kappa\bar\psi^i(X_{\eta,\beta}(s))
	 +(1-\kappa)h^4(X_{\eta,\beta}(s))\right]
	 e^{-(\tau_-(\eta,\beta)-s)}\,\mathrm ds.
	\end{align}
	The bounds on $h$ and $\psi_b$ imply inductively that
	\[
	0\leq\psi^i\leq a_2^4.
	\]
	Moreover,
	\[
	\|\psi^{i+1}-\psi^i\|_{L^\infty(\Omega_B)}
	\leq\kappa\|\psi^i-\psi^{i-1}\|_{L^\infty(\Omega_B)}.
	\]
	Since $\kappa<1$, the sequence converges in $L^\infty(\Omega_B)$ to the unique solution of \eqref{4.3}. Formula \eqref{4.6} also gives continuity away from $\gamma_0$. Applying the same iteration simultaneously to two ordered pairs of data proves the comparison statement.
\end{proof}

\subsubsection{Existence and uniqueness}
Inspired by \cite{GHM-2023-ARMA}, we use the following notion for the monotone construction.
\begin{definition}\label{def1}
	Let $(\psi,h)\in L^\infty(\Omega_B)\times
	(H^1(0,B)\cap C([0,B]))$, with boundary traces
	$\psi(0,\cdot)$, $\psi(B,\cdot)$, and $\partial_\eta h(B)$.
	We call $(\psi,h)$ a subsolution of \eqref{4.1}--\eqref{4.2}
	if it satisfies the inequalities
	\begin{align}
	\beta_3\partial_\eta\psi+\psi-\kappa\bar\psi
	&\leq(1-\kappa)h^4,\nonumber\\
	-\partial_\eta^2h+(1-\kappa)4\pi h^4
	&\leq(1-\kappa)4\pi\bar\psi,
	\end{align}
	in the sense of distributions on $(0,B)\times\mathbb S^2$ and
	$(0,B)$, respectively, together with the boundary conditions
	\[
	\begin{gathered}
	0\leq\psi(0,\beta)\leq I_b(\beta),\qquad \beta_3>0,\\
	\psi(B,\beta_{\shortparallel},\beta_3)
	=\psi(B,\beta_{\shortparallel},-\beta_3),\qquad \beta_3<0,\\
	0\leq h(0)\leq T_b,\qquad \partial_\eta h(B)=0.
	\end{gathered}
	\]
\end{definition}

\begin{lemma}\label{lem4.3}
	Assume that
	\[
	0\leq I_b\leq a_2^4,\qquad 0<a_1\leq T_b\leq a_2
	\]
	with $a_1,a_2$ as in \eqref{bb}. Then \eqref{4.1}--\eqref{4.2} has a unique bounded solution $(\mathbf I^B,\mathbf T^B)$ satisfying
	\[
	0\leq\mathbf I^B\leq a_2^4\quad\text{a.e. on }\Omega_B,
	\qquad \frac{a_1}{1+4\pi a_2^3}\leq\mathbf T^B\leq a_2\quad\text{on }[0,B].
	\]
	If $I_b$ is continuous on $\{\beta\in\mathbb S^2:\beta_3>0\}$, then $\mathbf I^B$ is continuous on $\Omega_B\setminus\gamma_0$ and $\mathbf T^B\in C^2([0,B])$.
\end{lemma}
\begin{proof}
	The existence and uniqueness proof is similar to \cite[Proof of Theorem 4]{GHM-2023-ARMA}. We ignore it here and obtain a unique solution $(\mathbf I^B,\mathbf T^B)$ satisfying
	\eqref{4.1}--\eqref{4.2} with 
	\begin{align*}
		0\leq\mathbf I^B\leq a_2^4\quad\text{a.e. on }\Omega_B,
		\qquad 0\leq\mathbf T^B\leq a_2\quad\text{on }[0,B]. 
	\end{align*}

	Under the assumption $T_b\geq a_1>0$, we obtain a positive lower bound independent of $B$ and $\kappa$. Taking the angular integral of \eqref{4.1}$_1$ and using \eqref{4.1}$_2$ gives
\[
 \partial_\eta\left(\partial_\eta\mathbf T^B-\langle\beta_3\mathbf I^B\rangle\right)=0.
\]
The reflection condition and $\partial_\eta\mathbf T^B(B)=0$ imply
\begin{align}\label{TIE}
	\partial_\eta\mathbf T^B=\langle\beta_3\mathbf I^B\rangle.
\end{align}
Multiplying \eqref{4.1}$_1$ by $\beta_3$ and integrating over $\mathbb S^2$, we also obtain
\[
 \partial_\eta\langle\beta_3^2\mathbf I^B\rangle
 =-\langle\beta_3\mathbf I^B\rangle.
\]
Consequently,
\begin{align}\label{tmc}
 \mathbf T^B(\eta)+\langle\beta_3^2\mathbf I^B(\eta,\cdot)\rangle
 &=T_b+\langle\beta_3^2\mathbf I^B(0,\cdot)\rangle
 \geq T_b.
\end{align}
Since $\mathbf I^B\geq0$ and $\beta_3^2\leq1$, it follows that
$\mathbf T^B+\langle\mathbf I^B\rangle\geq T_b$.
Using \eqref{4.1}$_2$ and $0\leq\mathbf T^B\leq a_2$, we obtain
\begin{align*}
 &-\partial_\eta^2\mathbf T^B
 +(1-\kappa)(1+4\pi a_2^3)\mathbf T^B\\
 &\quad=(1-\kappa)\left[\langle\mathbf I^B\rangle+\mathbf T^B
 +4\pi\mathbf T^B\bigl(a_2^3-(\mathbf T^B)^3\bigr)\right]
 \geq(1-\kappa)a_1.
\end{align*}
Set $w=\mathbf T^B-a_1/(1+4\pi a_2^3)$. Then
\[
 -w''+(1-\kappa)(1+4\pi a_2^3)w\geq0,
 \qquad w(0)\geq0,\qquad w'(B)=0.
\]
If $w$ had a negative minimum at some $\eta_0\in(0,B)$, then
$w''(\eta_0)\geq0$, whereas the differential inequality gives
$w''(\eta_0)\leq(1-\kappa)(1+4\pi a_2^3)w(\eta_0)<0$.
If the negative minimum were attained at $B$, then $w''(B)<0$ and $w'(B)=0$ would give
\[
 w(B-t)=w(B)+\tfrac12w''(B)t^2+o(t^2)<w(B)
\]
for sufficiently small $t>0$, again a contradiction. Thus
\begin{align}\label{tpt}
 \mathbf T^B(\eta)\geq\frac{a_1}{1+4\pi a_2^3},
 \qquad 0\leq\eta\leq B.
\end{align}
Then Lemma \ref{lem4.3} holds.
\end{proof}

\subsubsection{Dependence on the boundary data}
\begin{lemma}[Stability with respect to the boundary data]\label{tstab}
Let $0\leq\kappa<1$ and $B>0$. Suppose that the two sets of boundary
data $(I_{b,i},T_{b,i})$, $i=1,2$, satisfy
\[
 0\leq I_{b,i}\leq a_2^4,\qquad 0<a_1\leq T_{b,i}\leq a_2,
\]
with the same constants $a_1,a_2$ as in~\eqref{bb}. Let
$(\mathbf I_i^B,\mathbf T_i^B)$ be the corresponding solutions of
\eqref{4.1}--\eqref{4.2} furnished by Lemma~\ref{lem4.3}. Then
\begin{align}\label{tshift}
 &\|\mathbf I_1^B-\mathbf I_2^B\|_{L^\infty(\Omega_B)}
  +\|\mathbf T_1^B-\mathbf T_2^B\|_{L^\infty(0,B)}\nonumber\\
 &\qquad\leq C\Bigl(
       \|I_{b,1}-I_{b,2}\|_{L^\infty(\{\beta_3>0\})}
       +|T_{b,1}-T_{b,2}|\Bigr),
\end{align}
where $C$ depends only on $a_2$ and is independent of $B$ and $\kappa$.
\end{lemma}
\begin{proof}
We first consider ordered boundary data:
$I_{b,1}\leq I_{b,2}$ and $T_{b,1}\leq T_{b,2}$.
To verify that the corresponding solutions are ordered, we apply
the comparison principles in Lemmas~\ref{lem4.1} and~\ref{lem4.2}
to the monotone construction used in Lemma~\ref{lem4.3}.
For clarity, denoting the iterates by
$(\mathbf I_i^{B,n},\mathbf T_i^{B,n})$, we start with
$(\mathbf I_i^{B,0},\mathbf T_i^{B,0})=(0,0)$ and solve
\begin{equation}\label{titer}
 \begin{cases}
 \beta_3\partial_\eta\mathbf I_i^{B,n+1}
  +\mathbf I_i^{B,n+1}-\kappa\overline{\mathbf I_i^{B,n+1}}
  =(1-\kappa)(\mathbf T_i^{B,n})^4,\\
 -\partial_\eta^2\mathbf T_i^{B,n+1}
  +4\pi(1-\kappa)(\mathbf T_i^{B,n+1})^4
  =(1-\kappa)\langle\mathbf I_i^{B,n}\rangle,
 \end{cases}
\end{equation}
with the boundary conditions~\eqref{4.2} and data $I_{b,i},T_{b,i}$
at every step. If the two iterates are ordered at step $n$, then
$(\mathbf T_1^{B,n})^4\leq(\mathbf T_2^{B,n})^4$ and
$\langle\mathbf I_1^{B,n}\rangle\leq\langle\mathbf I_2^{B,n}\rangle$.
Applying Lemma~\ref{lem4.2} to the first equation of~\eqref{titer}
and Lemma~\ref{lem4.1} to the second gives the ordering at step
$n+1$. Induction and passage to the limits
in the construction of Lemma~\ref{lem4.3} yield
\begin{equation}\label{tord}
 \mathbf I_1^B\leq\mathbf I_2^B,\qquad
 \mathbf T_1^B\leq\mathbf T_2^B.
\end{equation}
For bounded incoming data, these inequalities follow by ordered
approximation and the uniqueness in Lemma~\ref{lem4.3}.

It follows from \eqref{TIE} and \eqref{tmc} that
\begin{align}\label{tdflux}
 \partial_\eta(\mathbf T_2^B-\mathbf T_1^B)
 &=\langle\beta_3(\mathbf I_2^B-\mathbf I_1^B)\rangle
\end{align}
and
\begin{align}\label{tdmom}
 \mathbf T_2^B(\eta)-\mathbf T_1^B(\eta)
   +\langle\beta_3^2(\mathbf I_2^B-\mathbf I_1^B)(\eta,\cdot)\rangle
  =T_{b,2}-T_{b,1}
   +\langle\beta_3^2(\mathbf I_2^B-\mathbf I_1^B)(0,\cdot)\rangle.
\end{align}

Using~\eqref{tord} and the upper bound in Lemma~\ref{lem4.3}, we have
\begin{equation}\label{tquart}
 0\leq(\mathbf T_2^B)^4-(\mathbf T_1^B)^4
 \leq4a_2^3(\mathbf T_2^B-\mathbf T_1^B).
\end{equation}
Then 
\begin{align}\label{tdtemp}
 &-\partial_\eta^2(\mathbf T_2^B-\mathbf T_1^B)
 +16\pi(1-\kappa)a_2^3(\mathbf T_2^B-\mathbf T_1^B)\nonumber\\
 &\quad=(1-\kappa)\langle\mathbf I_2^B-\mathbf I_1^B\rangle
 +4\pi(1-\kappa)\Bigl[
 4a_2^3(\mathbf T_2^B-\mathbf T_1^B)
 -\bigl((\mathbf T_2^B)^4-(\mathbf T_1^B)^4\bigr)\Bigr]
 \nonumber\\
 &\quad\geq0,
\end{align}
where the last inequality follows from~\eqref{tord}
and~\eqref{tquart}. The boundary conditions~\eqref{4.2} give
\[
 (\mathbf T_2^B-\mathbf T_1^B)(0)=T_{b,2}-T_{b,1},\qquad
 \partial_\eta(\mathbf T_2^B-\mathbf T_1^B)(B)=0.
\]
Solving the homogeneous equation for the operator
$-\partial_\eta^2+16\pi(1-\kappa)a_2^3$ with these boundary
conditions gives the expression on the right below. Applying the
scalar comparison principle to~\eqref{tdtemp} yields
\begin{equation}\label{tlow}
 \mathbf T_2^B(\eta)-\mathbf T_1^B(\eta)
 \geq(T_{b,2}-T_{b,1})
 \frac{\cosh\bigl(4\sqrt{\pi(1-\kappa)a_2^3}\,(B-\eta)\bigr)}
      {\cosh\bigl(4\sqrt{\pi(1-\kappa)a_2^3}\,B\bigr)}.
\end{equation}
Since both sides of~\eqref{tlow} agree at $\eta=0$,
subtracting their common boundary value, dividing by
$\eta>0$, and letting $\eta\downarrow0$, we obtain
\begin{align}\label{tslope}
	\partial_\eta(\mathbf T_2^B-\mathbf T_1^B)(0)
	&\geq
	-4\sqrt{\pi(1-\kappa)a_2^3}\,
	\tanh\bigl(4\sqrt{\pi(1-\kappa)a_2^3}\,B\bigr)
	(T_{b,2}-T_{b,1})
	\nonumber\\
	&\geq
	-4\sqrt{\pi(1-\kappa)a_2^3}\,
	(T_{b,2}-T_{b,1}),
\end{align}
where we used $T_{b,2}-T_{b,1}\geq0$ and
$0\leq\tanh z\leq1$ for $z\geq0$.

Evaluating~\eqref{tdflux} at $\eta=0$, separating the
incoming and outgoing hemispheres, and using the incoming
condition in~\eqref{4.2} together with~\eqref{tslope},
we obtain
\begin{align}\label{tout}
	&\int_{\beta_3<0}|\beta_3|
	(\mathbf I_2^B-\mathbf I_1^B)(0,\beta)\,\mathrm d\beta
	\nonumber\\
	&\quad=
	\int_{\beta_3>0}\beta_3(I_{b,2}-I_{b,1})\,\mathrm d\beta
	-\partial_\eta(\mathbf T_2^B-\mathbf T_1^B)(0)
	\nonumber\\
	&\quad\leq
	\int_{\beta_3>0}\beta_3(I_{b,2}-I_{b,1})\,\mathrm d\beta
	+4\sqrt{\pi(1-\kappa)a_2^3}\,
	(T_{b,2}-T_{b,1}).
\end{align}
We now return to~\eqref{tdmom}. Dropping its nonnegative angular
moment at $\eta$, taking the supremum over $\eta\in[0,B]$, and
using $\beta_3^2\leq|\beta_3|$ in the boundary moment, we obtain
\begin{align}\label{ttbound}
 &\|\mathbf T_2^B-\mathbf T_1^B\|_{L^\infty(0,B)}\nonumber\\
 &\leq T_{b,2}-T_{b,1}
       +\langle\beta_3^2(\mathbf I_2^B-\mathbf I_1^B)(0,\cdot)\rangle
       \nonumber\\
 &\leq T_{b,2}-T_{b,1}
       +\int_{\beta_3>0}\beta_3(I_{b,2}-I_{b,1})\,\mathrm d\beta
       +\int_{\beta_3<0}|\beta_3|
          (\mathbf I_2^B-\mathbf I_1^B)(0,\beta)\,\mathrm d\beta
       \nonumber\\
 &\leq\bigl(1+4\sqrt{\pi a_2^3}\bigr)(T_{b,2}-T_{b,1})
       +2\pi\|I_{b,2}-I_{b,1}\|_{L^\infty(\{\beta_3>0\})}.
\end{align}
In the last inequality, we substituted~\eqref{tout} and used
$1-\kappa\leq1$ and
$\int_{\beta_3>0}\beta_3\,\mathrm d\beta=\pi$.
Having bounded the temperature difference by~\eqref{ttbound}, we
estimate $\mathbf{I}_2^B-\mathbf{I}_1^B$ which satisfies
\begin{align}\label{ttrans}
 &\beta_3\partial_\eta(\mathbf I_2^B-\mathbf I_1^B)
   +(\mathbf I_2^B-\mathbf I_1^B)
   -\kappa\overline{\mathbf I_2^B-\mathbf I_1^B}\nonumber\\
 &\qquad=(1-\kappa)\bigl((\mathbf T_2^B)^4-(\mathbf T_1^B)^4\bigr),
\end{align}
with incoming value $I_{b,2}-I_{b,1}$ and the reflection condition
at $B$. The constant
\begin{equation}\label{tsuper}
 \max\Bigl\{\|I_{b,2}-I_{b,1}\|_{L^\infty(\{\beta_3>0\})},
             \|(\mathbf T_2^B)^4-(\mathbf T_1^B)^4\|_{L^\infty(0,B)}
      \Bigr\}
\end{equation}
is independent of $\eta$ and $\beta$. 
Applying the comparison
principle of Lemma~\ref{lem4.2} to~\eqref{ttrans} with this
constant supersolution, and then using~\eqref{tquart}, we obtain
\begin{align}\label{tirad}
 \|\mathbf I_2^B-\mathbf I_1^B\|_{L^\infty(\Omega_B)}
 &\leq\|I_{b,2}-I_{b,1}\|_{L^\infty(\{\beta_3>0\})}
       +4a_2^3\|\mathbf T_2^B-\mathbf T_1^B\|_{L^\infty(0,B)}.
\end{align}
Substituting~\eqref{ttbound} into~\eqref{tirad} proves
\eqref{tshift} for ordered boundary data, with a constant
independent of $B$ and $\kappa$.

For general boundary data, define their lower and upper bounds by
\begin{align*}
 I_{b,-}(\beta)&=\min\{I_{b,1}(\beta),I_{b,2}(\beta)\},&
 T_{b,-}&=\min\{T_{b,1},T_{b,2}\},\\
 I_{b,+}(\beta)&=\max\{I_{b,1}(\beta),I_{b,2}(\beta)\},&
 T_{b,+}&=\max\{T_{b,1},T_{b,2}\}.
\end{align*}
The minimum and maximum of the radiation data are taken
pointwise in $\beta$. These data satisfy the same bounds,
so Lemma~\ref{lem4.3} gives corresponding solutions
$(\mathbf I_-^B,\mathbf T_-^B)$ and
$(\mathbf I_+^B,\mathbf T_+^B)$.
Applying~\eqref{tord} to each original data set and its
lower and upper bounds yields
\begin{align}\label{tabs}
	|\mathbf I_1^B-\mathbf I_2^B|
	&\leq\mathbf I_+^B-\mathbf I_-^B,
	\qquad
	|\mathbf T_1^B-\mathbf T_2^B|
	\leq\mathbf T_+^B-\mathbf T_-^B.
\end{align}
Taking the two $L^\infty$ norms in~\eqref{tabs},
adding them, and applying the ordered-data case
of~\eqref{tshift} to the upper and lower solutions
proves~\eqref{tshift}, since
\[
I_{b,+}-I_{b,-}=|I_{b,1}-I_{b,2}|,
\qquad
T_{b,+}-T_{b,-}=|T_{b,1}-T_{b,2}|.
\]
Then we have finished the proof of Lemma \ref{tstab}.
\end{proof}

\subsubsection{Weighted estimate}
\begin{lemma}\label{lem4.5}
	Let $(\mathbf I^B,\mathbf T^B)$ be the solution constructed in Lemma \ref{lem4.3}, and suppose
	\[
	\int_{\beta_3>0}\beta_3(I_b-T_b^4)^2\,\mathrm d\beta<\infty.
	\]
	Then, for every $0\leq\alpha<1-\kappa$ and every $\eta\in[0,B]$,
	\begin{align}\label{4.32}
	&\int_0^\eta e^{2\alpha z}4(\mathbf T^B)^3|\partial_z\mathbf T^B|^2\,\mathrm dz
	+\kappa\|e^{\alpha z}(\mathbf I^B-\bar{\mathbf I}^B)\|_{L^2([0,\eta]\times\mathbb S^2)}^2\nonumber\\
	&\quad +(1-\kappa-\alpha)
	\|e^{\alpha z}(\mathbf I^B-(\mathbf T^B)^4)\|_{L^2([0,\eta]\times\mathbb S^2)}^2\nonumber\\
	&\quad +\frac12\int_{\beta_3<0}|\beta_3|
	(\mathbf I^B(0,\beta)-T_b^4)^2\,\mathrm d\beta
	\leq\frac12\int_{\beta_3>0}\beta_3(I_b-T_b^4)^2\,\mathrm d\beta.
	\end{align}
\end{lemma}
\begin{proof}
	Multiplying \eqref{4.1}$_1$, \eqref{4.1}$_2$ by $\mathbf I^B$ and $(\mathbf T^B)^4$, integrating respectively over $[0,\eta]\times\mathbb S^2$, and $[0,\eta]$, one has 
	\begin{align}\label{4.26}
	&-\int_0^\eta(\mathbf T^B)^4\partial_z^2\mathbf T^B\,\mathrm dz
	+\frac12\int_0^\eta\int_{\mathbb S^2}\beta_3\partial_z(\mathbf I^B)^2\,\mathrm d\beta\,\mathrm dz\nonumber\\
	&\quad+\kappa\|\mathbf I^B-\bar{\mathbf I}^B\|_{L^2([0,\eta]\times\mathbb S^2)}^2
	+(1-\kappa)\|\mathbf I^B-(\mathbf T^B)^4\|_{L^2([0,\eta]\times\mathbb S^2)}^2=0,
	\end{align}
	since
	\[
	\int_{\mathbb S^2}(\mathbf I^B-\bar{\mathbf I}^B)\mathbf I^B\,\mathrm d\beta
	=\int_{\mathbb S^2}|\mathbf I^B-\bar{\mathbf I}^B|^2\,\mathrm d\beta.
	\]
	The remaining estimate are similar to \cite[Section 2.3]{GHM-2023-ARMA}. We just take the constant $1$ as $1-\kappa$ and add a new term about $\kappa\|\mathbf{I}^B-\bar{\mathbf{I}}^B\|_{L^2}$. Then we omit the details here. 
\end{proof}

\subsection{The half-space problem}
We now let $B\to\infty$. 
\begin{lemma}\label{lem4.h}
	Assume $0\leq I_b\leq a_2^4$ and $0<a_1\leq T_b\leq a_2$. Then \eqref{3.16}--\eqref{3.17} has a weak solution
	\[
	(\mathbf I,\mathbf T)\in L^2_{\mathrm{loc}}(\mathbb R_+\times\mathbb S^2)
	\times C_{\mathrm{loc}}(\mathbb R_+)
	\]
	satisfying $0\leq\mathbf I\leq a_2^4$ and $\frac{a_1}{1+4\pi a_2^3}\leq\mathbf T\leq a_2$. For every $0\leq\alpha<1-\kappa$ and $\eta>0$,
	\begin{align}\label{4.54}
	&\int_0^\eta e^{2\alpha z}4\mathbf T^3|\partial_z\mathbf T|^2\,\mathrm dz
	+\kappa\|e^{\alpha z}(\mathbf I-\bar{\mathbf I})\|_{L^2([0,\eta]\times\mathbb S^2)}^2\nonumber\\
	&\quad +(1-\kappa-\alpha)
	\|e^{\alpha z}(\mathbf I-\mathbf T^4)\|_{L^2([0,\eta]\times\mathbb S^2)}^2
	+\frac12\int_{\beta_3<0}|\beta_3|(\mathbf I(0,\beta)-T_b^4)^2\,\mathrm d\beta\nonumber\\
	&\leq\frac12\int_{\beta_3>0}\beta_3(I_b-T_b^4)^2\,\mathrm d\beta.
	\end{align}
	Moreover,
	\begin{align}\label{4.55}
	\partial_\eta\mathbf T(\eta)=\langle\beta_3\mathbf I(\eta,\cdot)\rangle
	\end{align}
	for almost every $\eta>0$. There exists $T_{0,\infty}\in\bigl[\frac{a_1}{1+4\pi a_2^3},a_2\bigr]$ such that, for every $0<\alpha<1-\kappa$,
	\begin{align}\label{4.56-1}
	|\mathbf T(\eta)-T_{0,\infty}|
	&\leq\frac1{\sqrt{2\alpha}}e^{-\alpha\eta}
	\left(\frac{2\pi}{3(1-\kappa-\alpha)}
	\int_{\beta_3>0}\beta_3(I_b-T_b^4)^2\,\mathrm d\beta\right)^{1/2},
	\end{align}
	and
	\begin{align}\label{4.60}
	|\partial_\eta\mathbf T(\eta)|
	&\leq(1-\kappa)\left(\frac{\pi}{\alpha(1-\kappa-\alpha)}\right)^{1/2}
	e^{-\alpha\eta}
	\left(\int_{\beta_3>0}\beta_3(I_b-T_b^4)^2\,\mathrm d\beta\right)^{1/2}.
	\end{align}
\end{lemma}
\begin{proof}
	The bounds in Lemma \ref{lem4.3} are uniform in $B$. A diagonal subsequence therefore satisfies
	\[
	\mathbf T^B\to\mathbf T\quad\text{locally in }C^1,
	\qquad
	\mathbf I^B\rightharpoonup\mathbf I
	\quad\text{locally in }L^2.
	\]
	Under the standing assumption $T_b\geq a_1>0$, the bound
\eqref{tpt} is uniform in $B$.
The local convergence of $\mathbf T^B$ therefore yields
\begin{align}\label{tph}
 \mathbf T(\eta)\geq\frac{a_1}{1+4\pi a_2^3}>0,
 \qquad \eta\geq0.
\end{align}

The strong convergence of $(\mathbf T^B)^4$, the weak convergence of
	$\mathbf I^B$, and the standard weak compactness of transport traces allow
	us to pass to the equations and to the lower boundary data. Thus the limit
	solves \eqref{3.16}--\eqref{3.17}. This is the compactness argument used in
	\cite[Sections 3.1--3.2]{GHM-2023-ARMA}; the scattering terms are linear and
	pass to the limit under the same weak convergence. Lower semicontinuity in
	\eqref{4.32} gives \eqref{4.54}, and \eqref{TIE} passes to the limit to give
	\eqref{4.55}.

	Fix $0<\alpha<1-\kappa$. Letting the upper endpoint in
	\eqref{4.54} tend to infinity and retaining the radiation dissipation, we obtain
	\begin{align}\label{hd}
	\|e^{\alpha z}(\mathbf I-\mathbf T^4)\|_{L^2(\mathbb R_+\times\mathbb S^2)}^2
	\leq\frac1{2(1-\kappa-\alpha)}
	\int_{\beta_3>0}\beta_3(I_b-T_b^4)^2\,\mathrm d\beta.
	\end{align}
	Since $\langle\beta_3\rangle=0$, identity \eqref{4.55} gives
	\[
	\partial_z\mathbf T(z)
	=\int_{\mathbb S^2}\beta_3
	\bigl(\mathbf I(z,\beta)-\mathbf T^4(z)\bigr)\,\mathrm d\beta.
	\]
	Applying Cauchy--Schwarz inequality in $\beta$ and using
	$\int_{\mathbb S^2}\beta_3^2\,\mathrm d\beta=4\pi/3$, we obtain
	\begin{align}\label{tde}
	\int_0^\infty e^{2\alpha z}|\partial_z\mathbf T|^2\,\mathrm dz
	&\leq\frac{4\pi}{3}
	\|e^{\alpha z}(\mathbf I-\mathbf T^4)\|_{L^2(\mathbb R_+\times\mathbb S^2)}^2\nonumber\\
	&\leq\frac{2\pi}{3(1-\kappa-\alpha)}
	\int_{\beta_3>0}\beta_3(I_b-T_b^4)^2\,\mathrm d\beta.
	\end{align}
	For $R>\eta\geq0$, the fundamental theorem of calculus
	and Cauchy--Schwarz inequality give
	\begin{align}\label{Reta}
	|\mathbf T(R)-\mathbf T(\eta)|
	\leq\frac{e^{-\alpha\eta}}{\sqrt{2\alpha}}
	\left(
	\int_\eta^R e^{2\alpha z}
	|\partial_z\mathbf T|^2\,\mathrm dz
	\right)^{1/2}.
	\end{align}
	By~\eqref{tde}, the right-hand side tends to zero as
	$\eta\to\infty$, uniformly in $R>\eta$.
	Thus $\mathbf T(\eta)$ converges to a limit
	$T_{0,\infty}$, which inherits the bounds
	\[
	\frac{a_1}{1+4\pi a_2^3}\leq T_{0,\infty}\leq a_2
	\]
	from~\eqref{tph} and $\mathbf T\leq a_2$.
	Letting $R\to\infty$ in \eqref{Reta}
	and using~\eqref{tde} proves~\eqref{4.56-1}.
	
	Applying Cauchy--Schwarz inequality in $\beta$ to the
	temperature equation in~\eqref{3.16}, and then using
	\eqref{hd}, we obtain
	\begin{align*}
	\int_0^\infty e^{2\alpha z}|\partial_z^2\mathbf T|^2\,\mathrm dz
	&\leq4\pi(1-\kappa)^2
	\|e^{\alpha z}(\mathbf I-\mathbf T^4)\|_{L^2(\mathbb R_+\times\mathbb S^2)}^2\\
	&\leq\frac{2\pi(1-\kappa)^2}{1-\kappa-\alpha}
	\int_{\beta_3>0}\beta_3(I_b-T_b^4)^2\,\mathrm d\beta.
	\end{align*}
	In particular, $\partial_z^2\mathbf T\in L^1(\mathbb R_+)$, since
	\[
	\int_0^\infty|\partial_z^2\mathbf T|\,\mathrm dz
	\leq\frac1{\sqrt{2\alpha}}
	\left(\int_0^\infty e^{2\alpha z}
	|\partial_z^2\mathbf T|^2\,\mathrm dz\right)^{1/2}<\infty.
	\]
	Hence $\partial_z\mathbf T$ has a finite limit at infinity.
	This limit is zero by \eqref{tde}.
	Integrating the temperature equation from $\eta$ to infinity
	and applying Cauchy--Schwarz inequality, we therefore obtain
	\begin{align*}
	|\partial_\eta\mathbf T(\eta)|
	&=\left|\int_\eta^\infty\partial_z^2\mathbf T(z)\,\mathrm dz\right|\leq\frac{e^{-\alpha\eta}}{\sqrt{2\alpha}}
	\left(\int_\eta^\infty e^{2\alpha z}
	|\partial_z^2\mathbf T|^2\,\mathrm dz\right)^{1/2}\\
	&\leq(1-\kappa)
	\left(\frac{\pi}{\alpha(1-\kappa-\alpha)}\right)^{1/2}
	e^{-\alpha\eta}
	\left(\int_{\beta_3>0}\beta_3(I_b-T_b^4)^2\,\mathrm d\beta\right)^{1/2},
	\end{align*}
	which proves \eqref{4.60}. Therefore we have finished the proof of Lemma \ref{lem4.h}.
\end{proof}

\subsection{Proof of Theorem~\ref{thm2}}\label{ntreg}
\begin{proof}
	Lemma~\ref{lem4.h} gives existence, the bounds~\eqref{4.00},
	and the temperature decay estimate in~\eqref{4.0}.
	We first prove the radiation decay estimate in~\eqref{4.0}
	and then the stability estimates~\eqref{hshift}--\eqref{hfarstab}.
	Set
	\[
	J(\eta,\beta)=\mathbf I(\eta,\beta)-\mathbf T^4(\eta).
	\]
	The first equation in \eqref{3.16} becomes
	\begin{align}\label{4.61}
	\beta_3\partial_\eta J+J-\kappa\bar J
	=-\beta_3\partial_\eta(\mathbf T^4),
	\end{align}
	with $J(0,\beta)=I_b(\beta)-T_b^4$ for $\beta_3>0$.

	Fix $0<\alpha<1-\kappa$ and define
	\[
	M_\alpha=\|e^{\alpha\eta}J\|_{L^\infty(\mathbb R_+\times\mathbb S^2)}.
	\]
	For $\beta_3>0$, the characteristic formula for \eqref{4.61} is
	\begin{align}\label{4.58-1}
		J(\eta,\beta)=(I_b-T_b^4)e^{-\f{\eta}{\beta_3}}+\int_0^\eta \left(\kappa\bar{J}-\beta_3\partial_z(\mathbf{T}^4)\right)\f{1}{\beta_3}e^{-\f{\eta-z}{\beta_3}}\,\mathrm{dz}.
	\end{align}
	For $\beta_3<0$, one has 
	\begin{align}\label{4.58}
		J(\eta,\beta)=J(\eta_1,\beta)e^{-\f{\eta_1-\eta}{|\beta_3|}}+\int_{\eta}^{\eta_1}\frac{1}{|\beta_3|}e^{-\f{z-\eta}{|\beta_3|}}\left( \kappa\bar{J}-\beta_3\partial_z(\mathbf{T}^4)\right)\mathrm{d}z ,
	\end{align}
	where $\eta<\eta_1$. Letting $\eta_1\to\infty$ removes the terminal term because
	$J$ is bounded. Since
	\[
	\int_0^\eta\frac1{\beta_3}e^{-(\eta-z)/\beta_3}e^{-\alpha z}\,\mathrm dz
	\leq\frac{e^{-\alpha\eta}}{1-\alpha},
	\qquad \beta_3>0,
	\]
	and
	\[
	\int_\eta^\infty\frac1{|\beta_3|}e^{-(z-\eta)/|\beta_3|}e^{-\alpha z}\,\mathrm dz
	\leq e^{-\alpha\eta},
	\qquad \beta_3<0,
	\]
	we get 
	\[
	M_\alpha
	\leq C|I_b-T_b^4|_{L^\infty(\{\beta_3>0\})}
	+\frac{\kappa}{1-\alpha}M_\alpha
	+\frac{C}{1-\alpha}
	\|e^{\alpha\eta}\partial_\eta(\mathbf T^4)\|_{L^\infty(\mathbb R_+)}.
	\]
	Since $0\leq\mathbf T\leq a_2$, estimate \eqref{4.60} and the finite measure of the incoming hemisphere imply
	\[
	\|e^{\alpha\eta}\partial_\eta(\mathbf T^4)\|_{L^\infty}
	\leq\frac{C(1-\kappa)}{\alpha^{1/2}(1-\kappa-\alpha)^{1/2}}
	|I_b-T_b^4|_{L^\infty(\{\beta_3>0\})}.
	\]
	
	Noting
	\[
	1-\frac{\kappa}{1-\alpha}=\frac{1-\kappa-\alpha}{1-\alpha}>0,
	\]
	one has
	\[
	M_\alpha\leq C\left(
	\frac{1-\alpha}{1-\kappa-\alpha}
	+\frac{1-\kappa}{\alpha^{1/2}(1-\kappa-\alpha)^{3/2}}
	\right)|I_b-T_b^4|_{L^\infty(\{\beta_3>0\})}.
	\]
	Substituting it into \eqref{4.58-1} and \eqref{4.58}, we obtain the same exponential decay pointwise
	for the characteristic representative of $J$.
	Finally,
	\[
	\mathbf I-T_{0,\infty}^4
	=J+\mathbf T^4-T_{0,\infty}^4.
	\]
	Using \eqref{4.56-1}, the bound $0\leq\mathbf T,T_{0,\infty}\leq a_2$,	we obtain
	\[
	e^{\alpha\eta}|\mathbf I(\eta,\beta)-T_{0,\infty}^4|
	\leq C\left(
	\frac{\kappa}{1-\kappa-\alpha}
	+\frac{1-\kappa}{\alpha^{1/2}(1-\kappa-\alpha)^{3/2}}
	\right)|I_b-T_b^4|_{L^\infty(\{\beta_3>0\})}.
	\]
	Together with~\eqref{4.56-1}, this proves~\eqref{4.0}.

	We now prove the stability estimates with respect to the boundary data.
	For $i=1,2$, let $(\mathbf I_i^B,\mathbf T_i^B)$ be the
	truncated solution corresponding to $(I_{b,i},T_{b,i})$.
	The compactness argument in Lemma~\ref{lem4.h} applies
	simultaneously to these two families. Thus a common subsequence
	$B_n\to\infty$ can be chosen such that, for every $R>0$ and
	both $i=1,2$,
	\[
	\mathbf T_i^{B_n}\longrightarrow\mathbf T_i
	\quad\text{in }C^1([0,R]),
	\qquad
	\mathbf I_i^{B_n}\rightharpoonup\mathbf I_i
	\quad\text{in }L^2((0,R)\times\mathbb S^2).
	\]
	Fix $R>0$. For all sufficiently large $n$, we have $B_n>R$.
	Restricting the two norms on the left-hand side of~\eqref{tshift}
	to $(0,R)\times\mathbb S^2$ and $(0,R)$, respectively, gives
	\begin{align*}
		&\|\mathbf I_1^{B_n}-\mathbf I_2^{B_n}\|_
		{L^\infty((0,R)\times\mathbb S^2)}
		+\|\mathbf T_1^{B_n}-\mathbf T_2^{B_n}\|_
		{L^\infty(0,R)}
		\\
		&\qquad\leq C\Bigl(
		\|I_{b,1}-I_{b,2}\|_{L^\infty(\{\beta_3>0\})}
		+|T_{b,1}-T_{b,2}|
		\Bigr).
	\end{align*}
	The constant is independent of $n$, $R$, and $\kappa$.
	
	Since $0\leq\mathbf I_i^{B_n}\leq a_2^4$, the local
	$L^2$ weak convergence also implies
	\[
	\mathbf I_i^{B_n}\overset{*}{\rightharpoonup}\mathbf I_i
	\quad\text{in }L^\infty((0,R)\times\mathbb S^2).
	\]
	Indeed, the convergence against $L^2$ test functions extends
	to all $L^1$ test functions by density and the uniform
	$L^\infty$ bounds. Consequently, weak-star lower
	semicontinuity gives
	\begin{align}\label{Il}
		\|\mathbf I_1-\mathbf I_2\|_
		{L^\infty((0,R)\times\mathbb S^2)}\leq
		\liminf_{n\to\infty}
		\|\mathbf I_1^{B_n}-\mathbf I_2^{B_n}\|_
		{L^\infty((0,R)\times\mathbb S^2)}.
	\end{align}
	On the other hand, uniform convergence of the temperatures
	on $[0,R]$ yields
	\begin{align}\label{Tl}
	\|\mathbf T_1-\mathbf T_2\|_{L^\infty(0,R)}
	=
	\lim_{n\to\infty}
	\|\mathbf T_1^{B_n}-\mathbf T_2^{B_n}\|_{L^\infty(0,R)}.
	\end{align}
	Combining \eqref{Il}, \eqref{Tl} and passing to the lower limit
	in the restricted estimate above, we obtain
	\begin{align*}
		&\|\mathbf I_1-\mathbf I_2\|_
		{L^\infty((0,R)\times\mathbb S^2)}
		+\|\mathbf T_1-\mathbf T_2\|_{L^\infty(0,R)}
		\\
		&\qquad\leq C\Bigl(
		\|I_{b,1}-I_{b,2}\|_{L^\infty(\{\beta_3>0\})}
		+|T_{b,1}-T_{b,2}|
		\Bigr).
	\end{align*}
	Letting $R\to\infty$ proves~\eqref{hshift}.
	
	Finally, applying~\eqref{4.56-1} to each half-space solution
	gives $\mathbf T_i(\eta)\to T_{0,\infty,i}$ as
	$\eta\to\infty$. Hence
	\[
	|T_{0,\infty,1}-T_{0,\infty,2}|
	=
	\lim_{\eta\to\infty}
	|\mathbf T_1(\eta)-\mathbf T_2(\eta)|
	\leq
	\|\mathbf T_1-\mathbf T_2\|_{L^\infty(\mathbb R_+)}.
	\]
	Using~\eqref{hshift} on the right-hand side proves
	\eqref{hfarstab}. This completes the proof of Theorem~\ref{thm2}.
\end{proof}

\section{Linear problem}\label{lin}
To study \eqref{5.1}--\eqref{5.1-2}, we first consider the case in which the source in the temperature equation is the angular integral of the source in the transport equation. In detail, we consider
\begin{align}\label{5.1-3}
	\begin{cases}
		\beta_3\partial_\eta\phi+\kappa(\phi-\bar{\phi})=(1-\kappa)(4\mathbf{T}^3g-\phi)+S,\\
		\partial_\eta^2g+(1-\kappa)4\pi(\bar{\phi}-4\mathbf{T}^3g)=\langle S\rangle,
	\end{cases}
\end{align}
with boundary conditions 
\begin{align}
	\phi(0,\beta)|_{\gamma_-}=p_b(\beta),\quad g(0)=0,\label{5.1-4}
\end{align}
and $p_b(\beta)$, $S(\eta,\beta)$ are given functions.

\begin{theorem}\label{thm4}
	Let $(\mathbf I,\mathbf T)$ be the nonlinear Milne profile from Theorem~\ref{thm2}. Suppose $\mathbf{T}$ satisfies the spectral assumption (A) for some constant $0<\vartheta<1-\kappa$ and $0<M<1$. Assume $p_b$ is bounded and 
	\[
	e^{\vartheta\eta}S\in L^\infty(\mathbb R_+\times\mathbb S^2)
	\cap L^2(\mathbb R_+\times\mathbb S^2).
	\]
	Then, there exists a unique bounded solution $\left(\phi,g\right) \in L^2_{\rm loc}(\mathbb{R}_{+}\times\mathbb S^2)\times L^2_{\rm loc}(\mathbb{R}_+) $
	to system \eqref{5.1-3}-\eqref{5.1-4} satisfying
	\begin{align}
		|g(\eta)|\leq N_\vartheta,\quad
		\|\phi\|_{L^\infty}\leq 4a_2^3N_\vartheta
		+|p_b|_{L^\infty(\gamma_-)}
		+\frac{\|S\|_{L^\infty(\mathbb R_+\times\mathbb S^2)}}{1-\kappa}.
	\end{align} 
	Moreover, there exists a constant $g_\infty \in \mathbb{R}$ satisfying $|g_\infty|\leq N_\vartheta$
	such that 
	\begin{equation}\label{5.00}
		|g(\eta)-g_\infty| \le N_\vartheta e^{-\vartheta\eta},\quad
			|\phi- 4T_{0,\infty}^3 g_\infty |\leq M_\vartheta e^{-\vartheta\eta},
	\end{equation}
	where $N_\vartheta$ and $M_\vartheta$ are defined in \eqref{5.00-1} and \eqref{5.32}, respectively. 
\end{theorem}

In this section, we assume that $\mathbf{T}$ satisfies the spectral assumption. Before proof, we derive some inequalities from the spectral assumption (A).
\begin{lemma}[\cite{GHM-2023-ARMA}]
	Assume $T\in C_{\text{loc}}^1(\R_+)$ satisfies the spectral assumption (A). Then for any function $f\in C_{\text{loc}}^1(\R_+)$ with $f(0)=0$, the following inequalities hold:
	\begin{align}
		&M\int_0^\infty(2T^{\f32})^2|\partial_\eta f|^2\,\mathrm{d}\eta\geq 4\int_0^\infty|\partial_\eta(2T^{\f32})|^2f^2\,\mathrm{d}\eta,\label{5.0-1}\\
		&\int_0^\infty 4T^3|\partial_\eta f|^2\,\mathrm{d}\eta+\int_0^\infty\partial_\eta(4T^3)f\partial_\eta f\,\mathrm{d}\eta\geq 0,\label{5.0-2}
	\end{align}
	and 
	\begin{align}\label{5.0-3}
		\int_0^\infty e^{2\vartheta \eta}4T^3|\partial_\eta f|^2\,\mathrm{d}\eta+\int_0^\infty e^{2\vartheta \eta}\partial_\eta(4T^3)f\partial_\eta f\,\mathrm{d}\eta\geq 0.
	\end{align}
\end{lemma}

\subsection{Existence in bounded domain}
We first construct solutions for the following truncated problem:
\begin{align}\label{5.2}
	\begin{cases}
		\beta_3\partial_\eta\phi^B+\kappa(\phi^B-\bar{\phi}^B)=(1-\kappa)(4\mathbf{T}^3g^B-\phi^B)+S,\\
		\partial_\eta^2g^B+(1-\kappa)4\pi(\bar{\phi}^B-4\mathbf{T}^3g^B)=\langle S\rangle,
	\end{cases}\quad (\eta,\beta)\in(0,B)\times\mathbb{S}^2
\end{align}
with boundary conditions 
\begin{align}
	&\phi^B(0,\beta)|_{\gamma_-}=p_b(\beta),\quad \phi^B(B,\beta)|_{\gamma_-}=\phi^B(B,\beta_{\shortparallel},-\beta_3),\label{5.3}\\
	&g^B(0)=0,\quad \partial_\eta g^B(B)=0.\label{5.4}
\end{align}
\begin{lemma}\label{lem5.1}
	Let $S:[0,B]\times\mathbb S^2\to\R$ be a continuous function and $p_b$ continuous on the boundary. Then there exists a unique solution $(\phi^B,g^B)\in C([0,B]\times\mathbb{S}^2\backslash\gamma_0)\times C^1([0,B])$ to system \eqref{5.2}--\eqref{5.4}, and the solution satisfies 
	\begin{align}\label{5.25}
			&\f12|\phi^B(0)|_{L^2(\gamma_+)}^2+\kappa\|e^{\vartheta \eta}(\phi^B-\bar{\phi}^B)\|_{L^2}^2+(1-\kappa-\vartheta-\xi)\|e^{\vartheta \eta}(\phi^B-4\mathbf{T}^3g^B)\|_{L^2}^2\nonumber\\
		&\leq \f{1}{2\xi}\|e^{\vartheta\eta} S\|_{L^2}^2+\f12|p_b|_{L^2(\gamma_-)}^2,
	\end{align}
	where $\vartheta< 1-\kappa$ and $\xi<1-\kappa-\vartheta$ are small positive constants.
\end{lemma}
The proof is similar to that of \cite[Lemma 3]{GHM-2023-ARMA},
with a few modifications due to the scattering term.
We give the details in Appendix~\ref{app:lem5.1}.

For $0<\vartheta<1-\kappa$, define
\begin{align}\label{5.00-1}
	N_\vartheta
	:=\frac1{\sqrt{2\vartheta}}
	\left(
	\frac{4\pi}{3(1-\kappa-\vartheta)}|p_b|_{L^2(\gamma_-)}^2
	+\frac{8\pi}{3(1-\kappa-\vartheta)^2}
	\|e^{\vartheta\eta}S\|_{L^2(\mathbb R_+\times\mathbb S^2)}^2
	\right)^{1/2},
\end{align}
and
\begin{align}\label{5.32}
	M_\vartheta:=\frac1{1-\kappa-\vartheta}\Biggl[
	&(1-\vartheta)\bigl(|p_b|_{L^\infty(\gamma_-)}
	+4a_2^3N_\vartheta\bigr)
	+\|e^{\vartheta\eta}S\|_{L^\infty(\mathbb R_+\times\mathbb S^2)}
	\nonumber\\
	&+4(1-\kappa)\left(
	a_2^3+\frac{3a_2^2}{\sqrt{2\vartheta}}
	\left(\frac{2\pi}{3(1-\kappa-\vartheta)}
	\int_{\beta_3>0}\beta_3(I_b-T_b^4)^2\,\mathrm d\beta
	\right)^{1/2}\right)N_\vartheta
	\Biggr].
\end{align}

\begin{lemma}\label{lem4.4}
	Let $(\phi^B,g^B)$ be the solution to \eqref{5.2}--\eqref{5.4} constructed in Lemma~\ref{lem5.1}. Then, for $0\leq\eta\leq B$,
	\begin{align*}
		|g^B(\eta)|\leq N_\vartheta, \quad |\partial_\eta g^B(\eta)|\leq 3(1-\kappa)N_\vartheta e^{-\vartheta\eta}, 
	\end{align*}
	and 
	\begin{align}
		|g^B(\eta)-g^B(B)|\leq N_\vartheta e^{-\vartheta\eta},
	\end{align} 
	where $N_\vartheta$ is defined in \eqref{5.00-1}. Moreover, 
		\begin{align}\label{S0}
		\|\phi^B\|_{L^\infty(\Omega_B)}\leq 4a_2^3N_\vartheta
		+|p_b|_{L^\infty(\{\beta_3>0\})}
		+\frac{\|S\|_{L^\infty(\Omega_B)}}{1-\kappa}.
	\end{align}
	The quantity $N_\vartheta$ uses the half-space source norm in~\eqref{5.00-1} and is independent of $B$.
\end{lemma}
\begin{proof}
By~\eqref{5.21-1} and $\langle\beta_3\rangle=0$, we have
$\partial_\eta g^B
=\langle\beta_3(\phi^B-4\mathbf T^3g^B)\rangle$.
Taking $\xi=(1-\kappa-\vartheta)/2$ in~\eqref{5.25}
and applying Cauchy--Schwarz inequality in $\beta$, we obtain
\begin{align}\label{5.29-1}
	\int_0^B e^{2\vartheta\eta}
	|\partial_\eta g^B|^2\,\mathrm d\eta
	&\leq\frac{4\pi}{3}
	\|e^{\vartheta\eta}(\phi^B-4\mathbf T^3g^B)\|
	_{L^2(\Omega_B)}^2\leq2\vartheta N_\vartheta^2.
\end{align}
The last inequality uses~\eqref{5.00-1} and the fact that
the source norm on $\Omega_B$ is bounded by its half-space norm. For $0\leq\eta_1\leq\eta_2\leq B$,
Cauchy--Schwarz inequality gives
\begin{align*}
	|g^B(\eta_2)-g^B(\eta_1)|
	\leq
	\left(\int_{\eta_1}^{\eta_2}
	e^{-2\vartheta z}\,\mathrm dz\right)^{1/2}
	\left(\int_{\eta_1}^{\eta_2}
	e^{2\vartheta z}|\partial_zg^B|^2\,\mathrm dz\right)^{1/2}\leq N_\vartheta e^{-\vartheta\eta_1}.
\end{align*}
Taking $(\eta_1,\eta_2)=(0,\eta)$ and $(\eta,B)$,
and using $g^B(0)=0$, proves the bounds for
$g^B(\eta)$ and $g^B(\eta)-g^B(B)$.

	To estimate $\partial_\eta g^B$ pointwise, we use
	$\partial_\eta g^B(B)=0$ and \eqref{5.2}$_2$ to obtain
	\begin{align*}
	\partial_\eta g^B(\eta)
	&=\partial_\eta g^B(B)-\int_\eta^B\partial_z^2g^B(z)\,\mathrm dz\\
	&=(1-\kappa)\int_\eta^B
	\langle\phi^B-4\mathbf T^3g^B\rangle(z)\,\mathrm dz
	-\int_\eta^B\langle S(z,\cdot)\rangle\,\mathrm dz.
	\end{align*}
	Applying the Cauchy--Schwarz inequality first in $\beta$ and then in $z$, we obtain
	\begin{align}\label{peg}
	|\partial_\eta g^B(\eta)|
	&\leq\sqrt{4\pi}\int_\eta^B
	\left[(1-\kappa)\|\phi^B(z,\cdot)-4\mathbf T^3(z)g^B(z)\|_{L^2_\beta}
	+\|S(z,\cdot)\|_{L^2_\beta}\right]\,\mathrm dz\nonumber\\
	&\leq\sqrt{4\pi}
	\left(\int_\eta^B e^{-2\vartheta z}\,\mathrm dz\right)^{1/2}
	\left[(1-\kappa)\|e^{\vartheta z}(\phi^B-4\mathbf T^3g^B)\|_{L^2}
	+\|e^{\vartheta z}S\|_{L^2}\right]\nonumber\\
	&\leq\sqrt{\frac{2\pi}{\vartheta}}e^{-\vartheta\eta}
	\left[(1-\kappa)\|e^{\vartheta z}(\phi^B-4\mathbf T^3g^B)\|_{L^2}
	+\|e^{\vartheta z}S\|_{L^2}\right],
	\end{align}
	where the $L^2$ norms in this calculation are taken over
	$[0,B]\times\mathbb S^2$. Taking
	$\xi=(1-\kappa-\vartheta)/2$ in \eqref{5.25} and using the
	definition of $N_\vartheta$, we have
	\begin{align*}
	\|e^{\vartheta z}(\phi^B-4\mathbf T^3g^B)\|_{L^2}^2
	&\leq\frac{|p_b|_{L^2(\gamma_-)}^2}{1-\kappa-\vartheta}
	+\frac{2\|e^{\vartheta z}S\|_{L^2}^2}{(1-\kappa-\vartheta)^2}
	\leq\frac{3\vartheta}{2\pi}N_\vartheta^2,\\
	\|e^{\vartheta z}S\|_{L^2}^2
	&\leq\frac{3\vartheta(1-\kappa-\vartheta)^2}{4\pi}N_\vartheta^2.
	\end{align*}
	Substituting these bounds into \eqref{peg}, we obtain
	\begin{align*}
	|\partial_\eta g^B(\eta)|
	\leq\left[\sqrt3(1-\kappa)
	+\sqrt{\frac32}(1-\kappa-\vartheta)\right]
	N_\vartheta e^{-\vartheta\eta}\leq3(1-\kappa)N_\vartheta e^{-\vartheta\eta},
	\end{align*}
	since $0<\vartheta<1-\kappa$ and $\sqrt3+\sqrt{3/2}<3$.
	Using the characteristic formulas~\eqref{5.5}--\eqref{5.5-1},
	$|g^B|\leq N_\vartheta$, and $0\leq\mathbf T\leq a_2$, we obtain
	\[
	\|\phi^B\|_{L^\infty(\Omega_B)}
	\leq 4a_2^3N_\vartheta
	+|p_b|_{L^\infty(\{\beta_3>0\})}
	+\frac{\|S\|_{L^\infty(\Omega_B)}}{1-\kappa},
	\]
	which proves~\eqref{S0}.
\end{proof}
\subsection{Existence on half space}\label{linhalf}
Using Lemmas~\ref{lem5.1} and~\ref{lem4.4}, the compactness and
trace argument in~\cite[Section 4.3]{GHM-2023-ARMA} gives a bounded
solution of~\eqref{5.1-3}--\eqref{5.1-4} as $B\to\infty$.
The additional term $\kappa(\phi^B-\bar\phi^B)$ passes to the limit
under local $L^2$ weak convergence.
Passing to the limit in~\eqref{5.21-1} and using lower
semicontinuity in~\eqref{5.29-1}, we obtain
\begin{align}
 \partial_\eta g=\langle\beta_3\phi\rangle,
\quad \text{and}\,
 \|e^{\vartheta\eta}\partial_\eta g\|_{L^2(\mathbb R_+)}^2
 \leq2\vartheta N_\vartheta^2.
 \label{5.31-1}
\end{align}
Integrating $\partial_\eta g$ from $0$ or toward infinity and
applying the Cauchy--Schwarz inequality, as
in~\cite[Section 4.4]{GHM-2023-ARMA}, gives
$g_\infty:=\lim_{\eta\to\infty}g(\eta)$,
$|g_\infty|\leq N_\vartheta$, and the temperature estimate
in~\eqref{5.00}.

We now give the radiation decay argument, retaining the contribution
of the scattering term. For $\beta_3>0$, integrating
\eqref{5.1-3}$_1$ along the characteristic from $\eta=0$ gives
\begin{align}\label{5.26-1}
 \phi(\eta,\beta)-4T_{0,\infty}^3g_\infty=&(p_b(\beta)-4T_{0,\infty}^3g_\infty)e^{-\eta/\beta_3}
 \nonumber\\
 &+\int_0^\eta\frac1{\beta_3}e^{-(\eta-z)/\beta_3}
 \Bigl[\kappa(\bar\phi(z)-4T_{0,\infty}^3g_\infty)\nonumber\\
 &\qquad\qquad
 +(1-\kappa)(4\mathbf T^3(z)g(z)-4T_{0,\infty}^3g_\infty)
 +S(z,\beta)\Bigr]\,\mathrm dz.
\end{align}
For $\beta_3<0$, the boundedness of $\phi$ allows the terminal
point on the backward characteristic to tend to infinity,
and we obtain
\begin{align}\label{5.26-2}
 &\phi(\eta,\beta)-4T_{0,\infty}^3g_\infty
 \nonumber\\
 &=\int_\eta^\infty\frac1{|\beta_3|}e^{-(z-\eta)/|\beta_3|}
 \Bigl[\kappa(\bar\phi(z)-4T_{0,\infty}^3g_\infty)
 +(1-\kappa)(4\mathbf T^3(z)g(z)-4T_{0,\infty}^3g_\infty)
 +S(z,\beta)\Bigr]\,\mathrm dz.
\end{align}

To estimate the temperature contribution in these formulas, we write
\[
 4\mathbf T^3g-4T_{0,\infty}^3g_\infty
 =4\mathbf T^3(g-g_\infty)
 +4(\mathbf T^3-T_{0,\infty}^3)g_\infty.
\]
Using $|g_\infty|\leq N_\vartheta$, the temperature estimate
in~\eqref{5.00}, and~\eqref{4.56-1} with $\alpha=\vartheta$,
we obtain
\begin{align}\label{lforcing}
 &\|e^{\vartheta\eta}
       (4\mathbf T^3g-4T_{0,\infty}^3g_\infty)\|_{L^\infty(\mathbb R_+)}
 \nonumber\\
 &\quad\leq4\left[
 a_2^3+\frac{3a_2^2}{\sqrt{2\vartheta}}
 \left(\frac{2\pi}{3(1-\kappa-\vartheta)}
 \int_{\beta_3>0}\beta_3(I_b-T_b^4)^2\,\mathrm d\beta
 \right)^{1/2}\right]N_\vartheta.
\end{align}
Here we used
$|\mathbf T^3-T_{0,\infty}^3|
\leq3a_2^2|\mathbf T-T_{0,\infty}|$.

With $g$ fixed, iteration of the linear characteristic
formulas~\eqref{5.26-1}--\eqref{5.26-2} gives the weighted
$L^\infty$ estimate below.
Multiplying \eqref{5.26-1} and \eqref{5.26-2} by $e^{\vartheta\eta}$, taking the
$L^\infty$ norm, and using~\eqref{lforcing}, we obtain
\begin{align}\label{lrad}
 &\|e^{\vartheta\eta}
       (\phi-4T_{0,\infty}^3g_\infty)\|_{L^\infty}
 \nonumber\\
 &\quad\leq |p_b|_{L^\infty(\gamma_-)}+4a_2^3N_\vartheta
 +\frac{\|e^{\vartheta\eta}S\|_{L^\infty}}{1-\vartheta}+\frac{\kappa}{1-\vartheta}
 \|e^{\vartheta\eta}
 (\phi-4T_{0,\infty}^3g_\infty)\|_{L^\infty}
 \nonumber\\
 &\qquad+\frac{4(1-\kappa)}{1-\vartheta}
 \left[
 a_2^3+\frac{3a_2^2}{\sqrt{2\vartheta}}
 \left(\frac{2\pi}{3(1-\kappa-\vartheta)}
 \int_{\beta_3>0}\beta_3(I_b-T_b^4)^2\,\mathrm d\beta
 \right)^{1/2}\right]N_\vartheta.
\end{align}
Since
\[
 1-\frac{\kappa}{1-\vartheta}
 =\frac{1-\kappa-\vartheta}{1-\vartheta}>0,
\]
absorbing the third term on the right hand side of \eqref{lrad} and using the
definition~\eqref{5.32}, we obtain
\[
 \|e^{\vartheta\eta}(\phi-4T_{0,\infty}^3g_\infty)\|_{L^\infty}
 \leq M_\vartheta,
\]
which proves the radiation estimate in~\eqref{5.00}.

\subsection{Uniqueness}

Let $(g_1,\phi_1)$ and $(g_2,\phi_2)$ be two bounded solutions
of~\eqref{5.1-3}--\eqref{5.1-4} with the same data.
Subtracting their equations and boundary conditions, we find that
$h:=g_1-g_2$ and $\varphi:=\phi_1-\phi_2$ satisfy
\begin{align}\label{5.33}
 \begin{cases}
 \beta_3\partial_\eta\varphi
 =-\kappa(\varphi-\bar\varphi)
  -(1-\kappa)(\varphi-4\mathbf T^3h),\\
 \partial_\eta^2h
 =-(1-\kappa)\langle\varphi-4\mathbf T^3h\rangle,
 \end{cases}
\end{align}
with
\begin{align}\label{5.34}
 \varphi(0,\beta)=0\quad\text{for }\beta_3>0,
 \qquad h(0)=0.
\end{align}

Integrating~\eqref{5.33}$_1$ over $\mathbb S^2$ and subtracting
\eqref{5.33}$_2$ gives the first identity below.
Multiplying~\eqref{5.33}$_1$ by $\beta_3$ and integrating over
$\mathbb S^2$, with $\langle\beta_3\rangle=0$, gives the second:
\begin{align}\label{5.35}
 \partial_\eta\bigl(\partial_\eta h-\langle\beta_3\varphi\rangle\bigr)
 &=0,
 &\partial_\eta\langle\beta_3^2\varphi\rangle
 &=-\langle\beta_3\varphi\rangle.
\end{align}
Differentiating the second identity and using the first, we obtain
$\partial_\eta^2(h+\langle\beta_3^2\varphi\rangle)=0$.
Since $h$ and $\varphi$ are bounded, this affine function is constant.
The second identity in~\eqref{5.35} therefore yields $\partial_\eta h=\langle\beta_3\varphi\rangle.$
The moment identities are understood in distributions and hold
for their locally absolutely continuous representatives.

Fix $R>0$. Multiplying~\eqref{5.33}$_1$ by $\varphi$ and integrating
over $(0,R)\times\mathbb S^2$, multiplying~\eqref{5.33}$_2$
by $-4\mathbf T^3h$ and integrating over $(0,R)$, and adding
the resulting identities, we obtain
\begin{align}\label{5.37}
 &-\int_0^R4\mathbf T^3h\partial_\eta^2h\,\mathrm d\eta
 +\frac12\int_0^R\partial_\eta
       \langle\beta_3\varphi^2\rangle\,\mathrm d\eta
 \nonumber\\
 &\quad +(1-\kappa)\int_0^R\int_{\mathbb S^2}
       |\varphi-4\mathbf T^3h|^2\,\mathrm d\beta\,\mathrm d\eta
 +\kappa\int_0^R\int_{\mathbb S^2}
       |\varphi-\bar\varphi|^2\,\mathrm d\beta\,\mathrm d\eta=0.
\end{align}
Here the scattering term was simplified using
$\langle(\varphi-\bar\varphi)\varphi\rangle
=\langle|\varphi-\bar\varphi|^2\rangle$.
To combine the boundary terms, we use
$\langle\beta_3\rangle=0$ and
$\langle\beta_3\varphi\rangle=\partial_\eta h$ to write
\[
 -4\mathbf T^3h\partial_\eta h
 +\frac12\langle\beta_3\varphi^2\rangle
 =\frac12\langle\beta_3(\varphi-4\mathbf T^3h)^2\rangle.
\]
Integrating the first two terms in~\eqref{5.37} by parts,
substituting this boundary identity, and using~\eqref{5.34}
at $\eta=0$, we obtain
\begin{align}\label{5.36}
 &-\int_0^R4\mathbf T^3h\partial_\eta^2h\,\mathrm d\eta
 +\frac12\int_0^R\partial_\eta
       \langle\beta_3\varphi^2\rangle\,\mathrm d\eta
 \nonumber\\
 &=\int_0^R4\mathbf T^3|\partial_\eta h|^2\,\mathrm d\eta
 +\int_0^R\partial_\eta(4\mathbf T^3)h\partial_\eta h\,\mathrm d\eta
 \nonumber\\
 &\quad+\frac12\langle\beta_3(\varphi-4\mathbf T^3h)^2\rangle(R)
 +\frac12\int_{\beta_3<0}|\beta_3|\varphi^2(0,\beta)\,\mathrm d\beta.
\end{align}
Extending $h|_{[0,R]}$ by the constant $h(R)$ for $\eta>R$
and applying~\eqref{5.0-2}, we obtain
\[
 \int_0^R4\mathbf T^3|\partial_\eta h|^2\,\mathrm d\eta
 +\int_0^R\partial_\eta(4\mathbf T^3)h\partial_\eta h\,\mathrm d\eta
 \geq0.
\]
Substituting~\eqref{5.36} into~\eqref{5.37} and discarding this
nonnegative sum gives
\begin{align}\label{lue}
 &\frac12\int_{\beta_3<0}|\beta_3|\varphi^2(0,\beta)\,\mathrm d\beta
 +(1-\kappa)\int_0^R\int_{\mathbb S^2}
       |\varphi-4\mathbf T^3h|^2\,\mathrm d\beta\,\mathrm d\eta
 \nonumber\\
 &\quad+\kappa\int_0^R\int_{\mathbb S^2}
       |\varphi-\bar\varphi|^2\,\mathrm d\beta\,\mathrm d\eta
 \leq-\frac12\langle\beta_3(\varphi-4\mathbf T^3h)^2\rangle(R).
\end{align}
Since $0\leq\mathbf T\leq a_2$, the boundary term at $R$ satisfies
\[
 \frac12\left|\langle\beta_3(\varphi-4\mathbf T^3h)^2\rangle(R)\right|
 \leq2\pi\left(
 \|\varphi\|_{L^\infty(\mathbb R_+\times\mathbb S^2)}
 +4a_2^3\|h\|_{L^\infty(\mathbb R_+)}
 \right)^2.
\]
This bound is independent of $R$. Using it in~\eqref{lue}
and applying the monotone convergence theorem to the integral
with coefficient $1-\kappa>0$, we deduce that
\[
 \int_0^\infty\int_{\mathbb S^2}
 |\varphi-4\mathbf T^3h|^2\,\mathrm d\beta\,\mathrm d\eta<\infty.
\]
Hence there exists a sequence $R_j\uparrow\infty$ such that
\[
 \left|\langle\beta_3(\varphi-4\mathbf T^3h)^2\rangle(R_j)\right|
 \leq\int_{\mathbb S^2}
 |\varphi(R_j,\beta)-4\mathbf T^3(R_j)h(R_j)|^2\,\mathrm d\beta
 \longrightarrow0.
\]
Taking $R=R_j$ in~\eqref{lue} and applying the
monotone convergence theorem to the nonnegative integrals, we obtain
\begin{align*}
 &\frac12\int_{\beta_3<0}|\beta_3|\varphi^2(0,\beta)\,\mathrm d\beta
 +(1-\kappa)\int_0^\infty\int_{\mathbb S^2}
       |\varphi-4\mathbf T^3h|^2\,\mathrm d\beta\,\mathrm d\eta\\
 &\quad+\kappa\int_0^\infty\int_{\mathbb S^2}
       |\varphi-\bar\varphi|^2\,\mathrm d\beta\,\mathrm d\eta
 \leq0.
\end{align*}
Since $1-\kappa>0$, it follows that $\varphi=4\mathbf T^3h$
almost everywhere. Substituting this identity into~\eqref{5.33}$_2$
gives $\partial_\eta^2h=0$.
The boundedness of $h$ and the boundary condition $h(0)=0$
then imply $h=0$, and consequently $\varphi=0$.
This proves uniqueness among bounded solutions and completes
the proof of Theorem~\ref{thm4}.$\hfill\Box$

\subsection{Proof of Theorem \ref{co}.}
We first estimate the solution of~\eqref{5.1-3}--\eqref{5.1-4}.
The orthogonal decomposition gives
\[
\|e^{\vartheta\eta}(\phi^B-4\mathbf T^3g^B)\|_2^2
=\|e^{\vartheta\eta}(\phi^B-\bar\phi^B)\|_2^2
+4\pi\|e^{\vartheta\eta}(\bar\phi^B-4\mathbf T^3g^B)\|_2^2.
\]
Substituting~\eqref{5.24-0}, \eqref{5.23}, and~\eqref{5.24}
into~\eqref{5.22}, and using $|\beta_3|\leq1$ and the Young inequality,
we obtain
\begin{align*}
&\frac12|\phi^B(0)|_{L^2(\gamma_+)}^2
+(1-\vartheta)\|e^{\vartheta\eta}(\phi^B-\bar\phi^B)\|_2^2\\
&\quad+4\pi(1-\kappa-\vartheta)
\|e^{\vartheta\eta}(\bar\phi^B-4\mathbf T^3g^B)\|_2^2\\
&\quad\leq\frac12|p_b|_{L^2(\gamma_-)}^2
+\int_0^B e^{2\vartheta\eta}
\langle(\phi^B-4\mathbf T^3g^B)S\rangle\,\mathrm d\eta\\
&\quad=\frac12|p_b|_{L^2(\gamma_-)}^2
+\int_0^B e^{2\vartheta\eta}
\langle(\phi^B-\bar\phi^B)(S-\bar S)\rangle\,\mathrm d\eta\\
&\qquad\quad+4\pi\int_0^B e^{2\vartheta\eta}
(\bar\phi^B-4\mathbf T^3g^B)\bar S\,\mathrm d\eta\\
&\quad\leq\frac12|p_b|_{L^2(\gamma_-)}^2
+\frac{1-\vartheta}{2}
\|e^{\vartheta\eta}(\phi^B-\bar\phi^B)\|_2^2\\
&\qquad\quad+2\pi(1-\kappa-\vartheta)
\|e^{\vartheta\eta}(\bar\phi^B-4\mathbf T^3g^B)\|_2^2\\
&\qquad\quad+\frac{\|e^{\vartheta\eta}(S-\bar S)\|_2^2}{2(1-\vartheta)}
+\frac{2\pi}{1-\kappa-\vartheta}\|e^{\vartheta\eta}\bar S\|_2^2,
\end{align*}
and this yields 
\begin{align*}
&(1-\vartheta)\|e^{\vartheta\eta}(\phi^B-\bar\phi^B)\|_2^2
+4\pi(1-\kappa-\vartheta)
\|e^{\vartheta\eta}(\bar\phi^B-4\mathbf T^3g^B)\|_2^2\\
&\quad\leq |p_b|_{L^2(\gamma_-)}^2
+\frac{\|e^{\vartheta\eta}(S-\bar S)\|_2^2}{1-\vartheta}
+\frac{4\pi}{1-\kappa-\vartheta}\|e^{\vartheta\eta}\bar S\|_2^2.
\end{align*}
Dividing by $1-\vartheta$ and discarding the second term on the left,
we pass to the half-line limit as in Section~\ref{linhalf}.
Weak lower semicontinuity on each fixed interval, followed by
exhaustion of $\mathbb R_+$, gives
	\begin{align}\label{mar}
		\|e^{\vartheta\eta}(\phi-\bar\phi)\|_2^2
		&\leq\frac{|p_b|_{L^2(\gamma_-)}^2}{1-\vartheta}
		+\frac{\|e^{\vartheta\eta}(S-\bar S)\|_2^2}{(1-\vartheta)^2}
		+\frac{4\pi\|e^{\vartheta\eta}\bar S\|_2^2}
		{(1-\vartheta)(1-\kappa-\vartheta)}.
	\end{align}
By~\eqref{5.31-1},
	\[
	\|e^{\vartheta\eta}\partial_\eta g\|_2
	=\|e^{\vartheta\eta}\langle\beta_3(\phi-\bar\phi)\rangle\|_2
	\leq\sqrt{\frac{4\pi}{3}}\,
	\|e^{\vartheta\eta}(\phi-\bar\phi)\|_2.
	\]
	Integrating $\partial_\eta g$ from $0$ or from infinity and applying the Cauchy--Schwarz inequality, we obtain
	\begin{align}\label{mtr}
		&|g_\infty|+\|g\|_\infty
		+\|e^{\vartheta\eta}(g-g_\infty)\|_\infty
		+\vartheta^{-1/2}\|e^{\vartheta\eta}\partial_\eta g\|_2\nonumber\\
		&\quad\leq\frac{C}{\sqrt{\vartheta}}\left[
		\frac{|p_b|_{L^2(\gamma_-)}}{\sqrt{1-\vartheta}}
		+\frac{\|e^{\vartheta\eta}(S-\bar S)\|_2}{1-\vartheta}
		+\frac{\|e^{\vartheta\eta}\bar S\|_2}
		{\sqrt{(1-\vartheta)(1-\kappa-\vartheta)}}\right].
	\end{align}
	In particular, when $\bar S=0$, the source contribution has no factor $(1-\kappa-\vartheta)^{-1/2}$.
	
	We first make the temperature boundary condition homogeneous by setting $\mathfrak{h}=g-g_b e^{-\sqrt{1-\kappa}\eta}$. Substituting this expression into \eqref{5.1}, we obtain
	\begin{align}\label{5.50}
		\begin{cases}
			\beta_3\partial_\eta\phi+\kappa(\phi-\bar{\phi})-(1-\kappa)(4\mathbf{T}^3\mathfrak{h}-\phi)=S_1+(1-\kappa)4\mathbf{T}^3g_be^{-\sqrt{1-\kappa}\eta}=:\tilde{S}_1,\\
			\partial_\eta^2\mathfrak{h}+(1-\kappa)4\pi(\bar{\phi}-4\mathbf{T}^3\mathfrak{h})=S_2+(1-\kappa)g_be^{-\sqrt{1-\kappa}\eta}(4\mathbf{T}^3\cdot4\pi-1)=:\tilde{S}_2.
		\end{cases}
	\end{align}
	The boundary conditions are
	\[
	\phi(0,\beta)|_{\beta_3>0}=\phi_b(\beta),\qquad \mathfrak h(0)=0.
	\]
	Since $\vartheta<1-\kappa\leq\sqrt{1-\kappa}$, the lifted sources belong to the weighted spaces in Theorem \ref{thm4}. In particular,
	\begin{align}\label{lsb}
		\|e^{\vartheta\eta}\tilde S_2\|_{L^2}
		&\leq \|e^{\vartheta\eta}S_2\|_{L^2}
		+\frac{C_{a_2}(1-\kappa)|g_b|}
		{\sqrt{\sqrt{1-\kappa}-\vartheta}},\\
		\langle\tilde S_1\rangle-\tilde S_2
		&=\langle S_1\rangle-S_2
		+(1-\kappa)g_b e^{-\sqrt{1-\kappa}\eta}.\label{lsd}
	\end{align}
	Here and below, $C_{a_2}$ depends only on the upper bound $a_2$ for $\mathbf T$.
	
	We first subtract the angularly constant part of the source difference. Theorem \ref{thm4} gives a unique bounded solution $(\mathfrak p_1,\mathfrak h_1)$ to
	\begin{align}\label{5.51}
		\begin{cases}
			\beta_3\partial_\eta\mathfrak{p}_1+\kappa(\mathfrak{p}_1-\bar{\mathfrak{p}}_1)-(1-\kappa)(4\mathbf{T}^3\mathfrak{h}_1-\mathfrak{p}_1)=\tilde{S}_1-\f{1}{4\pi}\bigl(\langle\tilde{S}_1\rangle-\tilde{S}_2\bigr),\\
			\partial_\eta^2\mathfrak{h}_1+(1-\kappa)4\pi(\bar{\mathfrak{p}}_1-4\mathbf{T}^3\mathfrak{h}_1)=\tilde{S}_2,
		\end{cases}
	\end{align}
	with $\mathfrak p_1(0,\beta)|_{\beta_3>0}=\phi_b(\beta)$ and $\mathfrak h_1(0)=0$.
	The transport source has angular average $\tilde S_2/(4\pi)$, and its part with zero angular average is $S_1-\bar S_1$, since the lifting term in $\tilde S_1$ is independent of $\beta$.
	Applying \eqref{mtr}, we obtain a far-field constant $\mathfrak h_{1,\infty}$ such that
	\begin{align}\label{gfp}
		&|\mathfrak h_{1,\infty}|+\|\mathfrak h_1\|_\infty
		+\|e^{\vartheta\eta}(\mathfrak h_1-\mathfrak h_{1,\infty})\|_\infty
		+\vartheta^{-1/2}\|e^{\vartheta\eta}\partial_\eta\mathfrak h_1\|_2\nonumber\\
		&\quad\leq\frac{C}{\sqrt{\vartheta}}\left[
		\frac{|\phi_b|_{L^2(\gamma_-)}}{\sqrt{1-\vartheta}}
		+\frac{\|e^{\vartheta\eta}(S_1-\bar S_1)\|_2}{1-\vartheta}
		+\frac{\|e^{\vartheta\eta}\tilde S_2\|_2}
		{\sqrt{(1-\vartheta)(1-\kappa-\vartheta)}}\right]
		\leq C_{a_2}\mathcal N_\vartheta.
	\end{align}
	For the last inequality, we used \eqref{lsb},
	$1-\kappa-\vartheta\leq1-\vartheta$, and
	$\sqrt{1-\kappa}-\vartheta\geq1-\kappa-\vartheta$.
	Here $\mathcal N_\vartheta$ is the quantity defined in \eqref{5.01-1}.
	
	To recover the subtracted source term, let $\mathfrak{q}(\eta,\beta)=-\beta_3G(\eta)$, where
	\begin{align}\label{gfun}
		G(\eta)
		&=\frac3{4\pi}\int_\eta^\infty
		\bigl(\langle\tilde S_1\rangle-\tilde S_2\bigr)(z)\,\mathrm dz\nonumber\\
		&=\frac3{4\pi}\int_\eta^\infty
		\bigl(\langle S_1\rangle-S_2\bigr)(z)\,\mathrm dz
		+\frac3{4\pi}\sqrt{1-\kappa}\,g_b e^{-\sqrt{1-\kappa}\eta}.
	\end{align}
	Since $\langle\beta_3\rangle=0$ and $\langle\beta_3^2\rangle=4\pi/3$, we have
	\begin{align*}
		\langle \mathfrak{q}\rangle&=0,\quad
		\left\langle\beta_3\partial_\eta \mathfrak{q}+\kappa(\mathfrak{q}-\bar{\mathfrak{q}})+(1-\kappa)\mathfrak{q}\right\rangle
		=-\frac{4\pi}{3}G'
		=\langle\tilde S_1\rangle-\tilde S_2.
	\end{align*}
	Minkowski's inequality, followed by integration of $e^{-\vartheta z}$, gives
	\begin{align}\label{ggb}
		\|e^{\vartheta\eta}G\|_2
		&\leq\frac3{4\pi\vartheta}
		\|e^{\vartheta\eta}(\langle S_1\rangle-S_2)\|_2
		+\frac{C\sqrt{1-\kappa}|g_b|}
		{\sqrt{\sqrt{1-\kappa}-\vartheta}},\nonumber\\
		|G(0)|
		&\leq\frac3{4\pi\sqrt{2\vartheta}}
		\|e^{\vartheta\eta}(\langle S_1\rangle-S_2)\|_2
		+\frac3{4\pi}\sqrt{1-\kappa}|g_b|,\nonumber\\
		\|e^{\vartheta\eta}G\|_\infty
		&\leq\frac3{4\pi\vartheta}
		\|e^{\vartheta\eta}(\langle S_1\rangle-S_2)\|_\infty
		+\frac3{4\pi}\sqrt{1-\kappa}|g_b|.
	\end{align}
	
	We next solve
	\begin{align}\label{5.52}
		\begin{cases}
			\beta_3\partial_\eta\mathfrak{p}_2+\kappa(\mathfrak{p}_2-\bar{\mathfrak{p}}_2)-(1-\kappa)(4\mathbf{T}^3\mathfrak{h}_2-\mathfrak{p}_2)=\f{1}{4\pi}\bigl(\langle\tilde{S}_1\rangle-\tilde{S}_2\bigr)-\left[\beta_3\partial_\eta \mathfrak{q}+\kappa(\mathfrak{q}-\bar{\mathfrak{q}})+(1-\kappa)\mathfrak{q}\right],\\
			\partial_\eta^2\mathfrak{h}_2+(1-\kappa)4\pi(\bar{\mathfrak{p}}_2-4\mathbf{T}^3\mathfrak{h}_2)=0,
		\end{cases}
	\end{align}
	with $\mathfrak p_2(0,\beta)|_{\beta_3>0}=\beta_3G(0)$ and $\mathfrak h_2(0)=0$.
	Using $G'=-3(\langle\tilde S_1\rangle-\tilde S_2)/(4\pi)$, the transport source becomes
	\[
	\frac{1-3\beta_3^2}{4\pi}
	(\langle\tilde S_1\rangle-\tilde S_2)+\beta_3G.
	\]
	It has zero angular average and belongs to the weighted $L^2$ and $L^\infty$ spaces by \eqref{lsd} and \eqref{ggb}. Thus Theorem \ref{thm4} applies. Then \eqref{mtr} gives
	\begin{align}\label{gsp}
		&|\mathfrak h_{2,\infty}|+\|\mathfrak h_2\|_\infty
		+\|e^{\vartheta\eta}(\mathfrak h_2-\mathfrak h_{2,\infty})\|_\infty
		+\vartheta^{-1/2}\|e^{\vartheta\eta}\partial_\eta\mathfrak h_2\|_2\nonumber\\
		&\quad\leq\frac{C}{\sqrt{\vartheta}}\left[
		\frac{|G(0)|}{\sqrt{1-\vartheta}}
		+\frac{\|e^{\vartheta\eta}(\langle\tilde S_1\rangle-\tilde S_2)\|_2
			+\|e^{\vartheta\eta}G\|_2}{1-\vartheta}\right]
		\leq C_{a_2}\mathcal N_\vartheta.
	\end{align}
	To verify the last inequality without losing an extra power of $1-\kappa$, we use
	\[
	\frac{\sqrt{1-\kappa}}{\sqrt{\sqrt{1-\kappa}-\vartheta}}
	\leq\frac{1-\kappa}{1-\kappa-\vartheta}.
	\]
	Indeed, $(1-\kappa-\vartheta)^2
	\leq(1-\kappa)(1-\kappa-\vartheta)
	\leq(1-\kappa)(\sqrt{1-\kappa}-\vartheta)$.
	Together with $\vartheta<1$ and \eqref{ggb}, this bounds the source-difference contribution by
	$C\vartheta^{-3/2}(1-\vartheta)^{-1}
	\|e^{\vartheta\eta}(\langle S_1\rangle-S_2)\|_2$
	and the lifting contribution by
	$C(1-\kappa)|g_b|/
	[\sqrt{\vartheta}(1-\vartheta)(1-\kappa-\vartheta)]$,
	as in \eqref{5.01-1}.
	
	Set
	\[
	\phi=\mathfrak p_1+\mathfrak p_2+\mathfrak{q},\qquad
	g=g_b e^{-\sqrt{1-\kappa}\eta}+\mathfrak h_1+\mathfrak h_2,\qquad
	g_\infty=\mathfrak h_{1,\infty}+\mathfrak h_{2,\infty}.
	\]
	Adding \eqref{5.51} and \eqref{5.52}, and using $\langle\mathfrak q\rangle=0$, shows that $(\phi,g)$ solves \eqref{5.1}--\eqref{5.1-2}.
	Equations \eqref{gfp} and \eqref{gsp} give the stated bounds for $g$ and $g_\infty$. For the derivative of the lifting,
	\[
	\|e^{\vartheta\eta}\partial_\eta(g_b e^{-\sqrt{1-\kappa}\eta})\|_2
	=\frac{\sqrt{1-\kappa}|g_b|}
	{\sqrt{2(\sqrt{1-\kappa}-\vartheta)}}
	\leq\frac{(1-\kappa)|g_b|}
	{\sqrt2(1-\kappa-\vartheta)}.
	\]
	This also proves \eqref{gg}.
	
	It remains to estimate the radiation. Applying the characteristic formula to the full solution, with source $S_1$, gives
	\[
	\|\phi\|_\infty
	\leq |\phi_b|_{L^\infty(\gamma_-)}
	+4a_2^3\|g\|_\infty+\frac{\|S_1\|_\infty}{1-\kappa}.
	\]
	Moreover,
	\[
	\|e^{\vartheta\eta}(4\mathbf T^3g-4T_{0,\infty}^3g_\infty)\|_\infty
	\leq C_{a_2}\left(a_2^3+
	\|e^{\vartheta\eta}(\mathbf T^3-T_{0,\infty}^3)\|_\infty\right)
	\mathcal N_\vartheta.
	\]
	By \eqref{4.56-1} with $\alpha=\vartheta$ and $0\leq\mathbf T,T_{0,\infty}\leq a_2$, we have
	\[
	\|e^{\vartheta\eta}(\mathbf T^3-T_{0,\infty}^3)\|_\infty
	\leq 3a_2^2\|e^{\vartheta\eta}(\mathbf T-T_{0,\infty})\|_\infty
	\leq\frac{3a_2^2}{\sqrt{2\vartheta}}
	\left(\frac{2\pi}{3(1-\kappa-\vartheta)}
	\int_{\beta_3>0}\beta_3(I_b-T_b^4)^2\,\mathrm d\beta\right)^{1/2}.
	\]
	Repeating the characteristic argument leading to~\eqref{lrad}, with source $S_1$ and the preceding temperature bounds, and absorbing the factor $\kappa/(1-\vartheta)$, we obtain the radiation bound in~\eqref{5.01} with $\mathcal M_\vartheta$ defined in~\eqref{5.01-2}.
	The difference of two bounded solutions with the same data satisfies
\eqref{5.33}--\eqref{5.34}, so the preceding uniqueness argument applies.

We now prove the difference estimates.
Since the two solutions have the same background
$\mathbf T$ and scattering parameter $\kappa$,
subtracting~\eqref{5.1}--\eqref{5.1-2} shows that
$(\phi_1-\phi_2,g_1-g_2)$ solves the same linear problem
with data
\[
(\phi_b,g_b,S_1,S_2)
=
(\phi_{b,1}-\phi_{b,2},\,
g_{b,1}-g_{b,2},\,
S_{1,1}-S_{1,2},\,
S_{2,1}-S_{2,2}).
\]
Its far-field temperature is
$g_{\infty,1}-g_{\infty,2}$.
Applying the estimates established above with this
data substitution gives~\eqref{ldt}--\eqref{ldd},
with the corresponding quantities
$\mathcal N_{\vartheta,12}$ and $\mathcal M_{\vartheta,12}$
defined in the statement.
This completes the proof of Theorem~\ref{co}.
	$\hfill\Box$

\subsection{Uniqueness of the nonlinear problem}Inspired by \cite{GHM-2023-ARMA}, we  have the following lemma.
\begin{lemma}[\cite{GHM-2023-ARMA}]
	For $T\in C^1(\mathbb{R}_+)$, if 
	\begin{align}\label{5.99}
		A_0:=\sup\limits_{r\in(0,\infty)}\left(\int_r^\infty e^{2\vartheta\eta}36T|\partial_\eta T|^2\,\mathrm{d}\eta\right)^{\f12}\cdot \left(\int_0^r\f{1}{e^{2\vartheta\eta}4T^3}\, \mathrm{d}\eta\right)^{\f12}<\f12,
	\end{align}
	then $T$ satisfies the spectral assumption (A). Moreover, if the spectral assumption (A) holds, then $A_0<1$.
\end{lemma}
The following sufficient condition for the spectral assumption follows by
adapting the argument in \cite{GHM-2023-ARMA}.
\begin{lemma}\label{lem-sp}
	Assume that, for some sufficiently small constant $C_b$,
	\begin{align}\label{4.62}
		\f12\int_{\beta_3>0}\beta_3(I_b-T_b^4)^2\,\mathrm{d}\beta
		\leq C_b^2.
	\end{align}
	Let $(\mathbf{I},\mathbf{T})$ be the
	corresponding solution obtained in Theorem \ref{thm2}. Then $\mathbf{T}$
	satisfies the spectral assumption (A).
\end{lemma}
\begin{proof}
	By \eqref{4.00} and \eqref{4.60}, for
	$0<\vartheta<\alpha<1-\kappa$,
	\begin{align}\label{5.100}
		\int_r^\infty e^{2\vartheta\eta}36\mathbf{T}
		|\partial_\eta \mathbf{T}|^2\,\mathrm{d}\eta
		\leq \frac{Ca_2}{\alpha-\vartheta}
		e^{-2(\alpha-\vartheta)r}
		\frac{(1-\kappa)^2}{\alpha(1-\kappa-\alpha)}C_b^2.
	\end{align}
	Meanwhile, 
	\begin{align*}
		|\mathbf{T}(\eta)-T_b|&=\left|\int_0^\eta \partial_z \mathbf{T} \mathrm{dz}\right|\leq \left(\int_0^\eta e^{2\vartheta z}|\partial_z \mathbf{T}|^2\,\mathrm{dz}\right)^{\f12}\cdot\left(\int_0^\eta e^{-2\vartheta z}\, dz\right)^{\f12}\nonumber\\
		&\leq
		\frac{C(1-\kappa)}
		{\sqrt{\alpha\vartheta(\alpha-\vartheta)
		(1-\kappa-\alpha)}}C_b.
	\end{align*}
	If
	\[
	C_b\leq c\,
	\frac{\sqrt{\alpha\vartheta(\alpha-\vartheta)
	(1-\kappa-\alpha)}}{1-\kappa}T_b
	\]
	for a sufficiently small universal constant $c>0$, then
	\begin{align*}
		|\mathbf{T}(\eta)-T_b|\leq \f12T_b,
	\end{align*}
	and therefore $
		\mathbf{T}(\eta)\geq \f12T_b.
$ Consequently,
	\begin{align}\label{5.101}
		\int_0^r\f{1}{e^{2\vartheta\eta}4\mathbf{T}^3}\, \mathrm{d}\eta\leq \frac{2}{T_b^3}\int_0^r e^{-2\vartheta\eta}\, \mathrm{d}\eta\leq \frac{1}{\vartheta T_b^3}.
	\end{align}
	
	Combining \eqref{5.100} and \eqref{5.101} gives
	\begin{align*}
		&\left(\int_r^\infty e^{2\vartheta\eta}36\mathbf{T}
		|\partial_\eta \mathbf{T}|^2\,\mathrm{d}\eta\right)
		\left(\int_0^r\f{1}{e^{2\vartheta\eta}4\mathbf{T}^3}\,
		\mathrm{d}\eta\right)\\
		&\qquad\leq
		\frac{Ca_2(1-\kappa)^2}
		{\vartheta(\alpha-\vartheta)\alpha
		(1-\kappa-\alpha)T_b^3}C_b^2.
	\end{align*}
	Thus \eqref{5.99} hold
	provided
	\[
	C_b\leq c_{a_2,T_b}
	\frac{\sqrt{\alpha\vartheta(\alpha-\vartheta)
	(1-\kappa-\alpha)}}{1-\kappa},
	\]
	where $c_{a_2,T_b}>0$ is independent of $\alpha$, $\vartheta$, and $\kappa$.
	Hence $\mathbf T$ satisfies the spectral assumption.
\end{proof}
\begin{theorem}[Uniqueness]\label{thm3}
	Assume the hypotheses of Theorem \ref{thm2} hold. Suppose $(\mathbf{I},\mathbf{T})$ is a solution to system \eqref{3.16}--\eqref{3.17}. Then if $\mathbf{T}$ satisfies the spectral assumption (A), the solution to system \eqref{3.16}--\eqref{3.17} is unique in the set 
	\begin{align*}
		\mathcal{V}_\v(\mathbf{T})=\left\{(\mathbf{I},G)\in L^2_{\text{loc}}(\R_+\times\mathbb{S}^2)\times L_{\text{loc}}^2(\R_+):\|G-\mathbf{T}\|_{L^2(\R_+)}\leq \v\right\}
	\end{align*}
	for $\v>0$ sufficiently small.
	
	Moreover, if the boundary data satisfy \eqref{4.62} with $C_b$
	sufficiently small, then the solution to
	\eqref{3.16}--\eqref{3.17} is unique.
\end{theorem}
%
\textbf{Proof of Theorem \ref{thm3}}. We prove Theorem \ref{thm3} for the case when boundary data is almost well-prepared. Let $(\mathbf{I},\mathbf{T} )$ and $(\mathbf{I}_1,\mathbf{T}_1)$ be two solutions to system \eqref{3.16} with boundary conditions \eqref{3.17} satisfying \eqref{4.62}. Then $\mathfrak{p}=\mathbf{I}_1-\mathbf{I}$ and $\mathfrak{g}=\mathbf{T}_1-\mathbf{T}$ satisfy
\begin{align}\label{4.64}
	\begin{cases}
		\beta_3\partial_\eta \mathfrak{p}+\kappa(\mathfrak{p}-\bar{\mathfrak{p}})=(1-\kappa)(\mathfrak{s}\mathfrak{g}-\mathfrak{p}),\\
		-\partial_\eta^2\mathfrak{g}=(1-\kappa)4\pi(\bar{\mathfrak{p}}-\mathfrak{s}\mathfrak{g}),
	\end{cases}
\end{align}
with boundary conditions 
\begin{align}\label{4.65}
	\mathfrak{p}(0,\beta)|_{\beta_3>0}=0,\quad \mathfrak{g}(0)=0,
\end{align}
where $\mathfrak{s}=(\mathbf{T}_1^4-\mathbf{T}^4)/(\mathbf{T}_1-\mathbf{T})=(\mathbf{T}+\mathbf{T}_1)(\mathbf{T}^2+\mathbf{T}_1^2)$. By \cite{GHM-2023-ARMA}, we know that $\mathfrak{s}^{\f13}$ also satisfies the spectral assumption by choosing suitably small $C_b$. In detail, it is sufficient to prove that  
\begin{align}
	\left(\int_r^\infty e^{2\vartheta\eta}\cdot\frac{1}{\mathfrak{s}}(\partial_\eta \mathfrak{s})^2\, \mathrm{d}\eta\right)\cdot\left(\int_0^r\frac{1}{e^{2\vartheta\eta}\mathfrak{s}}\mathrm{d}\eta\right)\, <\f14.
\end{align}

We choose $C_b$ small such that 
\begin{align*}
	|\mathbf{T}(\eta)-T_b|\leq\f12 T_b, \quad |\mathbf{T}_1(\eta)-T_b|\leq \f12 T_b.
\end{align*}
A direct calculation shows
\begin{align*}
	\partial_\eta \mathfrak{s}=(\partial_\eta \mathbf{T}+\partial_\eta \mathbf{T}_1)(\mathbf{T}^2+\mathbf{T}_1^2)+2(\mathbf{T}+\mathbf{T}_1)(\mathbf{T}\partial_\eta\mathbf{T}+\mathbf{T}_1\partial_\eta\mathbf{T}_1),
\end{align*}
Applying $(a+b)^2\leq2a^2+2b^2$ to the two terms in
$\partial_\eta\mathfrak s$ and then applying Cauchy--Schwarz inequality, we obtain
\begin{align*}
	\frac{|\partial_\eta\mathfrak s|^2}{\mathfrak s}
	&\leq
	2\frac{\mathbf T^2+\mathbf T_1^2}{\mathbf T+\mathbf T_1}
	(\partial_\eta\mathbf T+\partial_\eta\mathbf T_1)^2
	+8\frac{\mathbf T+\mathbf T_1}{\mathbf T^2+\mathbf T_1^2}
	(\mathbf T\partial_\eta\mathbf T+\mathbf T_1\partial_\eta\mathbf T_1)^2\\
	&\leq\left(
	4\frac{\mathbf T^2+\mathbf T_1^2}{\mathbf T+\mathbf T_1}
	+8(\mathbf T+\mathbf T_1)\right)
	\left(|\partial_\eta\mathbf T|^2+|\partial_\eta\mathbf T_1|^2\right)\\
	&\leq30T_b
	\left(|\partial_\eta\mathbf T|^2+|\partial_\eta\mathbf T_1|^2\right),
\end{align*}
where we used
$(\mathbf T^2+\mathbf T_1^2)/(\mathbf T+\mathbf T_1)
\leq\max\{\mathbf T,\mathbf T_1\}\leq3T_b/2$ and
$\mathbf T+\mathbf T_1\leq3T_b$.
Then similar to \eqref{5.100} and \eqref{5.101}, we choose $C_b<\delta \sqrt{\vartheta(\alpha-\vartheta)}\cdot\frac{\sqrt{\alpha(1-\kappa-\alpha)}}{1-\kappa}T_b$ with $\delta\ll1$ independent of $\alpha,\vartheta$ or $\kappa$ such that $\mathfrak{s}^{\f13}$ satisfies the spectral assumption. 

Applying the uniqueness argument for \eqref{5.33}--\eqref{5.34}, with $4\mathbf T^3$ replaced by $\mathfrak s$, gives $\mathfrak p=0$ and $\mathfrak g=0$ in \eqref{4.64}--\eqref{4.65}. We further deduce $\mathbf{I}_1=\mathbf{I}$ and $\mathbf{T}_1=\mathbf{T}$, then the solution to \eqref{3.16}--\eqref{3.17} is unique. 

On the other hand, if $\mathbf{T}$ satisfies the spectral assumption, by similar arguments in \cite{GHM-2023-ARMA}, the solution to \eqref{3.16}--\eqref{3.17} is unique in the set $\mathcal{V}_\v$. Therefore Theorem \ref{thm3} holds.
$\hfill\Box$
\section{Construction of the matched profiles}\label{prf}

Throughout this section, $D^w=D^w_{x_\shortparallel}$ and
$a=T_{0,\infty}$ denotes the far-field temperature of the lower
nonlinear Milne profile. For an integer $m\geq0$, fix decay parameters
\begin{align}\label{6.0-1}
 1>\lambda_{0,0}>\lambda_{0,1}>\cdots>\lambda_{0,m}>0,
\end{align}
independent of $\kappa$.

\begin{lemma}[Tangential regularity of the nonlinear Milne profile]\label{ffr}
Assume that the lower boundary data satisfy
\[
 \sum_{|w|\leq m}\left(
 \|D^w I_b(\cdot,0,\cdot)\|_{L^\infty(\mathbb T^2\times\{\beta_3>0\})}
 +\|D^w T_b(\cdot,0)\|_{L^\infty(\mathbb T^2)}\right)<\infty.
\]
Let $(\mathbf I,\mathbf T)$ be a solution of \eqref{3.16} furnished by
Theorem~\ref{thm2}.
Assume that, for every $x_\shortparallel\in\mathbb T^2$,
$\mathbf T(\cdot,x_\shortparallel)$ satisfies the spectral
assumption~\eqref{asA} for each weight
$\vartheta=\lambda_{0,r}(1-\kappa)$, $r=1,\ldots,m$, with
$M\leq M_*<1$ uniformly in $x_\shortparallel$ and $r$.
Then $a\in W^{m,\infty}(\mathbb T^2)$. For every multi-index $w$
with $1\leq |w|=r\leq m$,
    \begin{align}\label{6.0-0}
        &|D^wa|_{L^\infty}
        +\|e^{\lambda_{0,r}(1-\kappa)\eta}D^w(\mathbf T-a)\|_{L^\infty}
        \leq \frac{C_r}{(1-\kappa)^{\frac72r-2}},\\
        &\|e^{\lambda_{0,r}(1-\kappa)\eta}D^w(\mathbf I-a^4)\|_{L^\infty}
        \leq \frac{C_r}{(1-\kappa)^{\frac72r-1}}.\label{rad}
    \end{align}
    The constants $C_r$ do not depend on $\kappa$.
\end{lemma}
\begin{proof} We divide the proof into several steps.
	
    \noindent{\it Step 1. $|w|=0$.} Taking
    $\alpha=\lambda_{0,0}(1-\kappa)$ in Theorem \ref{thm2}, one has
	\begin{align}
		\begin{split}
			&e^{\lambda_{0,0}(1-\kappa)\eta}|\mathbf{I}-T_{0,\infty}^4|\leq \frac{C_{\lambda_{0,0}}}{1-\kappa}|I_b-T_b^4|,\\
			&e^{\lambda_{0,0}(1-\kappa)\eta}|\mathbf{T}-T_{0,\infty}|\leq\frac{C_{\lambda_{0,0}}}{1-\kappa}|I_b-T_b^4|,
		\end{split}
	\end{align} 
	where $C_{\lambda_{0,0}}$ does not depend on $\kappa$. Similarly, taking $\vartheta=\tilde{b}(1-\kappa)$ in Theorem \ref{co} where $0<\tilde{b}<\lambda_{0,0}$, we get
	\begin{align}\label{le}
&|g_\infty|+\|g\|_{L^\infty}
+\|e^{\tilde b(1-\kappa)\eta}(g-g_\infty)\|_{L^\infty}
\nonumber\\
&\qquad\leq\frac{C_{\tilde b}}{(1-\kappa)^{3/2}}
\left(|\phi_b|_{L^2(\gamma_-)}+|g_b|
+\|e^{\tilde b(1-\kappa)\eta}(S_1,S_2)\|_{L^2}\right),\\
&\|\phi\|_{L^\infty}
\leq\frac{C_{\tilde b}}{(1-\kappa)^{3/2}}
\left(|\phi_b|_{L^\infty(\gamma_-)}+|g_b|
+\|e^{\tilde b(1-\kappa)\eta}(S_1,S_2)\|_{L^2}\right)
+\frac{\|S_1\|_{L^\infty}}{1-\kappa},\nonumber\\
&\|e^{\tilde b(1-\kappa)\eta}
(\phi-4T_{0,\infty}^3g_\infty)\|_{L^\infty}
\nonumber\\
&\qquad\leq\frac{C_{\tilde b}}{(1-\kappa)^{5/2}}
\left(|\phi_b|_{L^\infty(\gamma_-)}+|g_b|
+\|e^{\tilde b(1-\kappa)\eta}(S_1,S_2)\|_{L^2}\right)
+\frac{C_{\tilde b}}{1-\kappa}
\|e^{\tilde b(1-\kappa)\eta}S_1\|_{L^\infty}.\nonumber
 \end{align}
 The constants are independent of $\kappa$. 

 \noindent{\it Step 2. $|w|= 1$.}  Set
	\begin{align*}
		a(x_{\shortparallel})=T_{0,\infty}(x_{\shortparallel}),\quad U(\eta,x_{\shortparallel},\beta)=\mathbf I(\eta,x_{\shortparallel},\beta)-a^4(x_{\shortparallel}),\quad V(\eta,x_{\shortparallel})=\mathbf T(\eta,x_{\shortparallel})-a(x_{\shortparallel}).
	\end{align*}
	The variables $(U,V)$ satisfy 
	\begin{align}\label{6.a}
		\begin{cases}
		\beta_3\partial_\eta U+\kappa(U-\bar U)=(1-\kappa)((a+V)^4-a^4-U),\\
		-\partial_\eta^2V=(1-\kappa)4\pi\left(\bar U-((a+V)^4-a^4)\right),
		\end{cases}
	\end{align}
	with boundary data
	\begin{align}\label{6.b}
		U(0,x_{\shortparallel},\beta)|_{\gamma_-}=I_b(x_{\shortparallel},\beta)-a^4(x_{\shortparallel}),\quad V(0,x_{\shortparallel})=T_b(x_{\shortparallel})-a(x_{\shortparallel}),
	\end{align}
and $U,V\to0$ as $\eta\to\infty$.
For $m\geq1$, fix $j\in\{1,2\}$ and define
\[
 \delta_h f(x_\shortparallel)
 :=\frac{f(x_\shortparallel+h e_j)-f(x_\shortparallel)}{h},
 \qquad h\neq0,
\]
where $e_j$ is the $j$-th coordinate vector. Applying~\eqref{hshift}--\eqref{hfarstab} at
$x_\shortparallel+h e_j$ and $x_\shortparallel$, and dividing by
$|h|$, we obtain
\begin{align}\label{ffr-dq1}
 \|\delta_h\mathbf I\|_\infty
 +\|\delta_h\mathbf T\|_\infty+\|\delta_h a\|_\infty
 \leq C\bigl(
 \|\partial_{x_j}I_b\|_\infty+\|\partial_{x_j}T_b\|_\infty\bigr).
\end{align}
Passing weak-star to the limit yields
$\partial_{x_j}\mathbf I,\partial_{x_j}\mathbf T,\partial_{x_j}a
\in L^\infty$.
Differentiating the weak formulation of~\eqref{3.16}--\eqref{3.17}
and applying the Sobolev chain rule, we obtain~\eqref{5.1} for
$(\partial_{x_j}\mathbf I,\partial_{x_j}\mathbf T)$, with
$S_1=S_2=0$ and boundary data
$(\partial_{x_j}I_b,\partial_{x_j}T_b)$.
Theorem~\ref{co} gives a uniform far-field limit for
$\partial_{x_j}\mathbf T$. For $\zeta\in C^\infty(\mathbb T^2)$,
integration by parts and $\mathbf T(\eta,\cdot)\to a$ give
\[
 \lim_{\eta\to\infty}\int_{\mathbb T^2}
       \partial_{x_j}\mathbf T(\eta,x_\shortparallel)
       \zeta(x_\shortparallel)\,\mathrm dx_\shortparallel
 =-\int_{\mathbb T^2}a\,\partial_{x_j}\zeta\,\mathrm dx_\shortparallel.
\]
Thus $\partial_{x_j}\mathbf T\to\partial_{x_j}a$ and
$\partial_{x_j}\mathbf I\to4a^3\partial_{x_j}a$, with the
decay supplied by Theorem~\ref{co}.
For $|w|=1$, applying~\eqref{le} with $S_1=S_2=0$ and boundary
data $(D^wI_b,D^wT_b)$, we obtain
\begin{align}\label{cfo}
	\begin{split}
 |D^wa|+\|e^{\lambda_{0,1}(1-\kappa)\eta}D^wV\|_{L^\infty}
 &\leq \frac{C_1}{(1-\kappa)^{3/2}},\\
 \|e^{\lambda_{0,1}(1-\kappa)\eta}D^wU\|_{L^\infty}
 &\leq \frac{C_1}{(1-\kappa)^{5/2}}.
\end{split}
\end{align}

\noindent{\it Step 3. $|w|\geq 2$.}  We prove the higher order estimates by induction on $|w|$. Let $2\leq s\leq m$, and assume that the derivatives, their far-field
limits, and~\eqref{6.0-0}--\eqref{rad} are established through order
$s-1$. For $|w|=s$, set formally
	\begin{align*}
		U_w=D^w_{x_{\shortparallel}}U,\quad V_w=D^w_{x_{\shortparallel}}V,
		\quad a_w=D^w_{x_{\shortparallel}}a.
	\end{align*}

		Formally applying $D^w$ to~\eqref{6.a}, we obtain
	\begin{align}\label{6.e}
		\begin{cases}
			\beta_3\partial_\eta U_w+\kappa(U_w-\bar U_w)+(1-\kappa)(U_w-4\mathbf T^3V_w)
			=(1-\kappa)\left(4(\mathbf T^3-a^3)a_w+\mathcal Q_w\right),\\
			-\partial_\eta^2V_w-(1-\kappa)4\pi(\bar U_w-4\mathbf T^3V_w)
			=-(1-\kappa)4\pi\left(4(\mathbf T^3-a^3)a_w+\mathcal Q_w\right),
		\end{cases}
	\end{align}
	where
	\begin{align}\label{6.c}
		\mathcal Q_w:=D^w_{x_{\shortparallel}}\left((a+V)^4-a^4\right)-4\mathbf T^3V_w-4(\mathbf T^3-a^3)a_w.
	\end{align}
	Applying $D^w$ to~\eqref{6.b} gives
	\begin{align}\label{6.f}
		U_w(0,\beta)|_{\gamma_-}=D^w_{x_{\shortparallel}}I_b-4a^3a_w-\mathcal R_w^a,
		\quad
		V_w(0)=D^w_{x_{\shortparallel}}T_b-a_w,
	\end{align}
	where
	\begin{align*}
		\mathcal R_w^a=D^w_{x_{\shortparallel}}(a^4)-4a^3a_w
	\end{align*}
	contains derivatives of $a$ only through order $s-1$. Set
	\begin{align*}
		\widehat U_w:=U_w+4a^3a_w,
		\quad
		\widehat V_w:=V_w+a_w.
	\end{align*}
	Substituting these definitions into~\eqref{6.e}--\eqref{6.f}, we obtain
	\begin{align}\label{6.g}
		\begin{cases}
			\beta_3\partial_\eta\widehat U_w+\kappa(\widehat U_w-\overline{\widehat U}_w)
			=(1-\kappa)(4\mathbf T^3\widehat V_w-\widehat U_w)+(1-\kappa)\mathcal Q_w,\\
			\partial_\eta^2\widehat V_w+(1-\kappa)4\pi(\overline{\widehat U}_w-4\mathbf T^3\widehat V_w)
			=(1-\kappa)4\pi\mathcal Q_w,
		\end{cases}
	\end{align}
	with boundary data
	\begin{align}\label{6.h}
		\widehat U_w(0,\beta)|_{\gamma_-}=D^w_{x_{\shortparallel}}I_b-\mathcal R_w^a,
		\quad
		\widehat V_w(0)=D^w_{x_{\shortparallel}}T_b.
	\end{align}
	
Since $V_w=D^w\mathbf T-D^wa$, the highest-order terms cancel exactly:
	\begin{align*}
		\mathcal Q_w
		=\left[D^w(\mathbf T^4)-4\mathbf T^3D^w\mathbf T\right]
		-\left[D^w(a^4)-4a^3D^wa\right].
	\end{align*}
	Define the finite index set
	\begin{align*}
		\mathcal A_w:=\left\{(w^1,w^2,w^3,w^4)\in(\mathbb N_0^2)^4:
		w^1+w^2+w^3+w^4=w,\quad |w^j|\leq s-1\ \text{for all }j\right\}.
	\end{align*}
	The multi-index Leibniz formula gives
	\begin{align}\label{cl}
		\mathcal Q_w
		=\sum_{(w^1,w^2,w^3,w^4)\in\mathcal A_w}
		\frac{w!}{w^1!w^2!w^3!w^4!}
		\left(\prod_{j=1}^4D^{w^j}\mathbf T-\prod_{j=1}^4D^{w^j}a\right).
	\end{align}
	Indeed, the four excluded terms, in which one multi-index equals $w$ and the other three vanish, are precisely the highest-order terms subtracted above. If $s=1$, the set $\mathcal A_w$ is empty and $\mathcal Q_w=0$.
	For each tuple in $\mathcal A_w$, the product difference has the exact decomposition
	\begin{align}\label{ct}
		\prod_{j=1}^4D^{w^j}\mathbf T-\prod_{j=1}^4D^{w^j}a
		=\sum_{\ell=1}^4
		\left(\prod_{j<\ell}D^{w^j}\mathbf T\right)
		D^{w^\ell}V
		\left(\prod_{j>\ell}D^{w^j}a\right),
	\end{align}
    where an empty product equals one.

For $s\geq2$, fix $|v|=s-1$. By the induction hypothesis, the pair
\[
 \widehat U_v=D^v\mathbf I-\mathcal R_v^a,\qquad
 \widehat V_v=D^v\mathbf T,\qquad
 \mathcal R_v^a=D^v(a^4)-4a^3D^va
\]
satisfies~\eqref{6.g}--\eqref{6.h} with $w=v$.
Set
\[
 \Phi_h=\delta_h\widehat U_v
       -4\delta_h(a^3)D^va(x_\shortparallel+h e_j),
 \qquad H_h=\delta_h\widehat V_v.
\]
Applying $\delta_h$ to~\eqref{6.g} and substituting these definitions,
we obtain
\begin{align}\label{dqs}
 \begin{cases}
 \beta_3\partial_\eta\Phi_h+\kappa(\Phi_h-\bar\Phi_h)
 =(1-\kappa)(4\mathbf T^3H_h-\Phi_h)+(1-\kappa)F_h,\\
 \partial_\eta^2H_h+(1-\kappa)4\pi(\bar\Phi_h-4\mathbf T^3H_h)
 =4\pi(1-\kappa)F_h,
 \end{cases}
\end{align}
where $\mathbf T=\mathbf T(\eta,x_\shortparallel)$ and
\begin{align}\label{dqf}
 F_h
 &=\delta_h\mathcal Q_v
   +4\delta_h(\mathbf T^3)\widehat V_v(\eta,x_\shortparallel+h e_j)
   -4\delta_h(a^3)D^va(x_\shortparallel+h e_j)
   \nonumber\\
 &=\delta_h\mathcal Q_v
   +4\delta_h(\mathbf T^3)D^vV(\eta,x_\shortparallel+h e_j)
   +4\delta_h(\mathbf T^3-a^3)D^va(x_\shortparallel+h e_j).
\end{align}
Applying $\delta_h$ to~\eqref{6.h} gives
\[
 \begin{aligned}
 \Phi_h(0,\beta)
 &=\delta_h(D^vI_b-\mathcal R_v^a)(\beta)
   -4\delta_h(a^3)D^va(x_\shortparallel+h e_j),
 &&\beta_3>0,\\
 H_h(0)&=\delta_hD^vT_b.
 \end{aligned}
\]
The product expansions~\eqref{cl}--\eqref{ct} at order $v$
contain derivatives only through order $s-2$.
Using the induction hypothesis and the identities
\[
 \begin{aligned}
 \delta_h(fg)
 &=(\delta_h f)\,g(x_\shortparallel+h e_j)+f\,\delta_h g,\\
 \partial_{x_j}(\mathbf T^3-a^3)
 &=3\mathbf T^2\partial_{x_j}V
   +3(\mathbf T^2-a^2)\partial_{x_j}a,
 \end{aligned}
\]
we obtain
\[
 |\Phi_h(0)|_{L^\infty(\gamma_-)}+|H_h(0)|
 +\|e^{\lambda_{0,s-1}(1-\kappa)\eta}F_h\|_\infty
 \leq C_{\kappa,s},
\]
with $C_{\kappa,s}$ independent of $h$.
Since $\lambda_{0,s}<\lambda_{0,s-1}$,
\[
 \sup_{x_\shortparallel}
 \|e^{\lambda_{0,s}(1-\kappa)\eta}F_h\|_{L^2(\mathbb R_+)}
 \leq
 \frac{\|e^{\lambda_{0,s-1}(1-\kappa)\eta}F_h\|_\infty}
 {\sqrt{2(\lambda_{0,s-1}-\lambda_{0,s})(1-\kappa)}}.
\]
For fixed $h$, the far-field temperature of $H_h$ is $\delta_hD^va$,
and
\[
 H_h-\delta_hD^va=\delta_hD^vV,\qquad
 \Phi_h-4a^3\delta_hD^va=\delta_hD^vU.
\]
Applying Theorem~\ref{co} to~\eqref{dqs} with
$\vartheta=\lambda_{0,s}(1-\kappa)$ therefore gives
\begin{align}\label{dqb}
 \|\delta_hD^va\|_\infty
 +\|e^{\lambda_{0,s}(1-\kappa)\eta}\delta_hD^vV\|_\infty
 +\|e^{\lambda_{0,s}(1-\kappa)\eta}\delta_hD^vU\|_\infty
 \leq C_{\kappa,s}.
\end{align}
Passing weak-star to the limit gives the order-$s$ weak derivatives
and their weighted bounds. In particular, for $w=v+e_j$,
$D^w\mathbf T\to D^wa$ and $D^w\mathbf I\to D^w(a^4)$.
Passing to distributions in~\eqref{dqs} and applying the
Sobolev product rule proves~\eqref{6.e}--\eqref{6.h} for these
order-$s$ weak derivatives.

We now estimate the dependence on $\kappa$. The zeroth-order estimates \eqref{4.0} and \eqref{4.00} give
    \begin{align}\label{czb}
        |\mathbf T|,\ |a|&\leq a_2,\qquad
        |V(\eta,x_{\shortparallel})|
        \leq \frac{C_0}{1-\kappa}e^{-\lambda_{0,0}(1-\kappa)\eta}.
    \end{align}
    For $1\leq |v|=r\leq s-1$, the induction hypothesis implies
    \begin{align}\label{clb}
        |D^va|+|D^v\mathbf T|
        \leq \frac{C_r}{(1-\kappa)^{\frac72r-2}},\quad
        |D^vV(\eta,x_{\shortparallel})|
        \leq \frac{C_r}{(1-\kappa)^{\frac72r-2}}
        e^{-\lambda_{0,r}(1-\kappa)\eta}.
    \end{align}
    For a fixed tuple in $\mathcal A_w$, put $r_j=|w^j|$ and let $n$ be the number of positive $r_j$. Since $\sum_jr_j=s$ and $r_j\leq s-1$, we have $2\leq n\leq4$. If $r_\ell>0$, \eqref{czb}--\eqref{clb} yield
    \begin{align}\label{cot}
        &\left|\left(\prod_{j<\ell}D^{w^j}\mathbf T\right)
        D^{w^\ell}V\left(\prod_{j>\ell}D^{w^j}a\right)\right|\leq C_s e^{-\lambda_{0,r_\ell}(1-\kappa)\eta}
        \prod_{j:r_j>0}(1-\kappa)^{-(\frac72r_j-2)}\nonumber\\
        &\quad=\frac{C_s e^{-\lambda_{0,r_\ell}(1-\kappa)\eta}}
        {(1-\kappa)^{\frac72s-2n}}
        \leq\frac{C_s e^{-\lambda_{0,s-1}(1-\kappa)\eta}}
        {(1-\kappa)^{\frac72s-4}}.
    \end{align}
    If $r_\ell=0$, the decaying factor is $V$ itself. Its zeroth-order bound contributes one additional power of $(1-\kappa)^{-1}$, so that
    \begin{align}\label{czt}
        \left|\left(\prod_{j<\ell}D^{w^j}\mathbf T\right)
        V\left(\prod_{j>\ell}D^{w^j}a\right)\right|\leq\frac{C_s e^{-\lambda_{0,0}(1-\kappa)\eta}}
        {(1-\kappa)^{\frac72s-2n+1}}
        \leq\frac{C_s e^{-\lambda_{0,s-1}(1-\kappa)\eta}}
        {(1-\kappa)^{\frac72s-3}}.
    \end{align}
    Summing \eqref{cot} and \eqref{czt} with the multinomial coefficients in \eqref{cl} gives
    \begin{align}\label{6.d}
        \|e^{\lambda_{0,s-1}(1-\kappa)\eta}\mathcal Q_w\|_{L^\infty}
        \leq \frac{C_s}{(1-\kappa)^{\frac72s-3}}.
    \end{align}
    Here the constants involve only lower-order induction bounds and the finite combinatorial sums. Since $\lambda_{0,s}<\lambda_{0,s-1}$,
    \begin{align}\label{dw}
        &\sup_{x_{\shortparallel}}\|e^{\lambda_{0,s}(1-\kappa)\eta}
        \mathcal Q_w(\cdot,x_{\shortparallel})\|_{L^2(\mathbb R_+)}\nonumber\\
        &\quad\leq\frac{C_s}{(1-\kappa)^{\frac72s-3}}
        \left(\int_0^\infty e^{-2(\lambda_{0,s-1}-\lambda_{0,s})(1-\kappa)\eta}\,d\eta\right)^{1/2}\nonumber\\
        &\quad=\frac{C_s}{\sqrt{2(\lambda_{0,s-1}-\lambda_{0,s})}}
        \frac{1}{(1-\kappa)^{\frac72s-\frac52}}.
    \end{align}
    The boundary remainder $\mathcal R_w^a$ contains products of derivatives of $a$ with the same index set. It has no zeroth-order decaying factor and hence
    \begin{align}\label{db}
        |\mathcal R_w^a|_{L^\infty}
        \leq C_s\max_{2\leq n\leq4}(1-\kappa)^{-(\frac72s-2n)}
        \leq\frac{C_s}{(1-\kappa)^{\frac72s-4}}.
    \end{align}
    The sources in \eqref{6.g} are $S_1=(1-\kappa)\mathcal Q_w$ and $S_2=4\pi(1-\kappa)\mathcal Q_w$. Therefore, absorbing the fixed weight gap into $C_s$,
    \begin{align}\label{ds}
        \|e^{\lambda_{0,s}(1-\kappa)\eta}(S_1,S_2)\|_{L^\infty}
        &\leq\frac{C_s}{(1-\kappa)^{\frac72s-4}},\\
        \sup_{x_{\shortparallel}}\|e^{\lambda_{0,s}(1-\kappa)\eta}(S_1,S_2)\|_{L^2_{\eta,\beta}}
        &\leq\frac{C_s}{(1-\kappa)^{\frac72s-\frac72}}.
    \end{align}
    
    The explicit linear estimates \eqref{le}, applied through Theorem \ref{co} to \eqref{6.g}--\eqref{6.h}, give the far-field limits $\widehat V_w\to a_w$ and $\widehat U_w\to4a^3a_w$. Substituting \eqref{db}--\eqref{ds} into the temperature estimate in \eqref{le} gives
    \begin{align}\label{dct}
        &|a_w|+\|e^{\lambda_{0,s}(1-\kappa)\eta}V_w\|_{L^\infty}\nonumber\\
        &\quad\leq\frac{C_s}{(1-\kappa)^{3/2}}
        \left(1+\frac{1}{(1-\kappa)^{\frac72s-4}}
        +\frac{1}{(1-\kappa)^{\frac72s-\frac72}}\right)\nonumber\\
        &\quad\leq\frac{C_s}{(1-\kappa)^{\frac72s-2}}.
    \end{align}
    Here the norms of $D^wI_b$ and $D^wT_b$ have been included in $C_s$. Similarly, the radiation estimate in \eqref{le} yields
    \begin{align}\label{dcr}
        \|e^{\lambda_{0,s}(1-\kappa)\eta}U_w\|_{L^\infty}
        &\leq\frac{C_s}{(1-\kappa)^{5/2}}
        \left(1+\frac{1}{(1-\kappa)^{\frac72s-4}}
        +\frac{1}{(1-\kappa)^{\frac72s-\frac72}}\right)
        +\frac{C_s}{(1-\kappa)^{\frac72s-3}}\nonumber\\
        &\leq\frac{C_s}{(1-\kappa)^{\frac72s-1}}.
    \end{align}
    Together with the difference-quotient argument above,
    equations~\eqref{cfo}, \eqref{dct}, and~\eqref{dcr} close the
    induction on both the existence of the weak derivatives \eqref{clb}.
\end{proof}

\begin{lemma}[Elliptic estimates for the leading interior profile]\label{lem1}
Let $q>3$ and $m\geq0$. Assume the hypotheses of
Lemma~\ref{ffr} with tangential order $m+2$ and
$T_b(1,\cdot)\in W^{m+2,\infty}(\mathbb T^2)$.
Then~\eqref{3.5} with \eqref{Tb} has a unique solution
$T_0\in W^{m+2,q}(\Omega)$ satisfying $a_*\leq T_0\leq a_2$,
where $a_*$ is given by~\eqref{4.00}. For $0\leq s\leq m$,
\begin{align}\label{ell}
 \|T_0\|_{W^{s+2,q}(\Omega)}
 +\|T_0\|_{W^{s+1,\infty}(\Omega)}
 &\leq C_{s,q}(1-\kappa)^{-\frac72(s+2)+2},\\
 \max_{|w|\leq s}
 \|D^w_{x_{\shortparallel}}\partial_{x_3}T_0|_{x_3=0}\|_\infty
 &\leq C_{s,q}(1-\kappa)^{-\frac72s-5}.
 \label{en}
\end{align}
The constants are independent of $\kappa$.
\end{lemma}
\begin{proof}
Set $\mathcal H(z)=z+\frac{4\pi}{3}z^4$ and $Y_0=\mathcal H(T_0)$.
Then
\begin{align}\label{eh}
 \Delta Y_0=0,\qquad
 Y_0|_{x_3=0}=\mathcal H(a),\qquad
 Y_0|_{x_3=1}=\mathcal H(T_b).
\end{align}
The boundary bounds give
$\mathcal H(a_*)\leq Y_0|_{\partial\Omega}\leq\mathcal H(a_2)$.
The maximum principle gives the same bounds in $\Omega$.
Since $\mathcal H'\geq1$ on the nonnegative range, inversion gives
existence and uniqueness, with $a_*\leq T_0\leq a_2$.

Here is the dependence of the boundary norms on $\kappa$.
For $n=s+2$, the Leibniz formula for $D^v(a^4)$, $|v|=r\leq n$,
contains products of positive-order derivatives whose orders add up
to $r$. If there are $\ell$ such factors, \eqref{6.0-0} bounds their
product by $C_r(1-\kappa)^{-\frac72r+2\ell}$.
Undifferentiated factors are uniformly bounded. Thus
\begin{align}\label{ebd}
 \|\mathcal H(a)\|_{W^{n,\infty}(\mathbb T^2)}
 +\|\mathcal H(T_b(\cdot,1))\|_{W^{n,\infty}(\mathbb T^2)}
 \leq C_n(1-\kappa)^{-\frac72n+2}.
\end{align}
The $W^{n,q}$ Dirichlet estimate on the fixed periodic slab uses the
boundary norm $W^{n-1/q,q}$, which is bounded by the norm in
\eqref{ebd}. The tame composition estimate for
$T_0=\mathcal H^{-1}(Y_0)$, on the uniformly bounded range of $Y_0$,
then gives
\[
 \|T_0\|_{W^{n,q}}
 \leq C_{n,q}(1+\|Y_0\|_{W^{n,q}})
 \leq C_{n,q}(1-\kappa)^{-\frac72n+2}.
\]
Finally $W^{s+2,q}(\Omega)\hookrightarrow C^{s+1,\alpha}(\overline\Omega)$
for $0<\alpha<1-3/q$, proving \eqref{ell}--\eqref{en}.
\end{proof}

Fix integers $N>5$ and $s_0\geq3N+5$, and fix $3<q<\infty$.

\begin{proposition}[Boundary matching]\label{mat}
Let $(I_0,T_0),...,(I_N,T_N)$ are internal solutions constructed in \eqref{3.5}, \eqref{3.4}, \eqref{3.10}, \eqref{3.9} with boundary conditions \eqref{Tb} and \eqref{ikm}. Let $(\mathscr{I}_0,\mathscr{T}_0),...,(\mathscr{I}_N,\mathscr{T}_N)$ are lower boundary solutions constructed in \eqref{Bd} and \eqref{Bdk}. Let $((\mathscr{I}_1^+,\mathscr{T}_1^+),...,(\mathscr{I}_N^+,\mathscr{T}_N^+)$ are upper boundary solutions constructed in \eqref{ukc}. Then the composite approximation defined in~\eqref{sc} satisfies

\begin{align}\label{eb}
 T^a|_{\partial\Omega}=T_b,\qquad I^a|_{\gamma_-}=I_b.
\end{align}

\end{proposition}
\begin{proof}
For $1\leq k\leq N$, the truncated lower boundary layers
vanish at $x_3=1$, and $\chi(0)=1$.
Using $g_k^+(0)=0$,
$T_k|_1=\mathfrak g_{k,\infty}^+$, and \eqref{uim}, we obtain
\begin{align*}
 (T_k+\mathscr T_k^+)|_1
 &=\mathfrak g_{k,\infty}^+-\mathfrak g_{k,\infty}^+=0,\nonumber\\
 (I_k+\mathscr I_k^+)|_1
 &=I_k|_1+(4T_b^3T_k-I_k)|_1
-4T_b^3\mathfrak g_{k,\infty}^+=0,
                      \qquad\beta_3<0.
\end{align*}
At $x_3=0$, the truncated upper boundary layers vanish.
For $1\leq k\leq N$, the lower boundary conditions
\eqref{lkb}, namely $g_k(0)=0$ and
$\phi_k(0,\beta)=4(T_0^0)^3T_k^0-I_k(0,\beta)$ give
\begin{align*}
 (T_k+\mathscr T_k)|_0
 &=\mathfrak g_{k,\infty}-\mathfrak g_{k,\infty}=0,\nonumber\\
 (I_k+\mathscr I_k)|_0
 &=I_k|_0+4(T_0^0)^3T_k^0-I_k|_0
 -4(T_0^0)^3\mathfrak g_{k,\infty}=0,
 \quad\beta_3>0.
\end{align*}
At order zero, the lower Milne data \eqref{Bd}, \eqref{Tb} and \eqref{3.17} give
$(T_0+\mathscr T_0)|_0=T_b|_0$ and
$(I_0+\mathscr I_0)|_0=I_b|_0$ on $\{\beta_3>0\}$.
At the upper boundary, \eqref{Tb} gives $T_0|_1=T_b|_1$ and $I_0|_1=T_b^4=I_b$ on
$\{\beta_3<0\}$. Summing the finite expansions proves
\eqref{eb}.
\end{proof}

Choose positive constants $\lambda_{k,s}$, independent of $\kappa$,
satisfying
\begin{align}\label{plh}
 1&>\lambda_{0,0}>\lambda_{0,1}>\cdots>\lambda_{0,s_0}>0,\nonumber\\
 0&<\lambda_{k,s+1}<\lambda_{k-1,s+2}<\lambda_{k,s}<\lambda_{k-1,s+1},
\end{align}
where $ 1\leq k\leq N$ and $0\leq s<s_0-3k.$ For $D^w=D^w_{x_\shortparallel}$ and $z=\eta$ or $\xi$, define
\begin{align}\label{pwn}
 \|(b,V,U)\|_{k,s;z}
 :=\sum_{|w|\leq s}\Bigl(&\|D^wb\|_{L^\infty(\mathbb T^2)}
 +\|e^{\lambda_{k,|w|}(1-\kappa)z}D^wV\|_{L^\infty_{z,x_\shortparallel}}
 \nonumber\\
 &+\|e^{\lambda_{k,|w|}(1-\kappa)z}D^wU\|_{L^\infty_{z,x_\shortparallel,\beta}}\Bigr),
\end{align}
where $b=b(x_\shortparallel)$, $V=V(z,x_\shortparallel)$, and
$U=U(z,x_\shortparallel,\beta)$.
Set $\mathfrak g_{0,\infty}=a$. For the profile families considered below, define
\begin{align}\label{pk}
 \mathfrak A_{k,s}(\kappa)
 &:=\|(\mathfrak g_{k,\infty},\mathcal T_k,\mathcal I_k)\|_{k,s;\eta},
 \qquad 0\leq k\leq N,\\
 \mathfrak A^+_{k,s}(\kappa)
 &:=\|(\mathfrak g_{k,\infty}^+,\mathcal T_k^+,\mathcal I_k^+)\|_{k,s;\xi},
 \qquad 1\leq k\leq N,
 \label{upn}\\
 \mathfrak E_{k,r}(\kappa)
 &:=\|T_k\|_{W^{r+2,q}(\Omega)}
 +\sup_{\beta\in\mathbb S^2}\|I_k(\cdot,\beta)\|_{W^{r+1,q}(\Omega)}.
 \label{pin}
\end{align}
Here $0\leq s\leq s_0-3k$ and $0\leq r\leq s_0-3k-2$.
We set $\mathscr I_0^+=\mathscr T_0^+=0$ and
$\mathfrak A^+_{0,s}=0$ for $0\leq s\leq s_0$.
In particular, the leading-layer norm is
\begin{align}\label{pz}
 \mathfrak A_{0,s}(\kappa)
 =\|(a,\mathcal T_0,\mathcal I_0)\|_{0,s;\eta},
 \qquad 0\leq s\leq s_0.
\end{align}

\begin{proposition}[Construction and high-order bounds]\label{prop}
Assume the boundary condition in~\eqref{ue},\eqref{bb} and the regularity condition~\eqref{pr}.
Let $(\mathbf I,\mathbf T)$ be a lower nonlinear Milne profile
furnished by Theorem~\ref{thm2}. Assume that, for every
$x_\shortparallel\in\mathbb T^2$, $\mathbf T(\cdot,x_\shortparallel)$
satisfies the spectral assumption~\eqref{asA} for each weight
\[
 \vartheta=\lambda_{k,r}(1-\kappa),\qquad
 0\leq k\leq N,\quad 0\leq r\leq s_0-3k,\quad k+r\geq1,
\]
with $M\leq M_*<1$ uniformly in $x_\shortparallel,k,r$.
Then there exist unique interior profiles $(I_k,T_k)$, lower boundary-layer
profiles $(\mathcal I_k,\mathcal T_k)$, and upper boundary-layer
profiles $(\mathcal I_k^+,\mathcal T_k^+)$ for $0\leq k\leq N$ constructed in \eqref{3.10}, \eqref{3.9}, \eqref{cen-1} and \eqref{cen}.
The following derivative bounds hold.

\medskip\noindent\emph{(i) Lower boundary layers.}
The lower-layer norms satisfy
\begin{align*}
 \mathfrak A_{0,0}(\kappa)&\leq C(1-\kappa)^{-1},\\
 \mathfrak A_{0,s}(\kappa)&\leq C_s(1-\kappa)^{-\frac72s+1},
 &&1\leq s\leq s_0,\\
 \mathfrak A_{1,s}(\kappa)&\leq C_s(1-\kappa)^{-\frac72s-9},
 &&0\leq s\leq s_0-3.
\end{align*}
For $2\leq k\leq N$,
\begin{align}\label{par}
 \mathfrak A_{k,s}(\kappa)
 &\leq C_{k,s}(1-\kappa)^{-(k+\frac72s+5\cdot4^k)},
 \qquad 0\leq s\leq s_0-3k.
\end{align}
The lower temperature profiles also satisfy
\begin{align}\label{nlb}
 (1-\kappa)|\mathcal T_0|+|\partial_\eta\mathcal T_0|
 &\leq C e^{-\lambda_{0,0}(1-\kappa)\eta},\nonumber\\
 |\mathcal T_k|+|\partial_\eta\mathcal T_k|
 &\leq C_k e^{-\lambda_{k,0}(1-\kappa)\eta}
 \begin{cases}
 (1-\kappa)^{-9},&k=1,\\
 (1-\kappa)^{-\left(k+\frac72k+5\cdot4^{k-1}+32\right)},
       &2\leq k\leq N.
 \end{cases}
\end{align}
The truncated lower layers satisfy
\begin{align}\label{pc}
 |D^w\mathscr T_k(\eta,x_\shortparallel)|
 +|D^w\mathscr I_k(\eta,x_\shortparallel,\beta)|
 &\leq\mathfrak A_{k,s}(\kappa)e^{-\lambda_{k,s}(1-\kappa)\eta},
 \quad
 \begin{gathered}
 0\leq k\leq N,\\
 |w|=s\leq s_0-3k.
 \end{gathered}
\end{align}

\medskip\noindent\emph{(ii) Upper boundary layers.}
The leading upper layer vanishes:
$\mathcal T_0^+=\mathcal I_0^+=0$. The higher layers satisfy
\begin{align}\label{uqb}
 \mathfrak A^+_{1,s}(\kappa)
 &\leq C_s^+(1-\kappa)^{-(\frac72s+9)},
 &&0\leq s\leq s_0-3,\nonumber\\
 \mathfrak A^+_{k,s}(\kappa)
 &\leq C_{k,s}^+(1-\kappa)^{-(k+\frac72s+\frac72k+5\cdot4^{k-1}+32)},
 &&
 \begin{gathered}
 2\leq k\leq N,\\
 0\leq s\leq s_0-3k.
 \end{gathered}
\end{align}

The normal derivatives of the upper temperature profiles satisfy
\begin{align}\label{unb}
 \|e^{\lambda_{k,0}(1-\kappa)\xi}\partial_\xi\mathcal T_k^+\|_\infty
 &\leq C_k
 \begin{cases}
 (1-\kappa)^{-15/2},& k=1,\\
 (1-\kappa)^{-\left(k+\frac72k+5\cdot4^{k-1}+32\right)},& 2\leq k\leq N.
 \end{cases}
\end{align}
The truncated upper layers satisfy
\begin{align}\label{upc}
 |D^w\mathscr T_k^+(\xi,x_\shortparallel)|
 +|D^w\mathscr I_k^+(\xi,x_\shortparallel,\beta)|
 &\leq\mathfrak A^+_{k,s}(\kappa)e^{-\lambda_{k,s}(1-\kappa)\xi},
 \quad
 \begin{gathered}
 0\leq k\leq N,\\
 |w|=s\leq s_0-3k.
 \end{gathered}
\end{align}

\medskip\noindent\emph{(iii) Interior profiles.}
The leading temperature $T_0$ satisfies
\eqref{ell}--\eqref{en} with $m=s_0-2$,
and $I_0=T_0^4$. For the first-order profiles,
\begin{align*}
 \mathfrak E_{1,r}(\kappa)
 &\leq C_{r,q}(1-\kappa)^{-\frac72r-18},
 \qquad 0\leq r\leq s_0-5.
\end{align*}
For $2\leq k\leq N$,
\begin{align}\label{per}
 \mathfrak E_{k,r}(\kappa)
 &\leq C_{k,r,q}(1-\kappa)^{-(k+\frac72r+5\cdot4^k)},
 \qquad 0\leq r\leq s_0-3k-2.
\end{align}
All constants in these estimates are independent of $\kappa$ and $\v$
under the stated uniform hypotheses.
\end{proposition}
\begin{remark}[Spectral weights and boundary data]\label{spw}
The spectral assumptions in Lemma~\ref{ffr} and Proposition~\ref{prop} can all be covered by
the boundary-mismatch condition~\eqref{4.62}.
Indeed, let $\lambda$ be any of the required parameters, so that
$0<\lambda<\lambda_{0,0}$. Applying Lemma~\ref{lem-sp} with
$\alpha=\lambda_{0,0}(1-\kappa)$ and
$\vartheta=\lambda(1-\kappa)$, it suffices to require
\[
 C_b<c_{a_1,a_2,M_*}(1-\kappa)
 \sqrt{\lambda_{0,0}\lambda
       (\lambda_{0,0}-\lambda)(1-\lambda_{0,0})}
\]
uniformly in $x_\shortparallel$.
There are only finitely many required weights, so choosing $C_b$
below the minimum of these thresholds ensures all the spectral
bounds with a common $M_*<1$.
For fixed decay parameters and $0\leq\kappa\leq\kappa_0<1$,
the smallness threshold can be chosen independently of $\kappa$.

At the upper boundary, the coefficient is $4T_b^3$, with $T_b$
evaluated at $x_3=1$. It is positive by~\eqref{bb} and independent
of $\xi$, so $\partial_\xi(2T_b^{3/2})=0$.
Thus the upper Milne problems satisfy the spectral inequality
automatically for every required weight and any fixed $M\in(0,1)$.
\end{remark}

\begin{proof}
We induct on the layer order $k$. At each order, both Milne problems
use only profiles of orders below $k$; their far-field values then
determine the two Dirichlet data for $T_k$, and~\eqref{3.9} determines
$I_k$. Within each Milne problem we induct separately on the
tangential derivative order $s$. 

\underline{Step 1. Zeroth-order profiles and tangential source estimates.}

\medskip\noindent\emph{Step 1.1. Zeroth-order boundary layers and interior profiles.}
Set $\mathcal T_0=\mathbf T-a$ and $\mathcal I_0=\mathbf I-a^4$.
By~\eqref{4.0} and Lemma~\ref{ffr}, the norm~\eqref{pz} satisfies
\begin{align}\label{pzr}
 \mathfrak A_{0,0}(\kappa)&\leq C_{0,0}(1-\kappa)^{-1},\nonumber\\
 \mathfrak A_{0,s}(\kappa)&\leq C_{0,s}(1-\kappa)^{-\frac72s+1},
 \qquad 1\leq s\leq s_0.
\end{align}
Since $\chi$ is independent of $x_\shortparallel$,
\eqref{Bd} gives~\eqref{pc} for $k=0$.
Lemma~\ref{lem1}, applied with $m=s_0-2$, constructs $T_0$
with the bounds~\eqref{ell}--\eqref{en}; set $I_0=T_0^4$.
The leading upper profiles are
$\mathcal T_0^+=\mathcal I_0^+=0$ by~\eqref{ue}.

\medskip\noindent\emph{Step 1.2. Tangential estimates for the linear Milne problems.}

\emph{Lower boundary.} For any $k\geq1$, consider the linear problem
\eqref{3.18} with prescribed sources and incoming data.
For $|w|=s\geq1$, set
\[
 \widehat U_{k,w}:=D^w\mathcal I_k+4a^3D^w\mathfrak g_{k,\infty},
 \qquad \widehat V_{k,w}:=D^wg_k.
\]
The additional coefficient term is
\begin{align}\label{lc}
 \mathcal Q_{k,w}:=
 \sum_{0<u\leq w}\binom wu\Bigl[
D^u(4\mathbf T^3)D^{w-u}\mathcal T_k+D^u\bigl(4\mathbf T^3-4a^3\bigr)
       D^{w-u}\mathfrak g_{k,\infty}\Bigr].
\end{align}
Since $\mathcal T_k=g_k-\mathfrak g_{k,\infty}$, $\mathcal Q_{k,w}$ contain
only tangential derivatives of the $k$-th temperature profile of order
strictly below $s$. Formally applying $D^w$ to both equations in~\eqref{3.18} and using
the definitions of $\widehat U_{k,w}$, $\widehat V_{k,w}$, and
$\mathcal Q_{k,w}$ in the resulting system, we obtain
\begin{align}\label{ld}
 \begin{cases}
 \beta_3\partial_\eta\widehat U_{k,w}
  +\kappa(\widehat U_{k,w}-\overline{\widehat U}_{k,w})
  =(1-\kappa)(4\mathbf T^3\widehat V_{k,w}-\widehat U_{k,w})
     +D^w\mathcal S_{k,1}+(1-\kappa)\mathcal Q_{k,w},\\
 \partial_\eta^2\widehat V_{k,w}
  +(1-\kappa)4\pi(\overline{\widehat U}_{k,w}-4\mathbf T^3\widehat V_{k,w})
  =-D^w\mathcal S_{k,2}+4\pi(1-\kappa)\mathcal Q_{k,w},
 \end{cases}
\end{align}
with $\widehat V_{k,w}(0)=0$ and, for $\beta_3>0$,
\begin{align}\label{ldb}
 \widehat U_{k,w}(0)
 =D^w\phi_k(0)
 -\sum_{0<u\leq w}\binom wu
       D^u(4a^3)D^{w-u}\mathfrak g_{k,\infty}.
\end{align}
Its far-field values are
$\widehat V_{k,w}(\infty)=D^w\mathfrak g_{k,\infty}$ and
$\widehat U_{k,w}(\infty)=4a^3D^w\mathfrak g_{k,\infty}$.


For $1\leq |u|\leq s$, estimates
\eqref{6.0-0} and~\eqref{czb} give
\begin{align}\label{4.53}
	\begin{split}
 \|D^u(4\mathbf T^3)\|_\infty
 +\|D^u(4a^3)\|_\infty
 &\leq C_s(1-\kappa)^{-\frac72|u|+2},\\
 \|e^{\lambda_{0,|u|}(1-\kappa)\eta}
       D^u(4\mathbf T^3-4a^3)\|_\infty
 &\leq C_s(1-\kappa)^{-\frac72|u|+1}.
 \end{split}
\end{align}
Using the two coefficient bounds above in~\eqref{lc},
we obtain
\begin{align}\label{lcs}
 &\|e^{\lambda_{k,s-1}(1-\kappa)\eta}\mathcal Q_{k,w}\|_\infty
 \nonumber\\
 &\quad\leq C_{k,s}\sum_{0<u\leq w}
 \Bigl[(1-\kappa)^{2-\frac72|u|}
 \|e^{\lambda_{k,s-1}(1-\kappa)\eta}D^{w-u}\mathcal T_k\|_\infty
 +(1-\kappa)^{1-\frac72|u|}
 \|D^{w-u}\mathfrak g_{k,\infty}\|_\infty\Bigr].
\end{align}
Here $\lambda_{k,s-1}<\lambda_{0,s}$.
Multiplying the estimate for \(\mathcal Q_{k,w}\) in~\eqref{lcs} by \(1-\kappa\) and applying \eqref{4.53} to the sum in~\eqref{ldb}, we obtain
\begin{align}\label{lcb}
 &\|e^{\lambda_{k,s-1}(1-\kappa)\eta}(1-\kappa)\mathcal Q_{k,w}\|_\infty
 +\left\|\sum_{0<u\leq w}\binom wu
       D^u(4a^3)D^{w-u}\mathfrak g_{k,\infty}\right\|_\infty
 \nonumber\\
 &\quad\leq C_{k,s}\sum_{0<u\leq w}(1-\kappa)^{2-\frac72|u|}
 \left(\|D^{w-u}\mathfrak g_{k,\infty}\|_\infty
 +\|e^{\lambda_{k,s-1}(1-\kappa)\eta}D^{w-u}\mathcal T_k\|_\infty\right).
\end{align}
In this estimate, the coefficient
$(1-\kappa)^{3-\frac72|u|}$ from the first term of
$(1-\kappa)\mathcal Q_{k,w}$ is bounded by
$(1-\kappa)^{2-\frac72|u|}$ since $0<1-\kappa\leq1$.
Integration of the exponential weight also gives
\begin{align}\label{lcl}
 \|e^{\lambda_{k,s}(1-\kappa)\eta}(1-\kappa)\mathcal Q_{k,w}\|_2\leq
 \frac{C}{\sqrt{(\lambda_{k,s-1}-\lambda_{k,s})(1-\kappa)}}
 \|e^{\lambda_{k,s-1}(1-\kappa)\eta}(1-\kappa)\mathcal Q_{k,w}\|_\infty .
\end{align}

\medskip\noindent\emph{Upper boundary.} Write $c^+=4T_b^3$. The data bounds give
\[
 4a_1^3\leq c^+\leq4a_2^3,\qquad
 \|D^u c^+\|_{L^\infty(\mathbb T^2)}\leq C_{\mathrm{data}}
 \quad(|u|\leq s\leq s_0).
\]
The coefficient is constant in $\xi$, so~\eqref{asA} holds with
$M=1/2$ at every admissible weight. 
For $|w|=s$, set
\[
 \widehat U_{k,w}^+=D^w\mathcal I_k^+
                    +c^+D^w\mathfrak g_{k,\infty}^+,
 \qquad \widehat V_{k,w}^+=D^wg_k^+.
\]
For $s\geq1$, define
\begin{align}\label{ucm}
 \mathcal Q_{k,w}^+
 =\sum_{0<u\leq w}\binom wu D^u c^+D^{w-u}\mathcal T_k^+.
\end{align}
For $s=0$, the sum is empty. Differentiating \eqref{uks} gives
\begin{align}\label{uds}
 \begin{cases}
 -\beta_3\partial_\xi\widehat U_{k,w}^+
       +\widehat U_{k,w}^+-\kappa\overline{\widehat U_{k,w}^+}
   =(1-\kappa) c^+\widehat V_{k,w}^+
            +D^w\mathcal S_{k,1}^++(1-\kappa)\mathcal Q_{k,w}^+,\\
 -\partial_\xi^2\widehat V_{k,w}^+
   =(1-\kappa)\langle\widehat U_{k,w}^+-c^+\widehat V_{k,w}^+\rangle
            +D^w\mathcal S_{k,2}^+-4\pi(1-\kappa)\mathcal Q_{k,w}^+,\\
 \widehat V_{k,w}^+(0)=0.
 \end{cases}
\end{align}
The incoming value is
\begin{align}\label{udb}
 \widehat U_{k,w}^+(0,\beta)
 =D^w\phi_k^+(0,\beta)
   -\sum_{0<u\leq w}\binom wu
           D^u c^+D^{w-u}\mathfrak g_{k,\infty}^+,
 \qquad\beta_3<0.
\end{align}
The far-field values are
$(\widehat V_{k,w}^+,\widehat U_{k,w}^+)\to
(D^w\mathfrak g_{k,\infty}^+,c^+D^w\mathfrak g_{k,\infty}^+)$.

Suppose that the upper estimates at this same $k$ have been proved
through tangential order $s-1$. The uniform bounds for $D^u c^+$
give
\begin{align}\label{ucn}
 \|e^{\lambda_{k,s-1}(1-\kappa)\xi}(1-\kappa)\mathcal Q_{k,w}^+\|_\infty
 &\leq C_s(1-\kappa)\sum_{0<u\leq w}\mathfrak A^+_{k,s-|u|},\nonumber\\
 \|e^{\lambda_{k,s}(1-\kappa)\xi}(1-\kappa)\mathcal Q_{k,w}^+\|_2
 &\leq C_{k,s}(1-\kappa)^{1/2}
                       \sum_{0<u\leq w}\mathfrak A^+_{k,s-|u|},\nonumber\\
 \left\|\sum_{0<u\leq w}\binom wu
          D^u c^+D^{w-u}\mathfrak g_{k,\infty}^+\right\|_\infty
 &\leq C_s\sum_{0<u\leq w}\mathfrak A^+_{k,s-|u|}.
\end{align}
The second line follows by integrating
$e^{-2(\lambda_{k,s-1}-\lambda_{k,s})(1-\kappa)\xi}$; this is the
one-half-power loss from weighted $L^\infty_\xi$ to $L^2_\xi$.
These bounds use only strictly lower tangential derivatives at
the current layer order. 

\underline{Step 2. First-order boundary layers and interior profiles.}

\medskip\noindent\emph{Step 2.1. Lower boundary layer.} By \eqref{3.19} and \eqref{3.15-1}, using $T_0^0=a$, we have
\begin{align}\label{6.5}
	\mathcal{S}_{1,1}&=(1-\kappa)\left[\mathcal{P}(\mathcal{T},1)-4\mathbf{T}^3(T_1^0+\mathcal{T}_1)+4a^3T_1^0\right]-\beta_{\shortparallel}\cdot\nabla_{\shortparallel}\mathcal{I}_0\nonumber\\
	&=(1-\kappa)(4\mathbf{T}^3-4a^3)\eta\partial_{x_3} T_0^0-\beta_{\shortparallel}\cdot\nabla_{\shortparallel}\mathcal{I}_0.
\end{align}
In the second equation of \eqref{3.18}, the corresponding source is
\begin{align}\label{pfs}
 \mathcal S_{1,2}=-4\pi(1-\kappa)(4\mathbf T^3-4a^3)\eta\partial_{x_3}T_0^0.
\end{align}
For $|w|=s\leq s_0-3$, Leibniz' rule and
$\sup_{\eta\geq0}\eta e^{-h\eta}\leq h^{-1}$ give
\begin{align}\label{pfsb}
 &\sum_{j=1}^2
 \|e^{\lambda_{0,s+1}(1-\kappa)\eta}D^w\mathcal S_{1,j}\|_\infty
 \nonumber\\
 &\quad\leq
 \frac{C_s}{\lambda_{0,s}-\lambda_{0,s+1}}
 \sum_{u\leq w}\binom wu
 \|e^{\lambda_{0,s}(1-\kappa)\eta}D^u(\mathbf T^3-a^3)\|_\infty
 \|D^{w-u}\partial_{x_3}T_0^0\|_\infty
 +C_s(1-\kappa)^{-\frac72s-\frac52}\nonumber\\
 &\quad\leq C_{s,q}(1-\kappa)^{-\frac72s-6}.
\end{align}
The exponent follows from Lemma \ref{ffr}, Lemma \ref{lem1} and {\it Step 1}.

The incoming datum is
$\phi_1(0,\beta)=\beta\cdot\nabla I_0(0)$ for $\beta_3>0$, and $g_1(0)=0$.
Since $I_0=T_0^4$, the composition estimate used in Lemma \ref{lem1}
and Sobolev embedding imply
\begin{align}\label{pfb}
 \max_{|w|\leq s}\sup_{x_{\shortparallel},\,\beta_3>0}
 |D^w(\beta\cdot\nabla I_0)(0)|
 \leq C_{s,q}\|I_0\|_{W^{s+2,q}(\Omega)}
 \leq C_{s,q}(1-\kappa)^{-\frac72s-5}.
\end{align}

Using~\eqref{pfsb} and integrating
$e^{-2(\lambda_{0,s+1}-\lambda_{1,s})(1-\kappa)\eta}$ over $\eta\in(0,\infty)$, we obtain
\begin{equation}\label{pfsb-l2}
 \sum_{j=1}^2
 \|e^{\lambda_{1,s}(1-\kappa)\eta}D^w\mathcal S_{1,j}\|_{L^2}
 \leq C_{s,q}(1-\kappa)^{-\frac72s-\frac{13}{2}}.
\end{equation}
Theorem~\ref{co} now gives a unique bounded solution $(\phi_1,g_1)$
and its far-field temperature $\mathfrak g_{1,\infty}$.
Set $\mathcal T_1=g_1-\mathfrak g_{1,\infty}$ and
$\mathcal I_1=\phi_1-4a^3\mathfrak g_{1,\infty}$.
At $s=0$, substituting the source bounds~\eqref{pfsb}, \eqref{pfsb-l2}
and the incoming bound~\eqref{pfb} into~\eqref{le} gives
\begin{align}\label{6.7}
 \|\mathfrak g_{1,\infty}\|_\infty
 +\|e^{\lambda_{1,0}(1-\kappa)\eta}\mathcal T_1\|_\infty
 &\leq C(1-\kappa)^{-8},\nonumber\\
 \|e^{\lambda_{1,0}(1-\kappa)\eta}\mathcal I_1\|_\infty
 &\leq C(1-\kappa)^{-9}.
\end{align}

Equation~\eqref{6.7} gives the base case $s=0$.
Fix $1\leq s\leq s_0-3$. For $0\leq r<s$, assume the bound for
$\mathfrak A_{1,r}$ in part~(i) and the sharper temperature estimate
\begin{equation}\label{fth}
 \|D^v\mathfrak g_{1,\infty}\|_\infty
 +\|e^{\lambda_{1,r}(1-\kappa)\eta}D^v\mathcal T_1\|_\infty
 \leq C_{r,q}(1-\kappa)^{-\frac72r-8},
 \qquad |v|=r.
\end{equation}
For $|w|=s$ and $0<u\leq w$, we have $|w-u|=s-|u|<s$.
Substituting~\eqref{fth} with $v=w-u$ into~\eqref{lcs}, we obtain
\begin{align}
 &\|e^{\lambda_{1,s-1}(1-\kappa)\eta}\mathcal Q_{1,w}\|_\infty\nonumber\\
 &\quad\leq C_s\sum_{0<u\leq w}\binom wu
 \left[(1-\kappa)^{2-\frac72|u|}
       +(1-\kappa)^{1-\frac72|u|}\right]
 (1-\kappa)^{-\frac72(s-|u|)-8}\nonumber\\
 &\quad\leq C_s\left[(1-\kappa)^{-\frac72s-6}
                    +(1-\kappa)^{-\frac72s-7}\right]
 \leq C_s(1-\kappa)^{-\frac72s-7}.\label{fqi}
\end{align}
Applying~\eqref{lcl} with $k=1$ and using~\eqref{fqi}, we obtain
\begin{align}
 &\|e^{\lambda_{1,s}(1-\kappa)\eta}
             (1-\kappa)\mathcal Q_{1,w}\|_2\nonumber\\
 &\quad\leq
 \frac{C(1-\kappa)^{1/2}}{\sqrt{\lambda_{1,s-1}-\lambda_{1,s}}}
 \|e^{\lambda_{1,s-1}(1-\kappa)\eta}\mathcal Q_{1,w}\|_\infty
 \leq C_s(1-\kappa)^{-\frac72s-\frac{13}{2}}.\label{fql}
\end{align}
For $k=1$, we estimate the full source pair in~\eqref{ld} in weighted $L^2$
by~\eqref{pfsb-l2} and~\eqref{fql}. For the transport source
in weighted $L^\infty$, we use~\eqref{pfsb} and multiply~\eqref{fqi}
by $1-\kappa$. The triangle inequality then gives
\begin{align}
 &\left\|e^{\lambda_{1,s}(1-\kappa)\eta}
 \bigl(D^w\mathcal S_{1,1}+(1-\kappa)\mathcal Q_{1,w},
       -D^w\mathcal S_{1,2}+4\pi(1-\kappa)\mathcal Q_{1,w}\bigr)
       \right\|_2\nonumber\\
 &\quad\leq
 \sum_{j=1}^2\|e^{\lambda_{1,s}(1-\kappa)\eta}
                         D^w\mathcal S_{1,j}\|_2
 +C\|e^{\lambda_{1,s}(1-\kappa)\eta}
                         (1-\kappa)\mathcal Q_{1,w}\|_2\nonumber\\
 &\quad\leq C_{s,q}(1-\kappa)^{-\frac72s-\frac{13}{2}},\label{ffl}\\
 &\|e^{\lambda_{1,s}(1-\kappa)\eta}
       (D^w\mathcal S_{1,1}+(1-\kappa)\mathcal Q_{1,w})\|_\infty\nonumber\\
 &\quad\leq
 \|e^{\lambda_{1,s}(1-\kappa)\eta}D^w\mathcal S_{1,1}\|_\infty
 +\|e^{\lambda_{1,s}(1-\kappa)\eta}
                         (1-\kappa)\mathcal Q_{1,w}\|_\infty\nonumber\\
 &\quad\leq C_{s,q}(1-\kappa)^{-\frac72s-6}.\label{ffi}
\end{align}
In~\eqref{ldb} with $k=1$, we bound $D^w\phi_1(0)$ by~\eqref{pfb}
and the boundary correction by~\eqref{lcb} and~\eqref{fth}.
This gives
\begin{align}
 \|\widehat U_{1,w}(0)\|_{L^\infty(\gamma_-)}
 &\leq \|D^w\phi_1(0)\|_{L^\infty(\gamma_-)}
 +\sum_{0<u\leq w}\binom wu
   \|D^u(4a^3)\|_\infty
   \|D^{w-u}\mathfrak g_{1,\infty}\|_\infty\nonumber\\
 &\leq C_{s,q}(1-\kappa)^{-\frac72s-5}
 +C_s\sum_{0<u\leq w}\binom wu
 (1-\kappa)^{2-\frac72|u|}
 (1-\kappa)^{-\frac72(s-|u|)-8}\nonumber\\
 &\leq C_{s,q}(1-\kappa)^{-\frac72s-6}.\label{fbi}
\end{align}
We apply~\eqref{le} to~\eqref{ld}--\eqref{ldb} with $k=1$,
$\tilde b=\lambda_{1,s}$ and $\widehat V_{1,w}(0)=0$.
We use~\eqref{ffl} for the source pair,
\eqref{ffi} for the transport source, and~\eqref{fbi}
for the incoming datum, since
$|\widehat U_{1,w}(0)|_{L^2(\gamma_-)}
\leq C\|\widehat U_{1,w}(0)\|_{L^\infty(\gamma_-)}$.
Summing the resulting estimates over $|w|=s$ and including the
already proved lower derivative orders, we obtain
\begin{align}
 &\sum_{|w|\leq s}\left(\|D^w\mathfrak g_{1,\infty}\|_\infty
 +\|e^{\lambda_{1,|w|}(1-\kappa)\eta}D^w\mathcal T_1\|_\infty\right)
 \nonumber\\
 &\quad\leq \frac{C_{s,q}}{(1-\kappa)^{3/2}}
 \left[(1-\kappa)^{-\frac72s-6}
       +(1-\kappa)^{-\frac72s-\frac{13}{2}}\right]
 \leq C_{s,q}(1-\kappa)^{-\frac72s-8},
 \label{ltr}\\
 &\sum_{|w|\leq s}
 \|e^{\lambda_{1,|w|}(1-\kappa)\eta}D^w\mathcal I_1\|_\infty
 \nonumber\\
 &\quad\leq \frac{C_{s,q}}{(1-\kappa)^{5/2}}
 \left[(1-\kappa)^{-\frac72s-6}
       +(1-\kappa)^{-\frac72s-\frac{13}{2}}\right]
 +\frac{C_{s,q}}{1-\kappa}(1-\kappa)^{-\frac72s-6}
 \nonumber\\
 &\quad\leq C_{s,q}(1-\kappa)^{-\frac72s-9}.
 \label{6.8}
\end{align}
Together with the case $s=0$ in~\eqref{6.7}, estimates~\eqref{ltr}
and~\eqref{6.8} give the bound for $\mathfrak A_{1,s}$ in part~(i).
Since the cutoff in~\eqref{Bdk} depends only on $\eta$ and
$|\chi|\leq1$, the definition~\eqref{pk} also gives~\eqref{pc} for $k=1$.\\
\medskip\noindent\emph{Step 2.2. Upper boundary layer.}
At $k=1$, the sources and boundary data are
\[
 \mathcal S_{1,1}^+=\mathcal S_{1,2}^+=0,\qquad
 g_1^+(0)=0,\qquad
 \phi_1^+(0,\beta)=\beta\cdot\nabla I_0|_{x_3=1},
 \quad \beta_3<0.
\]
The incoming norms below are taken over $\beta_3<0$.
The trace estimate for $I_0=T_0^4$ gives
\[
 \max_{|w|\leq s}\|D^w\phi_1^+(0)\|_\infty
 \leq C_{s,q}\|I_0\|_{W^{s+2,q}(\Omega)}
 \leq C_{s,q}(1-\kappa)^{-\frac72s-5},
 \qquad 0\leq s\leq s_0-3.
\]
The last inequality follows from Lemma~\ref{lem1} and the tame
composition estimate for $T_0^4$, using the uniform bound for $T_0$.
After the change of angular variable \(\beta_3\mapsto-\beta_3\), Theorem~\ref{co} gives the unique
bounded solution $(\phi_1^+,g_1^+)$ and its far-field temperature
$\mathfrak g_{1,\infty}^+$. Set
$\mathcal T_1^+=g_1^+-\mathfrak g_{1,\infty}^+$ and
$\mathcal I_1^+=\phi_1^+-c^+\mathfrak g_{1,\infty}^+$.
At $s=0$, substituting the incoming bound into~\eqref{le} with zero sources gives
\begin{align*}
 \|\mathfrak g_{1,\infty}^+\|_\infty
 +\|e^{\lambda_{1,0}(1-\kappa)\xi}\mathcal T_1^+\|_\infty
 &\leq \frac{C}{(1-\kappa)^{3/2}}\|\phi_1^+(0)\|_\infty
 \leq C(1-\kappa)^{-13/2},\\
 \|e^{\lambda_{1,0}(1-\kappa)\xi}\mathcal I_1^+\|_\infty
 &\leq \frac{C}{(1-\kappa)^{5/2}}\|\phi_1^+(0)\|_\infty
 \leq C(1-\kappa)^{-15/2}.
\end{align*}
By~\eqref{upn}, these bounds give $\mathfrak A^+_{1,0}(\kappa)\leq C(1-\kappa)^{-9}$.

Fix $1\leq s\leq s_0-3$ and assume the first bound in~\eqref{uqb} at all lower orders,
\[
 \mathfrak A^+_{1,r}(\kappa)\leq C_r(1-\kappa)^{-\frac72r-9},
 \qquad 0\leq r<s.
\]
For $|w|=s$ and $0<u\leq w$, we have $|w-u|=s-|u|\leq s-1$.
Substituting the induction hypothesis with $r=s-|u|$ into the first
two inequalities of~\eqref{ucn} gives
\begin{align*}
 \|e^{\lambda_{1,s-1}(1-\kappa)\xi}(1-\kappa)\mathcal Q_{1,w}^+\|_\infty
 &\leq C_s(1-\kappa)\sum_{0<u\leq w}\mathfrak A^+_{1,s-|u|}
 \leq C_s(1-\kappa)^{-\frac72s-\frac92},\\
 \|e^{\lambda_{1,s}(1-\kappa)\xi}(1-\kappa)\mathcal Q_{1,w}^+\|_2
 &\leq C_s(1-\kappa)^{1/2}\sum_{0<u\leq w}\mathfrak A^+_{1,s-|u|}
 \leq C_s(1-\kappa)^{-\frac72s-5}.
\end{align*}
Since $\mathcal S_{1,1}^+=\mathcal S_{1,2}^+=0$, the source pair
in~\eqref{uds}, written in the form~\eqref{5.1}, is
$((1-\kappa)\mathcal Q_{1,w}^+,4\pi(1-\kappa)\mathcal Q_{1,w}^+)$.
In~\eqref{udb}, we use the trace bound for $D^w\phi_1^+(0)$ and
the last inequality of~\eqref{ucn} to obtain
\begin{align*}
 \|\widehat U_{1,w}^+(0)\|_\infty
 \leq C_{s,q}(1-\kappa)^{-\frac72s-5}
       +C_s\sum_{0<u\leq w}\mathfrak A^+_{1,s-|u|}\leq C_{s,q}(1-\kappa)^{-\frac72s-\frac{11}{2}}.
\end{align*}
After reflecting $\beta_3$, we insert these source and boundary bounds
into~\eqref{le} with $\tilde b=\lambda_{1,s}$ and $\widehat V_{1,w}^+(0)=0$.
Summing over $|w|=s$ and including the lower orders by the induction hypothesis gives
\begin{align}
 &\sum_{|w|\leq s}\left(\|D^w\mathfrak g_{1,\infty}^+\|_\infty
 +\|e^{\lambda_{1,|w|}(1-\kappa)\xi}D^w\mathcal T_1^+\|_\infty\right)
 \nonumber\\
 &\quad\leq \frac{C_{s,q}}{(1-\kappa)^{3/2}}
 (1-\kappa)^{-\frac72s-\frac{11}{2}}
 \leq C_{s,q}(1-\kappa)^{-\frac72s-8},\nonumber\\
 &\sum_{|w|\leq s}\|e^{\lambda_{1,|w|}(1-\kappa)\xi}D^w\mathcal I_1^+\|_\infty
 \nonumber\\
 &\quad\leq \frac{C_{s,q}}{(1-\kappa)^{5/2}}
 (1-\kappa)^{-\frac72s-\frac{11}{2}}
 +\frac{C_{s,q}}{1-\kappa}(1-\kappa)^{-\frac72s-\frac92}
 \leq C_{s,q}(1-\kappa)^{-\frac72s-9}.\label{ufb}
\end{align}
Adding the two estimates in~\eqref{ufb} gives the first bound
in~\eqref{uqb} at order $s$, which closes the induction.

\medskip\noindent\emph{Step 2.3. Interior profiles.}
The two first-order far-field values now give
$T_1|_{x_3=0}=\mathfrak g_{1,\infty}$ and
$T_1|_{x_3=1}=\mathfrak g_{1,\infty}^+$.
Put $c_0=\mathcal H'(T_0)=1+\frac{16\pi}{3}T_0^3$.
For $k=1$, $\mathcal C(T,1)=4T_0^3T_1$, while
$\int_{\mathbb S^2}(\beta\cdot\nabla)^3I_0\,d\beta=0$ by oddness.
Thus identity~\eqref{3.10} gives
\begin{align}\label{fi}
 \Delta(c_0T_1)=0,\qquad
 (c_0T_1)|_{x_3=0}=\mathcal H'(a)\mathfrak g_{1,\infty},
 \qquad (c_0T_1)|_{x_3=1}=\mathcal H'(T_b)\mathfrak g_{1,\infty}^+.
\end{align}
This Dirichlet problem determines $T_1$ uniquely since $c_0\geq1$;
then $I_1=4T_0^3T_1-\beta\cdot\nabla I_0$.
For $q>3$ and integers $0\leq r\leq s_0-5$,
\begin{align}\label{fie}
 \|c_0T_1\|_{W^{r+2,q}}
 &\leq C_{r,q}\bigl(
   \|\mathcal H'(a)\mathfrak g_{1,\infty}\|_{W^{r+2,\infty}}
  +\|\mathcal H'(T_b)\mathfrak g_{1,\infty}^+\|_{W^{r+2,\infty}}\bigr)
 \nonumber\\
 &\leq C_{r,q}\left[\mathfrak A_{1,r+2}(\kappa)
   +\sum_{\ell=1}^{r+2}(1-\kappa)^{-\frac72\ell+2}
          \mathfrak A_{1,r+2-\ell}(\kappa)
   +\mathfrak A^+_{1,r+2}(\kappa)\right]\nonumber\\
 &\leq C_{r,q}(1-\kappa)^{-\frac72r-16}.
\end{align}
The maximum principle also gives
\[
 \|c_0T_1\|_{L^q}
 \leq C_q\bigl(\|\mathcal H'(a)\mathfrak g_{1,\infty}\|_\infty
       +\|\mathcal H'(T_b)\mathfrak g_{1,\infty}^+\|_\infty\bigr)
 \leq C_q(1-\kappa)^{-9}.
\]
Since $c_0\geq1$, Lemma \ref{lem1} and the tame composition estimate
give
\[
 \|c_0^{-1}\|_{L^\infty}\leq1,
 \qquad
 \|c_0^{-1}\|_{W^{r+2,\infty}}
 \leq C_{r,q}(1-\kappa)^{-\frac72(r+3)+2}.
\]
The tame product estimate, with the high derivative placed on only one
factor in each term, therefore yields
\begin{align*}
 \|T_1\|_{W^{r+2,q}}
 \leq C_{r,q}\left(
 \|c_0T_1\|_{W^{r+2,q}}
 +\|c_0^{-1}\|_{W^{r+2,\infty}}\|c_0T_1\|_{L^q}
 \right)\leq  C_{r,q}(1-\kappa)^{-\frac72r-\frac{35}{2}}.
\end{align*}
Moreover, Lemma \ref{lem1} and the composition estimate give
\begin{align*}
 \|4T_0^3\|_{W^{r+1,\infty}}+\|I_0\|_{W^{r+2,q}}
 \leq C_{r,q}(1-\kappa)^{-\frac72(r+2)+2}.
\end{align*}
Taking the $W^{r+1,q}(\Omega)$ norm in the formula for $I_1$ in~\eqref{3.4},
applying the tame product estimate to $4T_0^3T_1$, and using
$\|T_1\|_{L^q}\leq\|c_0T_1\|_{L^q}$, we obtain
\begin{align}\label{fir}
 \sup_{\beta\in\mathbb S^2}\|I_1(\cdot,\beta)\|_{W^{r+1,q}(\Omega)}
 &\leq C_{r,q}\left[
 \|T_1\|_{W^{r+1,q}}
 +\|4T_0^3\|_{W^{r+1,\infty}}\|T_1\|_{L^q}
 +\|I_0\|_{W^{r+2,q}}\right]\nonumber\\
 &\leq C_{r,q}(1-\kappa)^{-\frac72r-\frac{35}{2}}
 \leq C_{r,q}(1-\kappa)^{-\frac72r-18}.
\end{align}
Since $\frac{35}{2}\leq18$, these estimates imply the first-layer
interior bound in the statement.

\underline{Step 3. Construction of the $k$-th layer for $k\geq 2$.}
Fix $2\leq k\leq N$. Assume that, for each $0\leq i<k$, the interior
profile and both boundary layers of order $i$ have been constructed and satisfy
\begin{align}\label{6.11}
 \mathfrak A_{i,s}(\kappa)+\mathfrak A^+_{i,s}(\kappa)
 &\le C_{i,s}(1-\kappa)^{-(i+\f72s+5\cdot 4^i)},
 &&0\le i\le k-1,\quad 0\le s\le s_0-3i,\nonumber\\
 \mathfrak E_{i,r}(\kappa)
 &\le C_{i,r,q}(1-\kappa)^{-(i+\f72r+5\cdot 4^i)},
 &&0\le i\le k-1,\quad 0\le r\le s_0-3i-2.
\end{align}
We construct both $k$-th boundary layers before solving
for the $k$-th interior profile, using the decay parameters fixed in~\eqref{plh}.

\medskip\noindent{\it Step 3.1. Lower boundary layer.}
Fix $0\leq s\leq s_0-3k$ and $|w|=s$, and set
$\vartheta=\lambda_{k,s}(1-\kappa)$.
We use~\eqref{le} with $\tilde b=\lambda_{k,s}$,
applied to~\eqref{3.18} for $s=0$ and
to~\eqref{ld}--\eqref{ldb} for $s\geq1$.
In the notation of~\eqref{5.1}, the temperature datum
is $g_b=0$, while the sources and incoming datum are
\begin{align*}
	S_1&=D^w\mathcal S_{k,1}
	+(1-\kappa)\mathcal Q_{k,w},\\
	S_2&=-D^w\mathcal S_{k,2}
	+4\pi(1-\kappa)\mathcal Q_{k,w},\\
	\phi_b&=D^w\phi_k(0,\beta)
	-\sum_{0<u\leq w}\binom wu
	D^u(4a^3)D^{w-u}\mathfrak g_{k,\infty},
	\qquad \beta_3>0.
\end{align*}
For $s=0$, the sums defining $\mathcal Q_{k,w}$
and the boundary correction are empty.
We also use
$|\phi_b|_{L^2(\gamma_-)}
\leq C\|\phi_b\|_{L^\infty(\gamma_-)}$.

Applying $D^w$ to~\eqref{3.19} gives
\begin{align}
 D^w\mathcal S_{k,1}
 &=(1-\kappa)D^w\bigl[\mathcal P(\mathcal T,k)
       -4\mathbf T^3(\mathcal T_k+T_k^0)+4a^3T_k^0\bigr]
-\beta_\shortparallel\cdot\nabla_\shortparallel D^w\mathcal I_{k-1},
 \nonumber\\
 D^w\mathcal S_{k,2}
 &=-4\pi(1-\kappa)D^w\bigl[\mathcal P(\mathcal T,k)
       -4\mathbf T^3(\mathcal T_k+T_k^0)+4a^3T_k^0\bigr]
+\Delta_\shortparallel D^w\mathcal T_{k-2}.
 \label{pds}
\end{align}
Substituting~\eqref{3.9} into~\eqref{lkb} cancels the term $4a^3T_k^0$.
Applying $D^w$ to the remaining terms gives
\begin{align}\label{pmd}
 D^w\phi_k(0,\beta)
 &=D^w\left[\beta\cdot\nabla I_{k-1}(0)
 +\frac{\kappa}{4\pi(1-\kappa)}\Delta T_{k-2}(0)\right.
 \nonumber\\
 &\qquad\left.{}-\sum_{\substack{i_1+\cdots+i_4=k\\0\leq i_j<k}}
       T_{i_1}^0T_{i_2}^0T_{i_3}^0T_{i_4}^0\right],
 \qquad\beta_3>0 .
\end{align}

To estimate the quartic term in~\eqref{pds}, expand $\mathcal P(\mathcal T,k)$
by~\eqref{3.15-0}. The products containing $\mathcal T_k$ or $T_k^0$ sum to
$4\mathbf T^3\mathcal T_k+4(\mathbf T^3-a^3)T_k^0$ and cancel in~\eqref{pds}. Substituting \eqref{tr} into the quartic remainder in \eqref{pds} and applying \(D^w\), we obtain a sum of terms of the form
\[
 F=\prod_{a=1}^{p}
 \left(\frac{\eta^{\ell_a}}{\ell_a!}
       D^{w_a}\partial_{x_3}^{\ell_a}T_{j_a}^0\right)
 \prod_{b=p+1}^{4}D^{w_b}\mathcal T_{j_b},
\]
up to numerical coefficients. Here $0\leq p\leq3$ and $\ell_a\geq0$
for $1\leq a\leq p$. For $1\leq i\leq4$, we have $0\leq j_i<k$ and
$r_i=|w_i|$, with
\[
 \sum_{i=1}^4 j_i+\sum_{a=1}^p\ell_a=k,\qquad
 \sum_{i=1}^4 w_i=w.
\]
Each product thus contains a decaying layer factor.
For the Taylor factors with $\ell\geq1$ and $|u|=r$, the embedding
$W^{r+\ell+1,q}(\Omega)\hookrightarrow W^{r+\ell,\infty}(\Omega)$
and~\eqref{6.11} at order $r+\ell-1$ give
\begin{align}\label{itn}
 \left\|D_{x_\shortparallel}^u\partial_{x_3}^\ell T_j\right\|_{L^\infty(\Omega)}
 &\leq C_{j,r,\ell}(1-\kappa)^{-\left(j+\frac72(r+\ell-1)+5\cdot4^j\right)},
 \qquad 0\leq j<k,
\end{align}
whenever $r+\ell-1\leq s_0-3j-2$.
By~\eqref{6.11} and~\eqref{itn},
\begin{align}\label{6.12}
 |F|
 &\leq C_{k,s}\eta^{\sum_{a=1}^{p}\ell_a}
 \prod_{a=1}^{p}
 \bigl[1+\mathfrak E_{j_a,r_a+(\ell_a-1)_+}(\kappa)\bigr]
 \prod_{b=p+1}^{4}\mathfrak A_{j_b,r_b}(\kappa)
 e^{-\lambda_{k-1,s}(1-\kappa)\eta},
\end{align}
where $(z)_+=\max\{z,0\}$.
The interior derivative orders are available: for $\ell_a=0$,
$r_a\leq s\leq s_0-3j_a-2$; for $\ell_a\geq1$, the relation
$j_a+\ell_a\leq k$ gives
\[
 r_a+\ell_a-1
 \leq s_0-3k+\ell_a-1
 \leq s_0-2k-j_a-1
 \leq s_0-3j_a-2.
\]
For each layer factor, $r_b\leq s\leq s_0-3j_b$ for $p+1\leq b\leq4$.

The gap between the two weights controls the polynomial factor:
\[
 \sup_{\eta\geq0}
 \eta^{\sum_{a=1}^{p}\ell_a}
 e^{-(\lambda_{k-1,s}-\lambda_{k-1,s+1})(1-\kappa)\eta}
 \leq C_{k,s}(1-\kappa)^{-\sum_{a=1}^{p}\ell_a}.
\]
Here the constant may depend on the fixed gap between the decay rates,
but is independent of $\kappa$. Using~\eqref{6.11} for each factor
in~\eqref{6.12}, we therefore obtain
\begin{align}\label{qwt}
 &(1-\kappa)
 \|e^{\lambda_{k-1,s+1}(1-\kappa)\eta}F\|_\infty
 \nonumber\\
 &\quad\leq C_{k,s}(1-\kappa)^{1-\sum_{a=1}^{p}\ell_a}
 \prod_{a=1}^{p}
 \bigl[1+\mathfrak E_{j_a,r_a+(\ell_a-1)_+}(\kappa)\bigr]
 \prod_{b=p+1}^{4}\mathfrak A_{j_b,r_b}(\kappa)
 \nonumber\\
 &\quad\leq C_{k,s}
 (1-\kappa)^{-\left(k-1+\frac72s+
       \frac72\sum_{a=1}^{p}(\ell_a-1)_+
       +5\sum_{i=1}^{4}4^{j_i}\right)}
 \nonumber\\
 &\quad\leq C_{k,s}
 (1-\kappa)^{-\left(k+\frac72s+\frac72k+5\cdot4^{k-1}+29\right)}.
\end{align}
In the second inequality we used
$\sum_{i=1}^{4}j_i+\sum_{a=1}^{p}\ell_a=k$ and
$\sum_{i=1}^{4}r_i=s$. For the last inequality,
$\sum_{a=1}^{p}(\ell_a-1)_+\leq k$, while
\[
 \sum_{i=1}^{4}4^{j_i}\leq4^{k-1}+6.
\]
Summing~\eqref{qwt} over the finitely many products and derivative
splits bounds the quartic term in~\eqref{pds}. By~\eqref{pk}, the
weighted norms of $\beta_\shortparallel\cdot\nabla_\shortparallel D^w\mathcal I_{k-1}$
and $\Delta_\shortparallel D^w\mathcal T_{k-2}$ are bounded by
$C_s\mathfrak A_{k-1,s+1}$ and $C_s\mathfrak A_{k-2,s+2}$, respectively.
Substituting these three bounds into~\eqref{pds} and using~\eqref{6.11} gives
\begin{align}\label{6.13}
 &\|e^{\lambda_{k-1,s+1}(1-\kappa)\eta}
       D^w\mathcal S_{k,j}\|_\infty
 \nonumber\\
 &\quad\leq C_{k,s}\Bigl[
 (1-\kappa)^{-\left(k+\frac72s+\frac72k+5\cdot4^{k-1}+29\right)}
 +\mathfrak A_{k-1,s+1}(\kappa)
 +\mathfrak A_{k-2,s+2}(\kappa)\Bigr]
 \nonumber\\
 &\quad\leq C_{k,s}\Bigl[
 (1-\kappa)^{-\left(k+\frac72s+\frac72k+5\cdot4^{k-1}+29\right)}
 +(1-\kappa)^{-\left(k+\frac72s+\frac52+5\cdot4^{k-1}\right)}
 +(1-\kappa)^{-\left(k+\frac72s+5+5\cdot4^{k-2}\right)}
 \Bigr]
 \nonumber\\
 &\quad\leq C_{k,s}
 (1-\kappa)^{-\left(k+\frac72s+\frac72k+5\cdot4^{k-1}+29\right)},
 \qquad j=1,2.
\end{align}
Both lower-order layer bounds are available: $s\leq s_0-3k$ implies
$s+1\leq s_0-3(k-1)$ and $s+2\leq s_0-3(k-2)$.
Moreover,
because $\lambda_{k,s}<\lambda_{k-1,s+1}$,
\begin{align*}
 \|e^{\lambda_{k,s}(1-\kappa)\eta}D^w\mathcal S_{k,j}\|_{L^2_\eta}
 &\leq
 \frac{\|e^{\lambda_{k-1,s+1}(1-\kappa)\eta}D^w\mathcal S_{k,j}\|_{L^\infty_\eta}}
 {\sqrt{2(\lambda_{k-1,s+1}-\lambda_{k,s})(1-\kappa)}}.
\end{align*}
Together with~\eqref{6.13}, this gives
\begin{align}\label{swt}
 \|e^{\lambda_{k-1,s+1}(1-\kappa)\eta}
 D^w\mathcal S_{k,j}\|_{L^\infty}
 &\leq C_{k,s}(1-\kappa)^{-\left(k+\frac72s+\frac72k
 +5\cdot4^{k-1}+29\right)},\nonumber\\
 \|e^{\lambda_{k,s}(1-\kappa)\eta}
 D^w\mathcal S_{k,j}\|_{L^2}
 &\leq C_{k,s}(1-\kappa)^{-\left(k+\frac72s+\frac72k
 +5\cdot4^{k-1}+\frac{59}{2}\right)}.
\end{align}
Applying Leibniz' rule to \eqref{pmd} gives the
boundary contribution
\begin{align}\label{pbr}
 &\max_{|w|\leq s}\|D^w\phi_k(0)\|_{L^\infty(\gamma_-)}
 \nonumber\\
 &\quad\leq C_{k,s}\Bigl[
 \mathfrak E_{k-1,s+1}(\kappa)
 +\frac{\kappa}{1-\kappa}\mathfrak E_{k-2,s+2}(\kappa)\nonumber\\
 &\qquad+
 \sum_{\substack{i_1+\cdots+i_4=k\\0\le i_j<k}}
 \sum_{\substack{r_1+\cdots+r_4\le s\\r_j\ge0}}
 \prod_{j=1}^4\mathfrak E_{i_j,r_j}(\kappa)\Bigr]\nonumber\\
 &\quad\leq C_{k,s}\Bigl[
 (1-\kappa)^{-(k-1)-\f72(s+1)-5\cdot4^{k-1}}\nonumber\\
 &\qquad +(1-\kappa)^{-(k-2)-\f72(s+2)-5\cdot4^{k-2}-1}
 +(1-\kappa)^{-\left(k+\f72s+5\cdot4^{k-1}+30\right)}
 \Bigr]\nonumber\\
 &\quad\leq C_{k,s}(1-\kappa)^{-\left(k+\f72s+5\cdot4^{k-1}+30\right)}.
\end{align}
In the quartic term we used
$\sum_{j=1}^4 i_j=k$, $i_j\leq k-1$, and hence
$\sum_{j=1}^4 4^{i_j}\leq4^{k-1}+6$; this is the origin of the constant
$30$ in the last exponent.
The interior induction bounds in \eqref{6.11} also cover all boundary
terms: for $s\leq s_0-3k$,
\[
 s+1\leq s_0-3(k-1)-2,
 \qquad
 s+2\leq s_0-3(k-2)-2,
\]
and $r_j\leq s\leq s_0-3i_j-2$ whenever $i_j\leq k-1$.

Combining \eqref{swt} and
\eqref{pbr} gives, for $j=1,2$,
\begin{align}\label{6.14}
 &\|e^{\lambda_{k-1,s+1}(1-\kappa)\eta}D^w\mathcal S_{k,j}\|_{L^\infty}
 +\|e^{\lambda_{k,s}(1-\kappa)\eta}D^w\mathcal S_{k,j}\|_{L^2}+\max_{|w|\le s}\|D^w\phi_k(0)\|_{L^\infty(\gamma_-)}
 \nonumber\\
 &\le C_{k,s}(1-\kappa)^{-\left(k+\frac72s+\frac72k
 +5\cdot4^{k-1}+\frac{59}{2}\right)}.
\end{align}
These source and boundary estimates verify the hypotheses of
Theorem~\ref{co}. It gives a unique bounded solution $(\phi_k,g_k)$
and a far-field temperature $\mathfrak g_{k,\infty}$. We set
$\mathcal T_k=g_k-\mathfrak g_{k,\infty}$ and
$\mathcal I_k=\phi_k-4a^3\mathfrak g_{k,\infty}$.

At $s=0$, substituting~\eqref{6.14} and the $L^\infty$ estimate
in~\eqref{swt} into~\eqref{le} gives
\begin{align*}
 &\|\mathfrak g_{k,\infty}\|_\infty
 +\|e^{\lambda_{k,0}(1-\kappa)\eta}\mathcal T_k\|_\infty\\
 &\quad\leq\frac{C_k}{(1-\kappa)^{3/2}}
 (1-\kappa)^{-\left(k+\frac72k+5\cdot4^{k-1}+\frac{59}{2}\right)}
 \leq C_k(1-\kappa)^{-\left(k+\frac72k+5\cdot4^{k-1}+31\right)},\\
 &\|e^{\lambda_{k,0}(1-\kappa)\eta}\mathcal I_k\|_\infty\\
 &\quad\leq\frac{C_k}{(1-\kappa)^{5/2}}
 (1-\kappa)^{-\left(k+\frac72k+5\cdot4^{k-1}+\frac{59}{2}\right)}
 +\frac{C_k}{1-\kappa}
 (1-\kappa)^{-\left(k+\frac72k+5\cdot4^{k-1}+29\right)}\\
 &\quad\leq C_k(1-\kappa)^{-\left(k+\frac72k+5\cdot4^{k-1}+32\right)}.
\end{align*}

Fix $1\leq s\leq s_0-3k$ and assume that, for every $|v|=r<s$,
\begin{align*}
 \|D^v\mathfrak g_{k,\infty}\|_\infty
 +\|e^{\lambda_{k,r}(1-\kappa)\eta}D^v\mathcal T_k\|_\infty
 &\leq C_{k,r}(1-\kappa)^{-\left(k+\frac72r+\frac72k+5\cdot4^{k-1}+31\right)},\\
 \|e^{\lambda_{k,r}(1-\kappa)\eta}D^v\mathcal I_k\|_\infty
 &\leq C_{k,r}(1-\kappa)^{-\left(k+\frac72r+\frac72k+5\cdot4^{k-1}+32\right)}.
\end{align*}
For $|w|=s$ and $0<u\leq w$, we have $|w-u|=s-|u|<s$.
Substituting the temperature bound at $v=w-u$ into~\eqref{lcb} gives
\begin{align*}
 &\|e^{\lambda_{k,s-1}(1-\kappa)\eta}
                      (1-\kappa)\mathcal Q_{k,w}\|_\infty\\
 &\quad\leq C_{k,s}\sum_{0<u\leq w}(1-\kappa)^{2-\frac72|u|}
 \left(\|D^{w-u}\mathfrak g_{k,\infty}\|_\infty
 +\|e^{\lambda_{k,s-1}(1-\kappa)\eta}
                         D^{w-u}\mathcal T_k\|_\infty\right)\\
 &\quad\leq C_{k,s}\sum_{0<u\leq w}
 (1-\kappa)^{2-\frac72|u|}
 (1-\kappa)^{-\left(k+\frac72(s-|u|)+\frac72k
                              +5\cdot4^{k-1}+31\right)}\\
 &\quad\leq C_{k,s}
 (1-\kappa)^{-\left(k+\frac72s+\frac72k+5\cdot4^{k-1}+29\right)}.
\end{align*}
Since $\lambda_{k,s}<\lambda_{k,s-1}$,
\eqref{lcl} gives
\begin{align*}
 &\|e^{\lambda_{k,s}(1-\kappa)\eta}
                      (1-\kappa)\mathcal Q_{k,w}\|_2\\
 &\quad\leq
 \frac{C}{\sqrt{(\lambda_{k,s-1}-\lambda_{k,s})(1-\kappa)}}
 \|e^{\lambda_{k,s-1}(1-\kappa)\eta}
                      (1-\kappa)\mathcal Q_{k,w}\|_\infty\\
 &\quad\leq C_{k,s}
 (1-\kappa)^{-\left(k+\frac72s+\frac72k+5\cdot4^{k-1}+\frac{59}{2}\right)}.
\end{align*}
Combining these inequalities with~\eqref{swt},
we bound the full sources in~\eqref{ld} by
\begin{align*}
 &\left\|e^{\lambda_{k,s}(1-\kappa)\eta}
 \bigl(D^w\mathcal S_{k,1}+(1-\kappa)\mathcal Q_{k,w},
       -D^w\mathcal S_{k,2}+4\pi(1-\kappa)\mathcal Q_{k,w}\bigr)
 \right\|_2\\
 &\quad\leq
 \sum_{j=1}^2\|e^{\lambda_{k,s}(1-\kappa)\eta}
                         D^w\mathcal S_{k,j}\|_2
 +C\|e^{\lambda_{k,s}(1-\kappa)\eta}
                         (1-\kappa)\mathcal Q_{k,w}\|_2\\
 &\quad\leq C_{k,s}
 (1-\kappa)^{-\left(k+\frac72s+\frac72k+5\cdot4^{k-1}+\frac{59}{2}\right)},\\
 &\|e^{\lambda_{k,s}(1-\kappa)\eta}
       (D^w\mathcal S_{k,1}+(1-\kappa)\mathcal Q_{k,w})\|_\infty\\
 &\quad\leq
 \|e^{\lambda_{k,s}(1-\kappa)\eta}D^w\mathcal S_{k,1}\|_\infty
 +\|e^{\lambda_{k,s}(1-\kappa)\eta}
                         (1-\kappa)\mathcal Q_{k,w}\|_\infty\\
 &\quad\leq C_{k,s}
 (1-\kappa)^{-\left(k+\frac72s+\frac72k+5\cdot4^{k-1}+29\right)}.
\end{align*}
For the incoming datum, \eqref{ldb},
\eqref{pbr}, and
\eqref{lcb} yield
\begin{align*}
 &\|\widehat U_{k,w}(0)\|_{L^\infty(\gamma_-)}\\
 &\quad\leq\|D^w\phi_k(0)\|_{L^\infty(\gamma_-)}
 +\left\|\sum_{0<u\leq w}\binom wu
      D^u(4a^3)D^{w-u}\mathfrak g_{k,\infty}\right\|_\infty\\
 &\quad\leq C_{k,s}
 (1-\kappa)^{-\left(k+\frac72s+5\cdot4^{k-1}+30\right)}\\
 &\qquad+C_{k,s}\sum_{0<u\leq w}(1-\kappa)^{2-\frac72|u|}
 (1-\kappa)^{-\left(k+\frac72(s-|u|)+\frac72k
                              +5\cdot4^{k-1}+31\right)}\\
 &\quad\leq C_{k,s}
 (1-\kappa)^{-\left(k+\frac72s+\frac72k+5\cdot4^{k-1}+29\right)}.
\end{align*}
Substituting the bounds for $S_1,S_2$ and
$\phi_b=\widehat U_{k,w}(0)|_{\gamma_-}$ into~\eqref{le},
summing over $|w|=s$, and adding the lower orders from the induction
hypothesis gives
\begin{align}
 &\sum_{|w|\leq s}\left(\|D^w\mathfrak g_{k,\infty}\|_\infty
 +\|e^{\lambda_{k,|w|}(1-\kappa)\eta}D^w\mathcal T_k\|_\infty\right)
 \nonumber\\
 &\quad\leq\frac{C_{k,s}}{(1-\kappa)^{3/2}}
 (1-\kappa)^{-\left(k+\frac72s+\frac72k+5\cdot4^{k-1}+\frac{59}{2}\right)}
 \nonumber\\
 &\quad\leq C_{k,s}
 (1-\kappa)^{-\left(k+\frac72s+\frac72k+5\cdot4^{k-1}+31\right)},
 \label{lht}\\
 &\sum_{|w|\leq s}
 \|e^{\lambda_{k,|w|}(1-\kappa)\eta}D^w\mathcal I_k\|_\infty
 \nonumber\\
 &\quad\leq\frac{C_{k,s}}{(1-\kappa)^{5/2}}
 (1-\kappa)^{-\left(k+\frac72s+\frac72k+5\cdot4^{k-1}+\frac{59}{2}\right)}
 \nonumber\\
 &\qquad+\frac{C_{k,s}}{1-\kappa}
 (1-\kappa)^{-\left(k+\frac72s+\frac72k+5\cdot4^{k-1}+29\right)}
 \nonumber\\
 &\quad\leq C_{k,s}
 (1-\kappa)^{-\left(k+\frac72s+\frac72k+5\cdot4^{k-1}+32\right)}.
 \label{lhr}
\end{align}
Together with the case $s=0$, this closes the tangential induction
for $0\leq s\leq s_0-3k$. By the definition~\eqref{pk},
\begin{align}\label{6.15}
 \mathfrak A_{k,s}(\kappa)
 \leq C_{k,s}
 (1-\kappa)^{-\left(k+\frac72s+\frac72k+5\cdot4^{k-1}+32\right)}
 \leq C_{k,s}
 (1-\kappa)^{-\left(k+\frac72s+5\cdot4^k\right)},
 \qquad k\geq2.
\end{align}

The estimate~\eqref{6.15} holds for
$0\leq s\leq s_0-3k$.
By the definitions~\eqref{pwn} and~\eqref{pk},
for every $|w|=s$ in this range, we have
\begin{align}\label{6.20}
	|D^w\mathcal T_k(\eta,x_\shortparallel)|
	+|D^w\mathcal I_k(\eta,x_\shortparallel,\beta)|
	\leq
	\mathfrak A_{k,s}(\kappa)
	e^{-\lambda_{k,s}(1-\kappa)\eta}.
\end{align}
Since $\chi$ is independent of $x_\shortparallel$
and $|\chi|\leq1$, the definition~\eqref{Bdk}
gives the truncated estimate~\eqref{pc}.
The available tangential orders include $r+2$
for $0\leq r\leq s_0-3k-2$, as required for the interior Dirichlet data.

\medskip\noindent\emph{Step 3.2. Upper boundary layer.}
Fix $|w|=s\leq s_0-3k$ and take $\vartheta=\lambda_{k,s}(1-\kappa)$.
After reflecting $\beta_3\mapsto-\beta_3$, the a priori estimates~\eqref{le}
apply to the upper profiles with $\xi$ in place of $\eta$ and
$c^+=4T_b^3$ in place of $4a^3$. In the original angular variable,
\eqref{uds}--\eqref{udb} give $g_b=0$ and
\begin{align}
 S_1&=D^w\mathcal S_{k,1}^++(1-\kappa)\mathcal Q_{k,w}^+,\nonumber\\
 S_2&=-D^w\mathcal S_{k,2}^++4\pi(1-\kappa)\mathcal Q_{k,w}^+,\nonumber\\
 \phi_b&=D^w\phi_k^+(0,\beta)
 -\sum_{0<u\leq w}\binom wu D^u c^+D^{w-u}\mathfrak g_{k,\infty}^+,
 \qquad \beta_3<0. \label{ukdata}
\end{align}
In~\eqref{ucm} and in the boundary sum in~\eqref{ukdata}, the derivatives
on $\mathcal T_k^+$ and $\mathfrak g_{k,\infty}^+$ have order less than $s$.
Applying $D^w$ to~\eqref{up} and~\eqref{uki}
gives the prescribed terms
\begin{align*}
 D^w\mathcal S_{k,1}^+
 &=(1-\kappa)D^w\bigl[\mathcal P^+(\mathcal T^+,k)-c^+\mathcal T_k^+\bigr]
 -\beta_\shortparallel\cdot\nabla_\shortparallel D^w\mathcal I_{k-1}^+,\\
 D^w\mathcal S_{k,2}^+
 &=-4\pi(1-\kappa)D^w\bigl[\mathcal P^+(\mathcal T^+,k)-c^+\mathcal T_k^+\bigr]
 +\Delta_\shortparallel D^w\mathcal T_{k-2}^+,\\
 D^w\phi_k^+(0,\beta)
 &=\left.D^w\left[\beta\cdot\nabla I_{k-1}
 +\frac{\kappa}{4\pi(1-\kappa)}\Delta T_{k-2}
 -\sum_{\substack{i_1+\cdots+i_4=k\\0\leq i_j<k}}\prod_{j=1}^4 T_{i_j}
 \right]\right|_1,\qquad \beta_3<0.
\end{align*}

Since $\mathcal T_0^+=0$, every nonzero product in~\eqref{ust} contains
a layer factor $\mathcal T_j^+$ with $1\leq j<k$. The Taylor coefficients
in~\eqref{uq} contribute $(-\xi)^\ell$, and the normal derivative
estimate~\eqref{itn} holds also at $x_3=1$. Summing the product
estimate~\eqref{qwt} over the terms in the Leibniz expansion, with
the upper-layer bounds from~\eqref{6.11}, gives
\begin{align}
 &\left\|e^{\lambda_{k-1,s+1}(1-\kappa)\xi}(1-\kappa)
 D^w\bigl[\mathcal P^+(\mathcal T^+,k)-c^+\mathcal T_k^+\bigr]\right\|_\infty
 \nonumber\\
 &\qquad\leq C_{k,s}(1-\kappa)^{-\left(k+\frac72s+\frac72k+5\cdot4^{k-1}+29\right)}.
 \label{uspb}
\end{align}
In~\eqref{up}, the tangential transport and Laplacian terms are bounded
by $C_s\mathfrak A^+_{k-1,s+1}$ and $C_s\mathfrak A^+_{k-2,s+2}$
in the weighted $L^\infty$ norm of~\eqref{uspb}. The Laplacian term vanishes at $k=2$.
Combining~\eqref{uspb} with~\eqref{6.11} as in~\eqref{6.13}, and
integrating the weight as in~\eqref{swt}, we obtain, for $j=1,2$,
\begin{align}\label{usb}
 \|e^{\lambda_{k-1,s+1}(1-\kappa)\xi}D^w\mathcal S_{k,j}^+\|_\infty
 &\leq C_{k,s}(1-\kappa)^{-\left(k+\frac72s+\frac72k+5\cdot4^{k-1}+29\right)},\nonumber\\
 \|e^{\vartheta\xi}D^w\mathcal S_{k,j}^+\|_2
 &\leq C_{k,s}(1-\kappa)^{-\left(k+\frac72s+\frac72k+5\cdot4^{k-1}+\frac{59}{2}\right)}.
\end{align}
Applying the trace and product estimates of~\eqref{pbr} at $x_3=1$ to
the incoming datum~\eqref{uki} gives
\begin{equation}\label{uib}
 \max_{|v|\leq s}\|D^v\phi_k^+(0)\|_{L^\infty(\{\beta_3<0\})}
 \leq C_{k,s}(1-\kappa)^{-\left(k+\frac72s+5\cdot4^{k-1}+30\right)} .
\end{equation}

At $s=0$, $\mathcal Q_{k,0}^+=0$ and $\phi_b=\phi_k^+(0,\beta)$.
Since $T_b$ is independent of $\xi$, the
spectral condition~\eqref{asA} holds. With the data bounds
\eqref{usb}--\eqref{uib}, Theorem~\ref{co} gives the unique bounded
solution $(\phi_k^+,g_k^+)$ and its far-field temperature
$\mathfrak g_{k,\infty}^+$. Substituting those bounds into~\eqref{le}
gives, for the centered profiles in~\eqref{cen},
\begin{align*}
 \|\mathfrak g_{k,\infty}^+\|_\infty
 +\|e^{\lambda_{k,0}(1-\kappa)\xi}\mathcal T_k^+\|_\infty
 &\leq C_k(1-\kappa)^{-\left(k+\frac72k+5\cdot4^{k-1}+31\right)},\\
 \|e^{\lambda_{k,0}(1-\kappa)\xi}\mathcal I_k^+\|_\infty
 &\leq C_k(1-\kappa)^{-\left(k+\frac72k+5\cdot4^{k-1}+32\right)}.
\end{align*}

Fix $1\leq s\leq s_0-3k$ and assume the second bound in~\eqref{uqb}
at all lower tangential orders,
\[
 \mathfrak A^+_{k,r}(\kappa)\leq C_{k,r}
 (1-\kappa)^{-\left(k+\frac72r+\frac72k+5\cdot4^{k-1}+32\right)},
 \qquad 0\leq r<s .
\]
For $0<u\leq w$, we have $r=|w-u|=s-|u|\leq s-1$.
Since $\|D^uc^+\|_\infty\leq C_s$, substituting this induction bound
at $r=s-|u|$ into the three inequalities of~\eqref{ucn} gives
\begin{align}
 &\|e^{\lambda_{k,s-1}(1-\kappa)\xi}(1-\kappa)\mathcal Q_{k,w}^+\|_\infty \nonumber\\
 &\quad\leq C_{k,s}(1-\kappa)\sum_{0<u\leq w}\mathfrak A^+_{k,s-|u|}
 \leq C_{k,s}(1-\kappa)^{-\left(k+\frac72s+\frac72k+5\cdot4^{k-1}+\frac{55}{2}\right)},\nonumber\\
 &\|e^{\vartheta\xi}(1-\kappa)\mathcal Q_{k,w}^+\|_2 \nonumber\\
 &\quad\leq C_{k,s}(1-\kappa)^{1/2}\sum_{0<u\leq w}\mathfrak A^+_{k,s-|u|}
 \leq C_{k,s}(1-\kappa)^{-\left(k+\frac72s+\frac72k+5\cdot4^{k-1}+28\right)},\nonumber\\
 &\left\|\sum_{0<u\leq w}\binom wu D^uc^+D^{w-u}\mathfrak g_{k,\infty}^+\right\|_\infty \nonumber\\
 &\quad\leq C_{k,s}\sum_{0<u\leq w}\mathfrak A^+_{k,s-|u|}
 \leq C_{k,s}(1-\kappa)^{-\left(k+\frac72s+\frac72k+5\cdot4^{k-1}+\frac{57}{2}\right)} .
 \label{ukcorr}
\end{align}
In~\eqref{ukdata}, we use~\eqref{usb} for the prescribed sources,
\eqref{uib} for the prescribed incoming datum, and~\eqref{ukcorr}
for the coefficient terms. Since
$\lambda_{k,s}<\min\{\lambda_{k,s-1},\lambda_{k-1,s+1}\}$ by~\eqref{plh},
the triangle inequality gives
\begin{align*}
 \|e^{\vartheta\xi}(S_1,S_2)\|_2
 &\leq C_{k,s}(1-\kappa)^{-\left(k+\frac72s+\frac72k+5\cdot4^{k-1}+\frac{59}{2}\right)},\\
 \|e^{\vartheta\xi}S_1\|_\infty
 &\leq C_{k,s}(1-\kappa)^{-\left(k+\frac72s+\frac72k+5\cdot4^{k-1}+29\right)},\\
 \|\phi_b\|_{L^\infty(\{\beta_3<0\})}
 &\leq C_{k,s}(1-\kappa)^{-\left(k+\frac72s+\frac72k+5\cdot4^{k-1}+\frac{57}{2}\right)} .
\end{align*}
Substituting these bounds into~\eqref{le}, summing over $|w|=s$, and
including the lower orders by the induction hypothesis, we obtain
\begin{align}
 \sum_{|w|\leq s}\left(\|D^w\mathfrak g_{k,\infty}^+\|_\infty
 +\|e^{\lambda_{k,|w|}(1-\kappa)\xi}D^w\mathcal T_k^+\|_\infty\right)
 &\leq C_{k,s}(1-\kappa)^{-\left(k+\frac72s+\frac72k+5\cdot4^{k-1}+31\right)},
 \label{uht}\\
 \sum_{|w|\leq s}\|e^{\lambda_{k,|w|}(1-\kappa)\xi}D^w\mathcal I_k^+\|_\infty
 &\leq C_{k,s}(1-\kappa)^{-\left(k+\frac72s+\frac72k+5\cdot4^{k-1}+32\right)} .
 \label{uhr}
\end{align}
Adding~\eqref{uht} and~\eqref{uhr} in~\eqref{upn} closes the induction.
Since $\frac72k+32\leq15\cdot4^{k-1}$ for $k\geq2$, we have
\begin{align}\label{uhb}
 \mathfrak A^+_{k,s}(\kappa)
 &\leq C_{k,s}(1-\kappa)^{-\left(k+\frac72s+\frac72k+5\cdot4^{k-1}+32\right)}
 \nonumber\\
 &\leq C_{k,s}(1-\kappa)^{-\left(k+\frac72s+5\cdot4^k\right)},
 \qquad 0\leq s\leq s_0-3k .
\end{align}
Finally, the cutoff in~\eqref{ukc} is bounded by one and independent
of $x_\shortparallel$, so~\eqref{upn} gives~\eqref{upc}.

\medskip\noindent\emph{Step 3.3. Interior profiles.}
Fix $0\leq r\leq s_0-3k-2$. The far-field temperatures
$\mathfrak g_{k,\infty}$ and $\mathfrak g_{k,\infty}^+$ from Steps~3.1
and~3.2 give the two Dirichlet values of $T_k$ in~\eqref{ikm}.
Separating the terms containing $T_k$ in the quartic coefficient
$\mathcal C(T,k)$, we obtain
\begin{align}\label{iq}
 \mathcal C(T,k)=4T_0^3T_k+
 \sum_{\substack{i_1+\cdots+i_4=k\\0\leq i_j<k}}
 T_{i_1}T_{i_2}T_{i_3}T_{i_4}.
\end{align}
Substituting~\eqref{iq} into~\eqref{ikm} and using
$c_0=\mathcal H'(T_0)=1+\frac{16\pi}{3}T_0^3$, we obtain
\begin{align}\label{ip}
 \Delta(c_0T_k)
 &=-\frac{4\pi}{3}\Delta
 \sum_{\substack{i_1+\cdots+i_4=k\\0\leq i_j<k}}
 T_{i_1}T_{i_2}T_{i_3}T_{i_4}
 +\int_{\mathbb S^2}(\beta\cdot\nabla)^3I_{k-1}\,d\beta
 +\frac{\kappa}{3(1-\kappa)}\Delta^2T_{k-2},\\
 (c_0T_k)|_{x_3=0}&=\mathcal H'(a)\mathfrak g_{k,\infty},
 \qquad (c_0T_k)|_{x_3=1}=\mathcal H'(T_b)\mathfrak g_{k,\infty}^+.\nonumber
\end{align}
The Dirichlet problem~\eqref{ip} determines $c_0T_k$ uniquely.
Since $c_0\geq1$, division by $c_0$ gives $T_k$, and~\eqref{3.9}
then determines $I_k$.

Applying the $W^{r+2,q}$ Dirichlet estimate to~\eqref{ip} and bounding
the $W^{r+2-1/q,q}$ norms of the boundary data by their $W^{r+2,\infty}$ norms,
we obtain
\begin{align}
 \|c_0T_k\|_{W^{r+2,q}}
 &\leq C_{k,r,q}\Biggl[
 \|\mathcal H'(a)\mathfrak g_{k,\infty}\|_{W^{r+2,\infty}}
 +\|\mathcal H'(T_b)\mathfrak g_{k,\infty}^+\|_{W^{r+2,\infty}}
 \nonumber\\
 &\qquad+\sum_{\substack{i_1+\cdots+i_4=k\\0\leq i_j<k}}
 \left\|\prod_{j=1}^4T_{i_j}\right\|_{W^{r+2,q}}
 +\sup_{\beta\in\mathbb S^2}\|I_{k-1}(\cdot,\beta)\|_{W^{r+3,q}}
 \nonumber\\
 &\qquad+\frac{\kappa}{1-\kappa}\|T_{k-2}\|_{W^{r+4,q}}\Biggr].
 \label{ikdir}
\end{align}
Both boundary-layer estimates~\eqref{6.15} and~\eqref{uhb} apply through
order $r+2$ since $r+2\leq s_0-3k$. At $x_3=0$, the factor
$\mathcal H'(a)$ is uniformly bounded. For $1\leq\ell\leq r+2$, the
composition estimate and~\eqref{6.0-0} give
\[
 \|D^\ell\mathcal H'(a)\|_{L^\infty}
 \leq C_\ell(1-\kappa)^{-\frac72\ell+2}.
\]
Applying Leibniz' rule to the lower boundary value and using the
first inequality in~\eqref{6.15}, we obtain
\begin{align*}
 \|\mathcal H'(a)\mathfrak g_{k,\infty}\|_{W^{r+2,\infty}}
 &\le C_{k,r}\left[\mathfrak A_{k,r+2}(\kappa)
 +\sum_{\ell=1}^{r+2}(1-\kappa)^{-\frac72\ell+2}
 \mathfrak A_{k,r+2-\ell}(\kappa)\right]\nonumber\\
 &\leq C_{k,r}(1-\kappa)^{-\left(k+\frac72r+\frac72k
 +5\cdot4^{k-1}+39\right)}.
\end{align*}

At $x_3=1$, the derivatives of $\mathcal H'(T_b)$ through order $r+2$
are bounded by the prescribed data. Leibniz' rule and the first
inequality in~\eqref{uhb} therefore give
\begin{align}\label{uitb}
 \|\mathcal H'(T_b)\mathfrak g_{k,\infty}^+\|_{W^{r+2,\infty}}
 \leq C_{k,r}\sum_{\ell=0}^{r+2}\mathfrak A^+_{k,r+2-\ell}(\kappa)\leq C_{k,r}(1-\kappa)^{-\left(k+\frac72r+\frac72k+5\cdot4^{k-1}+39\right)}.
\end{align}

For the sum in~\eqref{iq}, we use the tame product estimate with
$W^{2,q}(\Omega)\hookrightarrow L^\infty(\Omega)$. Applying~\eqref{6.11}
at order $r$ to one factor and at order zero to the other three gives
\begin{align}
 \sum_{\substack{i_1+\cdots+i_4=k\\0\leq i_j<k}}
 \left\|\prod_{j=1}^4T_{i_j}\right\|_{W^{r+2,q}}&\leq C_{k,r,q}
 \sum_{\substack{i_1+\cdots+i_4=k\\0\leq i_j<k}}
 \sum_{\nu=1}^4\|T_{i_\nu}\|_{W^{r+2,q}}
 \prod_{j\ne\nu}\|T_{i_j}\|_{W^{2,q}}\nonumber\\
 &\leq C_{k,r,q}(1-\kappa)^{-\left(k+\frac72r
 +5\cdot4^{k-1}+30\right)}. \label{ikprod}
\end{align}
Here $\sum_{j=1}^4 i_j=k$ and $\sum_{j=1}^4 4^{i_j}\leq4^{k-1}+6$.
For the last two terms of~\eqref{ikdir}, the definition~\eqref{pin}
gives the bounds $\mathfrak E_{k-1,r+2}$ and
$\frac{\kappa}{1-\kappa}\mathfrak E_{k-2,r+2}$, respectively.
Using~\eqref{6.11} at order $r+2$ for these two terms and inserting
the trace bounds from~\eqref{6.15} and~\eqref{uitb} and the product
bound~\eqref{ikprod} into~\eqref{ikdir}, we obtain
\begin{align}\label{ipe}
 \|c_0T_k\|_{W^{r+2,q}}
 &\leq C_{k,r,q}\Bigg\{
 (1-\kappa)^{-\left(k+\frac72r+\frac72k+5\cdot4^{k-1}+39\right)}
 +(1-\kappa)^{-\left(k+\frac72r+5\cdot4^{k-1}+30\right)}\nonumber\\
 &\qquad
 +(1-\kappa)^{-\left((k-1)+\frac72(r+2)+5\cdot4^{k-1}\right)}
 +(1-\kappa)^{-\left((k-2)+\frac72(r+2)+5\cdot4^{k-2}+1\right)}
 \Bigg\}\nonumber\\
 &\leq C_{k,r,q}(1-\kappa)^{-\left(k+\frac72r+\frac72k
 +5\cdot4^{k-1}+39\right)}.
\end{align}

Since $c_0\geq1$, we have $\|c_0^{-1}\|_{L^\infty}\leq1$.
Applying the tame composition estimate to $c_0^{-1}=1/\mathcal H'(T_0)$
and using~\eqref{ell} at order $\ell-1$, we obtain
\[
 \|D^\ell c_0^{-1}\|_{L^\infty}
 \leq C_{\ell,q}(1-\kappa)^{-\frac72\ell-\frac32},
 \quad 1\leq\ell\leq r+2.
\]
Applying Leibniz' rule to $T_k=c_0^{-1}(c_0T_k)$, we use~\eqref{ipe}
at order $r-\ell$ for $1\leq\ell\leq r$ and at order zero for
$\ell=r+1,r+2$. This gives
\begin{align}\label{pir}
 \|T_k\|_{W^{r+2,q}}
 &\leq C_{k,r,q}\|c_0T_k\|_{W^{r+2,q}}
 +C_{k,r,q}\sum_{\ell=1}^{r}
 (1-\kappa)^{-\frac72\ell-\frac32}
 \|c_0T_k\|_{W^{r+2-\ell,q}}\nonumber\\
 &\quad+C_{k,r,q}\sum_{\ell=r+1}^{r+2}
 (1-\kappa)^{-\frac72\ell-\frac32}
 \|c_0T_k\|_{W^{2,q}}\nonumber\\
 &\leq C_{k,r,q}(1-\kappa)^{-\left(k+\frac72r+\frac72k
 +5\cdot4^{k-1}+\frac{95}{2}\right)}\nonumber\\
 &\leq C_{k,r,q}(1-\kappa)^{-\left(k+\frac72r+5\cdot4^k\right)}.
\end{align}
The first sum is empty when $r=0$.

For the radiation profile, substituting~\eqref{iq} into~\eqref{3.9} gives
\[
 I_k=-\beta\cdot\nabla I_{k-1}
 -\frac{\kappa}{4\pi(1-\kappa)}\Delta T_{k-2}
 +4T_0^3T_k
 +\sum_{\substack{i_1+\cdots+i_4=k\\0\leq i_j<k}}\prod_{j=1}^4T_{i_j}.
\]
In $W^{r+1,q}$, the first two terms are bounded by
$C\mathfrak E_{k-1,r+1}$ and
$C\frac{\kappa}{1-\kappa}\mathfrak E_{k-2,r+1}$.
For the third term, \eqref{ell} and the composition estimate give
$\|T_0^3\|_{W^{r+1,\infty}}\leq C_{r,q}(1-\kappa)^{-\frac72r-5}$.
Using the tame product estimate and the intermediate bound in~\eqref{pir}
at orders $(r-1)_+$ and zero, we obtain
\begin{align*}
 \|4T_0^3T_k\|_{W^{r+1,q}}
 &\leq C_{k,r,q}\left[
 \|T_k\|_{W^{r+1,q}}
 +\|T_0^3\|_{W^{r+1,\infty}}\|T_k\|_{L^q}\right]\\
 &\leq C_{k,r,q}(1-\kappa)^{-\left(k+\frac72r+\frac72k
 +5\cdot4^{k-1}+\frac{105}{2}\right)}.
\end{align*}
For the quartic sum, \eqref{ikprod} at order $(r-1)_+$ controls
the $W^{r+1,q}$ norm, including $r=0$ since $W^{2,q}\hookrightarrow W^{1,q}$.
Taking the $W^{r+1,q}$ norm in~\eqref{3.9} and using~\eqref{6.11}
at order $r+1$ for the transport and Laplacian terms yields
\begin{align*}
 \sup_\beta\|I_k(\cdot,\beta)\|_{W^{r+1,q}}
 &\leq C_{k,r,q}\Bigg[
 (1-\kappa)^{-\left((k-1)+\frac72(r+1)+5\cdot4^{k-1}\right)}\\
 &\quad +(1-\kappa)^{-\left((k-2)+\frac72(r+1)+5\cdot4^{k-2}+1\right)}\\
 &\quad +(1-\kappa)^{-\left(k+\frac72r+\frac72k
 +5\cdot4^{k-1}+\frac{105}{2}\right)}\\
 &\quad +(1-\kappa)^{-\left(k+\frac72(r-1)_+
 +5\cdot4^{k-1}+30\right)}\Bigg]\\
 &\leq C_{k,r,q}(1-\kappa)^{-\left(k+\frac72r+5\cdot4^k\right)}.
\end{align*}
The last inequality uses $\frac72k+\frac{105}{2}\leq15\cdot4^{k-1}$
for $k\geq2$. Adding the radiation bound to~\eqref{pir} in the
definition~\eqref{pin} proves~\eqref{per} at order $k$.
This completes the induction on $k$.

\medskip\noindent{\it Step 4. Normal derivatives of the boundary-layer temperatures.}

\medskip\noindent{\it Step 4.1. Lower boundary layers.}
For $k\geq1$, taking the angular integral of $\mathcal S_{k,1}$ in~\eqref{3.19}
and adding $\mathcal S_{k,2}$ cancels the quartic terms and gives
\[
 \langle\mathcal S_{k,1}\rangle+\mathcal S_{k,2}
 =-\nabla_\shortparallel\cdot
       \langle\beta_\shortparallel\mathcal I_{k-1}\rangle
   +\Delta_\shortparallel\mathcal T_{k-2}.
\]
Here $\mathcal T_{-1}=0$. By~\eqref{pwn} and~\eqref{pk}, the tangential
bounds~\eqref{pzr} and~\eqref{6.11} make the right-hand side integrable in $\eta$.
Integrating the first equation in~\eqref{3.18} over $\mathbb S^2$
and adding the second gives
\[
 \partial_\eta\bigl(\langle\beta_3\phi_k\rangle-\partial_\eta g_k\bigr)
 =\langle\mathcal S_{k,1}\rangle+\mathcal S_{k,2}.
\]
Hence $\langle\beta_3\phi_k\rangle-\partial_\eta g_k$ has a finite limit
as $\eta\to\infty$. Since the far field of $\phi_k$ is isotropic,
$\langle\beta_3\phi_k\rangle=\langle\beta_3\mathcal I_k\rangle\to0$.
Thus $\partial_\eta g_k$ has a finite limit, which is zero because $g_k$
is bounded. Integrating the moment identity over $(\eta,\infty)$ and
using $\partial_\eta\mathcal T_k=\partial_\eta g_k$, we obtain
\[
 \partial_\eta\mathcal T_k(\eta)
 =\langle\beta_3\mathcal I_k(\eta)\rangle
 +\int_\eta^\infty
   \bigl(\langle\mathcal S_{k,1}\rangle+\mathcal S_{k,2}\bigr)(z)\,dz.
\]
For $k\geq2$, applying~\eqref{pwn} and~\eqref{pk} to the two terms
in the source sum and using the ordering~\eqref{plh}, we obtain
\begin{align*}
 &\left|\int_\eta^\infty
   \bigl(\langle\mathcal S_{k,1}\rangle+\mathcal S_{k,2}\bigr)(z)\,dz\right|\\
 &\quad\leq C\mathfrak A_{k-1,1}(\kappa)
       \int_\eta^\infty e^{-\lambda_{k-1,1}(1-\kappa)z}\,dz
 +C\mathfrak A_{k-2,2}(\kappa)
       \int_\eta^\infty e^{-\lambda_{k-2,2}(1-\kappa)z}\,dz\\
 &\quad\leq\frac{C_k}{1-\kappa}
 \bigl(\mathfrak A_{k-1,1}(\kappa)+\mathfrak A_{k-2,2}(\kappa)\bigr)
 e^{-\lambda_{k,0}(1-\kappa)\eta}.
\end{align*}
For $k=1$, the same integration contains only the term with
$\mathfrak A_{0,1}$. Using~\eqref{pzr}
and~\eqref{6.7}, we obtain
\begin{align*}
 \|e^{\lambda_{1,0}(1-\kappa)\eta}\partial_\eta\mathcal T_1\|_\infty
 &\leq C\left[\mathfrak A_{1,0}+(1-\kappa)^{-1}\mathfrak A_{0,1}\right]\\
 &\leq C\left[(1-\kappa)^{-9}+(1-\kappa)^{-7/2}\right]
 \leq C(1-\kappa)^{-9}.
\end{align*}
For $k\geq2$, the bound~\eqref{6.15} at order zero and the common
lower-order bounds~\eqref{6.11} give
\begin{align*}
 &\|e^{\lambda_{k,0}(1-\kappa)\eta}\partial_\eta\mathcal T_k\|_\infty\\
 &\quad\leq C_k\left[\mathfrak A_{k,0}
 +(1-\kappa)^{-1}
       \bigl(\mathfrak A_{k-1,1}+\mathfrak A_{k-2,2}\bigr)\right]\\
 &\quad\leq C_k\Bigl[
 (1-\kappa)^{-\left(k+\frac72k+5\cdot4^{k-1}+32\right)}
 +(1-\kappa)^{-\left(k+\frac72+5\cdot4^{k-1}\right)}
 +(1-\kappa)^{-\left(k+6+5\cdot4^{k-2}\right)}\Bigr]\\
 &\quad\leq C_k
 (1-\kappa)^{-\left(k+\frac72k+5\cdot4^{k-1}+32\right)}.
\end{align*}
For $k=0$, $\partial_\eta\mathcal T_0=\partial_\eta\mathbf T$.
Taking $\alpha=\lambda_{0,0}(1-\kappa)$ in~\eqref{4.60} gives
\[
 \|e^{\lambda_{0,0}(1-\kappa)\eta}\partial_\eta\mathcal T_0\|_\infty
 \leq C\frac{1-\kappa}{\sqrt{\alpha(1-\kappa-\alpha)}}
 =\frac{C}{\sqrt{\lambda_{0,0}(1-\lambda_{0,0})}}.
\]
Combining these derivative bounds with the temperature bounds
in~\eqref{4.56-1}, \eqref{6.7}, and~\eqref{6.15} proves~\eqref{nlb}.

\medskip\noindent{\it Step 4.2. Upper boundary layers.}
For $k\geq1$, taking the angular integral of $\mathcal S_{k,1}^+$ in~\eqref{up}
and adding $\mathcal S_{k,2}^+$ gives
\[
 \langle\mathcal S_{k,1}^+\rangle+\mathcal S_{k,2}^+
 =-\nabla_\shortparallel\cdot
       \langle\beta_\shortparallel\mathcal I_{k-1}^+\rangle
   +\Delta_\shortparallel\mathcal T_{k-2}^+.
\]
The leading upper layer vanishes, and~\eqref{upn} and~\eqref{uqb}
give exponential decay for the positive-order terms. Thus the sum is integrable in $\xi$.
Integrating the first equation in~\eqref{uks} over $\mathbb S^2$,
adding the second, and using~\eqref{cen}, we obtain
\[
 \partial_\xi\bigl(-\langle\beta_3\mathcal I_k^+\rangle
                   -\partial_\xi\mathcal T_k^+\bigr)
 =\langle\mathcal S_{k,1}^+\rangle+\mathcal S_{k,2}^+.
\]
As at the lower boundary, integrability of the right-hand side and
$\langle\beta_3\mathcal I_k^+\rangle\to0$ imply that
$\partial_\xi\mathcal T_k^+$ has a finite limit. Boundedness of
$\mathcal T_k^+$ forces this limit to be zero. Integrating over
$(\xi,\infty)$ gives
\begin{align}\label{usp}
 \partial_\xi\mathcal T_k^+
 &=-\langle\beta_3\mathcal I_k^+\rangle
   +\int_\xi^\infty
      \bigl(\langle\mathcal S_{k,1}^+\rangle+\mathcal S_{k,2}^+\bigr)(z)\,dz.
\end{align}
At $k=1$, the source sum vanishes since $\mathcal I_0^+=\mathcal T_{-1}^+=0$.
Thus $\partial_\xi\mathcal T_1^+=-\langle\beta_3\mathcal I_1^+\rangle$.
The weighted $L^\infty$ estimate for $\mathcal I_1^+$ in Step~2.2 gives
\[
 \|e^{\lambda_{1,0}(1-\kappa)\xi}\partial_\xi\mathcal T_1^+\|_\infty
 \leq C\|e^{\lambda_{1,0}(1-\kappa)\xi}\mathcal I_1^+\|_\infty
 \leq C(1-\kappa)^{-15/2}.
\]
For $k\geq2$, applying the same exponential tail estimate to~\eqref{usp}
and using~\eqref{uhb} and~\eqref{6.11}, we obtain
\begin{align*}
 &\|e^{\lambda_{k,0}(1-\kappa)\xi}\partial_\xi\mathcal T_k^+\|_\infty\\
 &\quad\leq C_k\left[\mathfrak A^+_{k,0}
 +(1-\kappa)^{-1}
   \bigl(\mathfrak A^+_{k-1,1}+\mathfrak A^+_{k-2,2}\bigr)\right]\leq C_k(1-\kappa)^{-\left(k+\frac72k+5\cdot4^{k-1}+32\right)},
\end{align*}
where $\mathfrak A^+_{0,s}=0$. Together with the first-order bound,
this proves~\eqref{unb}.
\end{proof}

\section{Remainder estimate}\label{res}
\subsection{Exact residual decomposition}

For the matched approximate solution, substituting the interior and
layer equations into~\eqref{rd} gives
\begin{align}\label{fr}
 R_1={}&\v^{N+1}\beta\cdot\nabla I_N
       -(1-\kappa)\sum_{k=N+1}^{4N}\v^k\mathcal C(T,k)
       +\mathcal E_1+\mathcal E_1^+,\nonumber\\
 R_2={}&-\v^{N+1}\Delta T_{N-1}-\v^{N+2}\Delta T_N
       +4\pi(1-\kappa)\sum_{k=N+1}^{4N}\v^k\mathcal C(T,k)
       +\mathcal E_2+\mathcal E_2^+.
\end{align}
For the lower contributions, writing $\chi=\chi(\v\eta)$, we have
\begin{align}\label{lr}
 \mathcal E_1={}&\v^{N+1}\chi\,
             \beta_\shortparallel\cdot\nabla_\shortparallel\mathcal I_N
       +\partial_\eta\chi\sum_{k=0}^N\v^k\beta_3\mathcal I_k
       +(1-\kappa)\mathcal{R}_N,\nonumber\\
 \mathcal E_2={}&-\v^{N+1}\chi\Delta_\shortparallel\mathcal T_{N-1}
       -\v^{N+2}\chi\Delta_\shortparallel\mathcal T_N\nonumber\\
       &-2\partial_\eta\chi\sum_{k=0}^N\v^k\partial_\eta\mathcal T_k
       -\partial_\eta^2\chi\sum_{k=0}^N\v^k\mathcal T_k
       -4\pi(1-\kappa)\mathcal{R}_N.
\end{align}
For the upper contributions, writing $\chi=\chi(\v\xi)$, we have
\begin{align}\label{ur}
 \mathcal E_1^+
 ={}&\v^{N+1}\chi\,
               \beta_\shortparallel\cdot\nabla_\shortparallel\mathcal I_N^+
       -\partial_\xi\chi\sum_{k=1}^N\v^k\beta_3\mathcal I_k^+
       +(1-\kappa)\mathcal{R}_N^+,\nonumber\\
 \mathcal E_2^+
 ={}&-\v^{N+1}\chi\Delta_\shortparallel\mathcal T_{N-1}^+
       -\v^{N+2}\chi\Delta_\shortparallel\mathcal T_N^+\nonumber\\
    &-2\partial_\xi\chi\sum_{k=1}^N\v^k\partial_\xi\mathcal T_k^+
       -\partial_\xi^2\chi\sum_{k=1}^N\v^k\mathcal T_k^+
       -4\pi(1-\kappa)\mathcal{R}_N^+.
\end{align}
Here $\mathcal{R}_N$ and $\mathcal{R}_N^+$ collect the quartic remainders associated
with the two layers. To write them explicitly, set
\begin{align*}
 U^N&=\sum_{k=0}^N\v^kT_k,\qquad
 Q^N=\sum_{k=0}^N\v^kQ_k,\qquad
 V^N=\sum_{k=0}^N\v^k\mathcal T_k,\\
 Q_+^N&=\sum_{k=0}^N\v^kQ_k^+,\qquad
 V_+^N=\sum_{k=1}^N\v^k\mathcal T_k^+,
 \qquad \mathfrak F(X,Y)=(X+Y)^4-X^4.
\end{align*}

By~\eqref{cut}, $\chi(\v\eta)\chi(\v\xi)=0$, so all terms containing both
boundary layers vanish in the expansion of $(T^a)^4$. The lower and upper quartic remainders are therefore given by
\begin{align}\label{qu}
 \mathcal{R}_N
 &=\chi(\v\eta)\sum_{k=0}^N\v^k\mathcal P(\mathcal T,k)
               -\mathfrak F(U^N,\chi(\v\eta)V^N),\nonumber\\
 \mathcal{R}_N^+
 &=\chi(\v\xi)\sum_{k=1}^N\v^k\mathcal P^+(\mathcal T^+,k)
               -\mathfrak F(U^N,\chi(\v\xi)V_+^N).
\end{align}

\subsection{Smallness of the residual sources}

\begin{proposition}[Residual source estimate]\label{src}
Let the profiles be those constructed in Proposition~\ref{prop},
and define $(R_1,R_2)$ by~\eqref{rd}.
Then there exist constants $c_N>0$ and $C_N(\kappa)$, independent of
$\v,\rho$, such that, for $0<\v\leq1$ and $0<\rho<1$,
\begin{align}\label{7.0}
 \|R_1\|_\infty+\|R_2\|_\infty
 \leq C_N(\kappa)\left[\v^{N+1}
 +\left(1+\frac{\v}{\rho}+\frac{\v^2}{\rho^2}\right)
                  e^{-c_N(1-\kappa)\rho/\v}\right].
\end{align}
The coefficient may be chosen so that
\begin{align}\label{rpb}
 C_N(\kappa)\leq C_N(1-\kappa)^{-(20\cdot4^N+8N+28)}.
\end{align}
Here $c_N$ and $C_N$ are independent of $\kappa,\v$ and $\rho$.
\end{proposition}

\begin{proof}
\emph{Step 1. Profile bounds and interior truncation.}
Fix $q>3$ as in Proposition~\ref{prop}. Choose
$0<4c_N<\min_{0\leq j\leq N,\,0\leq r\leq2}\lambda_{j,r}$ and set
\begin{align}\label{pts}
 M=C_N(1-\kappa)^{-(N+7+5\cdot4^N)},\qquad
 \sigma=4c_N(1-\kappa),
\end{align}
with $C_N\geq1$ chosen below.

Since $s_0\geq3N+5$, the orders $r=N-j$ and $r=1$ satisfy
$0\leq r\leq s_0-3j-2$ for $0\leq j\leq N$.
Taking these orders in the second line of~\eqref{6.11} for $j=0,1$
and in~\eqref{per} for $2\leq j\leq N$, we obtain, for $0\leq j\leq N$,
\begin{align}\label{rib}
 \mathfrak E_{j,N-j}
 &\leq C_{N,q}(1-\kappa)^{-(\frac72N-\frac52j+5\cdot4^j)}
 \leq C_{N,q}(1-\kappa)^{-(N+5\cdot4^N)},\nonumber\\
 \mathfrak E_{j,1}
 &\leq C_{N,q}(1-\kappa)^{-(j+\frac72+5\cdot4^j)}
 \leq C_{N,q}(1-\kappa)^{-(N+\frac72+5\cdot4^N)}.
\end{align}
Applying $W^{m+1,q}(\Omega)\hookrightarrow W^{m,\infty}(\Omega)$
to $T_j$ with $m=N-j+1$ and $m=2$, and to $I_N(\cdot,\beta)$ with $m=1$,
then using the definition~\eqref{pin} and the bounds~\eqref{rib},
we obtain
\begin{align}
 &\|\nabla I_N\|_\infty
 +\sum_{j=0}^N\left(\|\Delta T_j\|_\infty
   +\sum_{\ell=0}^{N-j+1}\|\partial_{x_3}^\ell T_j\|_\infty\right)
 \nonumber\\
 &\quad\leq C_{N,q}\left[\mathfrak E_{N,1}
        +\sum_{j=0}^N(\mathfrak E_{j,1}+\mathfrak E_{j,N-j})\right]
 \nonumber\\
 &\quad\leq C_{N,q}(1-\kappa)^{-(N+\frac72+5\cdot4^N)}\leq M.
 \label{isb}
\end{align}

For $2\leq j\leq N$, applying the tangential estimates
\eqref{par} and~\eqref{uqb} with $s=2$, and the normal derivative estimates
\eqref{nlb} and~\eqref{unb} with $k=j$, we obtain
\begin{align*}
 &\mathfrak A_{j,2}+\mathfrak A^+_{j,2}
 +\|e^{\sigma\eta}\partial_\eta\mathcal T_j\|_\infty
 +\|e^{\sigma\xi}\partial_\xi\mathcal T_j^+\|_\infty
 \leq C_j(1-\kappa)^{-(j+7+5\cdot4^j)},
\end{align*}
where $\frac72j+32\leq15\cdot4^{j-1}$ for $j\geq2$.
For $j=0$, \eqref{pzr} with $s=2$ and the first line of~\eqref{nlb}
bound the same sum by $C(1-\kappa)^{-6}$, since the upper layer vanishes.
For $j=1$, \eqref{ltr}--\eqref{6.8} and~\eqref{uqb} with $s=2$,
together with the $k=1$ cases of~\eqref{nlb} and~\eqref{unb},
bound the sum by $C(1-\kappa)^{-16}$.
By~\eqref{pts}, $\sigma<\lambda_{j,r}(1-\kappa)$ for $0\leq r\leq2$.
Thus the definitions~\eqref{pwn}--\eqref{upn} bound the weighted tangential
terms below by $\mathfrak A_{j,2}+\mathfrak A^+_{j,2}$.
Summing the tangential and normal derivative bounds over $j$ and choosing
$C_N$ in~\eqref{pts} sufficiently large, we obtain
\begin{align}\label{bls}
 \sum_{j=0}^N\Bigl[&
 \|e^{\sigma\eta}(\mathcal T_j,\partial_\eta\mathcal T_j,
            \Delta_\shortparallel\mathcal T_j,
            \mathcal I_j,\nabla_\shortparallel\mathcal I_j)\|_\infty
 \nonumber\\
 &+\|e^{\sigma\xi}(\mathcal T_j^+,\partial_\xi\mathcal T_j^+,
            \Delta_\shortparallel\mathcal T_j^+,
            \mathcal I_j^+,\nabla_\shortparallel\mathcal I_j^+)\|_\infty
 \Bigr]\leq M.
\end{align}

Applying~\eqref{isb} to the differential terms in~\eqref{fr}
and expanding the powers of $\v$ in $(U^N)^4$, we obtain, for $0<\v\leq1$,
\begin{align}\label{itb}
 &\|\v^{N+1}\beta\cdot\nabla I_N\|_\infty
  +\|\v^{N+1}\Delta T_{N-1}+\v^{N+2}\Delta T_N\|_\infty
  \leq C_N M\v^{N+1},\nonumber\\
 &\left\|\sum_{k=N+1}^{4N}\v^k\mathcal C(T,k)\right\|_\infty
  \leq C_N M^4\v^{N+1}.
\end{align}
Indeed, each product in the second line has four factors from
$T_0,\ldots,T_N$ and carries a power of $\v$ of order at least $N+1$.

\medskip\noindent{\it Step 2. Quartic remainders.}
One has
\begin{align}\label{lqs}
 \mathcal R_N
 =&\chi\left[\sum_{k=0}^N\v^k\mathcal P(\mathcal T,k) -\mathfrak F(Q^N,V^N)\right]+\chi\left[\mathfrak F(Q^N,V^N)-\mathfrak F(U^N,V^N)\right]\nonumber\\
 &+\chi\mathfrak F(U^N,V^N)-\mathfrak F(U^N,\chi V^N).
\end{align}
For $\mathcal R_N^+$, the second term of~\eqref{qu} gives the same
decomposition with $\xi,Q_+^N,V_+^N,\mathcal P^+(\mathcal T^+,k)$
in place of $\eta,Q^N,V^N,\mathcal P(\mathcal T,k)$.
The sum starts at $k=1$ since $\mathcal T_0^+=0$.
The factors $(-\xi)^\ell$ in~\eqref{uq} have absolute value $\xi^\ell$,
so the following three estimates apply at both boundaries.

By~\eqref{3.15-0} (and~\eqref{up} at the upper boundary), the coefficients
through order $N$ cancel in the first bracket of~\eqref{lqs}.
Each remaining product contains at least one layer factor, since
$\mathfrak F(Q^N,V^N)=(Q^N+V^N)^4-(Q^N)^4$.
By~\eqref{3.12-1}, \eqref{isb}, and~\eqref{bls},
$|Q_j|\leq C_NM(1+\eta)^j$ and $|\mathcal T_j|\leq M e^{-\sigma\eta}$
for $0\leq j\leq N$. Hence
\begin{align}\label{qtb}
 \left|\chi\left[\sum_{k=0}^N\v^k\mathcal P(\mathcal T,k)-\mathfrak F(Q^N,V^N)\right]\right|
 \leq C_N M^4\v^{N+1}(1+\eta)^{4N}e^{-\sigma\eta}.
\end{align}

For the second term in \eqref{lqs}, note that $\v\eta\leq3\rho/8<1$
on the support of $\chi$, and each monomial of $Q^N$ contains a factor
$\v^j(\v\eta)^\ell$. Applying the Taylor formula~\eqref{tr}
(and~\eqref{utr} at the upper boundary), together with
the bounds~\eqref{isb} and~\eqref{bls}, we obtain
\begin{align}\label{utb}
 |U^N|+|Q^N|&\leq C_NM,\qquad |V^N|\leq M e^{-\sigma\eta},\nonumber\\
 |U^N-Q^N|&\leq C_NM\v^{N+1}\sum_{j=0}^N\eta^{N-j+1}.
\end{align}
Factoring $X-Z$ and $Y$ from $\mathfrak F(X,Y)-\mathfrak F(Z,Y)$
gives
\[
 |\mathfrak F(X,Y)-\mathfrak F(Z,Y)|
 \leq C|X-Z|\,|Y|\,\bigl(|X|+|Z|+|Y|\bigr)^2.
\]
Taking $(X,Z,Y)=(Q^N,U^N,V^N)$ and using~\eqref{utb}, we obtain
\begin{align}\label{qtf}
 \left|\chi\left[\mathfrak F(Q^N,V^N)-\mathfrak F(U^N,V^N)\right]\right|
 \leq C_N M^4\v^{N+1}e^{-\sigma\eta}\sum_{j=0}^N\eta^{N-j+1}.
\end{align}

For the third term of~\eqref{lqs}, expansion of the fourth powers gives
\[
 \chi\mathfrak F(X,Y)-\mathfrak F(X,\chi Y)
 =6(\chi-\chi^2)X^2Y^2+4(\chi-\chi^3)XY^3+(\chi-\chi^4)Y^4.
\]
Each term contains at least two factors of $Y$. By~\eqref{cut},
$\chi-\chi^j$ for $j=2,3,4$ is supported where
$\rho/4\leq\v\eta\leq3\rho/8$. Taking $(X,Y)=(U^N,V^N)$
and using~\eqref{utb} on this support, we obtain
\begin{align}\label{uqc}
 |\chi\mathfrak F(U^N,V^N)-\mathfrak F(U^N,\chi V^N)|
 \leq C_N M^4 e^{-\sigma\rho/(2\v)}.
\end{align}

Substituting~\eqref{qtb}, \eqref{qtf}, and~\eqref{uqc} into~\eqref{lqs}
and its upper-boundary counterpart, we take the supremum in each layer.
Since $N+1\leq4N$ and
$\sup_{z\geq0}(1+z)^{4N}e^{-\sigma z}\leq C_N(1+\sigma^{-4N})$,
we obtain
\begin{align}\label{qlb}
 \|\mathcal R_N\|_\infty+\|\mathcal R_N^+\|_\infty
 \leq C_N M^4\left[(1+\sigma^{-4N})\v^{N+1}
                    +e^{-\sigma\rho/(2\v)}\right].
\end{align}

\medskip\noindent{\it Step 3. Cutoff derivatives and dependence on $\kappa$.}
For $z=\eta$ or $z=\xi$, the chain rule gives
$\partial_z\chi(\v z)=\v\chi'(\v z)$ and
$\partial_z^2\chi(\v z)=\v^2\chi''(\v z)$.
Both derivatives are supported where $z\geq\rho/(4\v)$. By~\eqref{cut},
\begin{align}\label{ucb}
 \sup_{z\geq0}e^{-\sigma z}
 \bigl(|\partial_z\chi(\v z)|+|\partial_z^2\chi(\v z)|\bigr)
 \leq C\left(\frac{\v}{\rho}+\frac{\v^2}{\rho^2}\right)e^{-\sigma\rho/(4\v)}.
\end{align}
Using~\eqref{ucb} and~\eqref{bls}, we bound the sum of the absolute values
of the derivative terms in~\eqref{lr}--\eqref{ur} by
\[
 C_N M\left[\v^{N+1}
 +\left(\frac{\v}{\rho}+\frac{\v^2}{\rho^2}\right)
 e^{-\sigma\rho/(4\v)}\right].
\]

Substituting the interior estimate~\eqref{itb}, the quartic estimate
\eqref{qlb}, and the preceding derivative bounds
into~\eqref{fr}--\eqref{ur}, we obtain
\begin{align}\label{ss}
 &\|R_1\|_\infty+\|R_2\|_\infty\nonumber\\
 &\quad\leq C_N M^4(1+\sigma^{-4N})
 \left[\v^{N+1}
     +\left(1+\frac{\v}{\rho}+\frac{\v^2}{\rho^2}\right)e^{-\sigma\rho/(4\v)}\right].
\end{align}
By~\eqref{pts}, we have
\begin{align*}
 C_N M^4(1+\sigma^{-4N})
 &\leq C_N(1-\kappa)^{-\left(4(N+7+5\cdot4^N)+4N\right)}\\
 &= C_N(1-\kappa)^{-(20\cdot4^N+8N+28)}.
\end{align*}
Since $\sigma/4=c_N(1-\kappa)$, estimate~\eqref{ss} gives
\eqref{7.0} and~\eqref{rpb}.
\end{proof}
\subsection{Linear problem}
Having established Theorems \ref{thm2} and \ref{co}, we have the following lemma.
\begin{lemma}\label{lem7.2}
Let $N>5$ be an integer. Assume \eqref{bb} and the
uniform spectral gap $M\le M_*<1$. Let $a_*$ be defined by
\eqref{4.00}, and let $T^a$ be the composite profile
constructed in Proposition \ref{prop}. Then there are constants
$c_{\rm coh}>0$, $C_N>0$, $\rho_*>0$, $\varepsilon_*>0$, and an exponent
$K_N>0$, all of which may depend on $N$ but are independent of $\kappa$, such that, if
$\rho=\rho_*(1-\kappa)^5$ and
$\v(1-\kappa)^{-K_N}\leq\varepsilon_*$, then
$a_*/2\leq T^a\leq C_N$ and
\begin{align}\label{7.6}
-\int_\Omega4(T^a)^3g\Delta g\,\mathrm dx
\ge c_{\rm coh}\|\nabla g\|_{L^2(\Omega)}^2
-C_N(1-\kappa)^{-K_N}\|g\|_{L^2(\Omega)}^2
\end{align}
for every $g\in H^2(\Omega)$ that is periodic in $x_\shortparallel$ and has
zero traces at $x_3=0$ and $x_3=1$.
\end{lemma}
\begin{proof}
Fix $q>3$ and set
\begin{equation}\label{KN}
 \Lambda_N:=N+\frac72+5\cdot4^N,\qquad
 K_N:=\max\{7,2\Lambda_N\}.
\end{equation}
We first choose $\rho_*>0$ and then $\varepsilon_*>0$
sufficiently small. In particular, we require
$\rho_*<1$ and $\varepsilon_*\leq\rho_*$, so that
\begin{equation}\label{7.21}
 \frac{\v}{\rho}
 =\rho_*^{-1}\v(1-\kappa)^{-5}\leq1.
\end{equation}
Constants depending on $\rho_*$ below are independent of $\kappa,\v$
once $\rho_*$ has been fixed.

\medskip\noindent\emph{Step 1. Uniform bounds for the approximate temperature.}
Taking $s=1$ in the lower tangential estimates~\eqref{pzr},
\eqref{ltr}, \eqref{par} and the upper estimate~\eqref{uqb},
and combining these with the interior bound~\eqref{rib}
and the normal derivative bounds~\eqref{nlb}, \eqref{unb}, we obtain
\begin{align}\label{7.22}
 &\sum_{k=0}^N\left(
 \|T_k\|_{W^{2,q}}
 +\|\mathcal T_k\|_\infty
 +\|\nabla_\shortparallel\mathcal T_k\|_\infty
 +\|\partial_\eta\mathcal T_k\|_\infty\right)\nonumber\\
 &\quad+\sum_{k=1}^N\left(
 \|\mathcal T_k^+\|_\infty
 +\|\nabla_\shortparallel\mathcal T_k^+\|_\infty
 +\|\partial_\xi\mathcal T_k^+\|_\infty\right)
 \leq C_N(1-\kappa)^{-\Lambda_N}.
\end{align}
Here and below, layer profiles are evaluated at
$\eta=x_3/\v$ and $\xi=(1-x_3)/\v$.
Writing $a=T_{0,\infty}$ and
\begin{equation}\label{7.23}
 W:=\sum_{k=1}^N\v^{k-1}
       (T_k+\mathscr T_k+\mathscr T_k^+),
 \qquad \|W\|_\infty\leq C_N(1-\kappa)^{-\Lambda_N},
\end{equation}
we have
\[
 T^a=(1-\chi(x_3))T_0+\chi(x_3)\mathbf T
       +\chi(x_3)(T_0-a)+\v W.
\]
By~\eqref{4.00} and the maximum principle applied to
\eqref{eh}, $a_*\leq T_0,\mathbf T,a\leq a_2$.
The first two terms therefore form a convex combination bounded
below by $a_*$. On the support of $\chi(x_3)$,
\[
 |\chi(x_3)(T_0-a)|
 \leq C\rho\|\partial_{x_3}T_0\|_\infty
 \leq C\rho(1-\kappa)^{-5},
\]
where we used the boundary value~\eqref{Tb} and~\eqref{ell} at $s=0$.
Combining this bound with~\eqref{7.23} gives, throughout $\Omega$,
\[
 T^a\geq a_*-C\rho(1-\kappa)^{-5}
                 -C_N\v(1-\kappa)^{-\Lambda_N}
 \geq a_*-C\rho_*-C_N\varepsilon_*.
\]
Choosing $C\rho_*+C_N\varepsilon_*\leq a_*/2$ also gives
\begin{equation}\label{7.24}
 a_*/2\leq T^a\leq C_N.
\end{equation}

\medskip\noindent\emph{Step 2. Regular derivative terms.}
For $k\geq1$, applying the chain rule to~\eqref{Bdk} and~\eqref{ukc} gives
\begin{align}\label{7.25}
 \partial_{x_3}(\v^k\mathscr T_k)
 &=\v^k\chi'(x_3)\mathcal T_k
   +\v^{k-1}\chi(x_3)\partial_\eta\mathcal T_k,\nonumber\\
 \partial_{x_3}(\v^k\mathscr T_k^+)
 &=-\v^k\chi'(1-x_3)\mathcal T_k^+
   -\v^{k-1}\chi(1-x_3)\partial_\xi\mathcal T_k^+.
\end{align}
Using~\eqref{7.25}, \eqref{7.22}, and~\eqref{7.21},
together with the cutoff derivative bounds following~\eqref{cut}
and $W^{2,q}\hookrightarrow W^{1,\infty}$, we obtain
\begin{equation}\label{7.26}
 \|\nabla(\v W)\|_\infty
 \leq C_N(1-\kappa)^{-\Lambda_N}(1+\v/\rho)
 \leq C_N(1-\kappa)^{-\Lambda_N}.
\end{equation}
Differentiating the leading terms yields
\begin{align}\label{7.27}
 \nabla_\shortparallel T^a
 &=\nabla_\shortparallel T_0
   +\chi(x_3)\nabla_\shortparallel(\mathbf T-a)
   +\nabla_\shortparallel(\v W),\nonumber\\
 \partial_{x_3}T^a-\frac{\chi(x_3)}{\v}\partial_\eta\mathbf T
 &=\partial_{x_3}T_0+\chi'(x_3)(\mathbf T-a)
   +\partial_{x_3}(\v W).
\end{align}
Substituting~\eqref{7.22} and~\eqref{7.26}
into~\eqref{7.27}, and using~\eqref{4.00} and the cutoff derivative
bound $|\chi'|\leq C/\rho$ stated after~\eqref{cut}, we obtain
\begin{equation}\label{7.28}
 \left\|\nabla T^a
       -\frac{\chi(x_3)}{\v}\partial_\eta\mathbf T\,e_3\right\|_\infty
 \leq C_N\bigl((1-\kappa)^{-\Lambda_N}+\rho^{-1}\bigr)
 \leq C_{N,\rho_*}(1-\kappa)^{-\Lambda_N},
\end{equation}
where $e_3=(0,0,1)$ and $\Lambda_N\geq5$.

Integrating by parts on $\Omega$ and substituting~\eqref{7.27},
using the zero Dirichlet traces and tangential periodicity, we obtain
\begin{align}\label{7.29}
 -\int_\Omega 4(T^a)^3 g\Delta g\,dx
 ={}&\int_\Omega 4(T^a)^3|\nabla g|^2\,dx\nonumber\\
 &+\frac{12}{\v}\int_\Omega
       \chi(x_3)(T^a)^2\partial_\eta\mathbf T\,
       g\partial_{x_3}g\,dx\nonumber\\
 &+12\int_\Omega(T^a)^2
       \left(\nabla T^a
       -\frac{\chi(x_3)}{\v}\partial_\eta\mathbf T\,e_3\right)
       \cdot g\nabla g\,dx.
\end{align}
For any $\theta>0$, applying Young's inequality to the last integral
in~\eqref{7.29}, with~\eqref{7.24} and~\eqref{7.28},
gives
\begin{align}\label{7.8}
 &\left|12\int_\Omega(T^a)^2
       \left(\nabla T^a
       -\frac{\chi(x_3)}{\v}\partial_\eta\mathbf T\,e_3\right)
       \cdot g\nabla g\,dx\right|\nonumber\\
 &\quad\leq\theta\|\nabla g\|_2^2
       +C_{N,\rho_*,\theta}(1-\kappa)^{-2\Lambda_N}\|g\|_2^2.
\end{align}

\medskip\noindent\emph{Step 3. The singular normal term.}
The remaining cross term in~\eqref{7.29} is supported in
\[
 \Omega_{\rm in}:=\{0<x_3<3\rho/8\}.
\]
Set $\Omega_{\rm out}:=\{3\rho/8<x_3<1\}$.
We compare $(T^a)^3$ with $\mathbf T^3$ only on
$\Omega_{\rm in}$, where $T_0-a$ is small. The corresponding normal
contribution is
\begin{align*}
 I_\perp
 :={}&\int_{\Omega_{\rm in}}\left[
      4(T^a)^3|\partial_{x_3}g|^2
      +\frac{12}{\v}\chi(x_3)(T^a)^2\partial_\eta\mathbf T\,
          g\partial_{x_3}g\right]dx\\
 ={}&\int_{\Omega_{\rm in}}\left[
      4\mathbf T^3|\partial_{x_3}g|^2
      +\frac{12}{\v}\mathbf T^2\partial_\eta\mathbf T\,
          g\partial_{x_3}g\right]dx\\
 &+4\int_{\Omega_{\rm in}}
       \bigl((T^a)^3-\mathbf T^3\bigr)|\partial_{x_3}g|^2\,dx\\
 &+\frac{12}{\v}\int_{\Omega_{\rm in}}
       \bigl(\chi(x_3)(T^a)^2-\mathbf T^2\bigr)
       \partial_\eta\mathbf T\,g\partial_{x_3}g\,dx.
\end{align*}

To estimate the first integral, set $L=3\rho/(8\v)$.
In integrals over $\eta$, we abbreviate
$g(\v\eta,x_\shortparallel)$ by $g$. The change of variables gives
$dx_3=\v\,d\eta$ and
$\partial_{x_3}g=\v^{-1}\partial_\eta g$.
For fixed $x_\shortparallel$, extend $g$ constantly beyond $L$.
Applying~\eqref{5.0-1} to this extension, by approximation preserving
the zero trace at $\eta=0$, gives
\[
 \int_0^L36\mathbf T|\partial_\eta\mathbf T|^2g^2\,d\eta
 \leq M_*\int_0^L4\mathbf T^3|\partial_\eta g|^2\,d\eta.
\]
Indeed, the derivative of the extension vanishes beyond $L$, and
the omitted part of the left-hand integral is nonnegative.
Young's inequality and $\mathbf T\geq a_*$ therefore yield
\begin{align}\label{5.9}
 &\int_{\Omega_{\rm in}}\left[
      4\mathbf T^3|\partial_{x_3}g|^2
      +\frac{12}{\v}\mathbf T^2\partial_\eta\mathbf T\,
          g\partial_{x_3}g\right]dx\nonumber\\
 &\quad\geq\frac{1}{2\v}\int_{\mathbb T^2}\int_0^L
       \left(4\mathbf T^3|\partial_\eta g|^2
            -36\mathbf T|\partial_\eta\mathbf T|^2g^2\right)
       d\eta\,dx_\shortparallel\nonumber\\
 &\quad\geq2(1-M_*)a_*^3
             \|\partial_{x_3}g\|_{L^2(\Omega_{\rm in})}^2.
\end{align}

For the coefficient difference, the identity
\[
 T^a-\mathbf T=(T_0-a)
       -(1-\chi(x_3))(\mathbf T-a)+\v W
\]
with the temperature bounds~\eqref{7.24} and~\eqref{4.00},
the estimates~\eqref{ell} and~\eqref{7.23},
and the temperature decay in~\eqref{4.0} at $\alpha=\lambda_{0,0}(1-\kappa)$ gives
\begin{align*}
 \|(T^a)^3-\mathbf T^3\|_{L^\infty(\Omega_{\rm in})}
 &\leq C_N\left[
       \rho(1-\kappa)^{-5}
       +\frac{\v}{\rho}(1-\kappa)^{-2}
       +\v(1-\kappa)^{-\Lambda_N}\right]\\
 &\leq C_N\left[\rho_*+(\rho_*^{-1}+1)\varepsilon_*\right].
\end{align*}
Since $(\v/\rho)(1-\kappa)^{-2}=\rho_*^{-1}\v(1-\kappa)^{-7}$
and $K_N\geq\max\{7,\Lambda_N\}$, the hypothesis
$\v(1-\kappa)^{-K_N}\leq\varepsilon_*$ gives the last inequality.
Here $1-\chi(x_3)$ can be nonzero only when $x_3\geq\rho/4$;
the second term follows from
\[
 (1-\kappa)^{-1}
 e^{-\lambda_{0,0}(1-\kappa)\rho/(4\v)}
 \leq C\frac{\v}{\rho}(1-\kappa)^{-2}.
\]
Choosing $\rho_*$ first and then $\varepsilon_*$ so that
\[
 C_N\left[\rho_*+(\rho_*^{-1}+1)\varepsilon_*\right]
 \leq\frac{1-M_*}{8}a_*^3,
\]
we obtain
\begin{align}\label{5.10}
 \left|4\int_{\Omega_{\rm in}}
       \bigl((T^a)^3-\mathbf T^3\bigr)|\partial_{x_3}g|^2\,dx\right|
 \leq\frac{1-M_*}{2}a_*^3
             \|\partial_{x_3}g\|_{L^2(\Omega_{\rm in})}^2.
\end{align}

For the last integral in $I_\perp$, \eqref{7.24}
and~\eqref{4.00} give
\begin{equation}\label{7.10-1}
 |\chi(x_3)(T^a)^2-\mathbf T^2|
 \leq C\bigl(|T_0-a|+\v|W|+1-\chi(x_3)\bigr).
\end{equation}
Using~\eqref{Tb}, the mean value theorem, and~\eqref{ell} at $s=0$, we obtain
\begin{equation}\label{7.10-2}
 |T_0(\v\eta,x_\shortparallel)-a(x_\shortparallel)|
 \leq \v\eta\|\partial_{x_3}T_0\|_\infty
 \leq C\v\eta(1-\kappa)^{-5}.
\end{equation}
Taking $\alpha=\lambda_{0,0}(1-\kappa)$ in~\eqref{4.60}, we have
\begin{equation}\label{7.10-3}
 |\partial_\eta\mathbf T|
 \leq C e^{-\lambda_{0,0}(1-\kappa)\eta}.
\end{equation}
Multiplying~\eqref{7.10-1} by $\v^{-1}|\partial_\eta\mathbf T|$,
we use~\eqref{7.10-2}, \eqref{7.23}, and
$\eta\geq\rho/(4\v)$ on $\{1-\chi(x_3)\ne0\}$ from~\eqref{cut}
for the three terms, respectively. Together with~\eqref{7.10-3}, these give
\begin{align}\label{7.10-4}
 &\frac{1}{\v}
   |\chi(x_3)(T^a)^2-\mathbf T^2|\,|\partial_\eta\mathbf T|\nonumber\\
 &\quad\leq C_N\left[
       (1-\kappa)^{-5}\eta
                  e^{-\lambda_{0,0}(1-\kappa)\eta}
       +(1-\kappa)^{-\Lambda_N}
       +\v^{-1}e^{-\lambda_{0,0}(1-\kappa)\rho/(4\v)}\right]\nonumber\\
 &\quad\leq C_N\left[
       (1-\kappa)^{-6}
       +(1-\kappa)^{-\Lambda_N}
       +\frac{1}{(1-\kappa)\rho}\right]
 \leq C_{N,\rho_*}(1-\kappa)^{-\Lambda_N}.
\end{align}
The second inequality uses $ze^{-z}\leq e^{-1}$ for $z\geq0$,
and the last uses $\rho=\rho_*(1-\kappa)^5$ and $\Lambda_N\geq6$.
Applying Young's inequality with~\eqref{7.10-4} gives
\begin{align}\label{7.11}
 &\left|\frac{12}{\v}\int_{\Omega_{\rm in}}
       \bigl(\chi(x_3)(T^a)^2-\mathbf T^2\bigr)
       \partial_\eta\mathbf T\,g\partial_{x_3}g\,dx\right|\nonumber\\
 &\quad\leq\frac{1-M_*}{2}a_*^3
             \|\partial_{x_3}g\|_{L^2(\Omega_{\rm in})}^2
       +C_{N,\rho_*}(1-\kappa)^{-2\Lambda_N}
             \|g\|_{L^2(\Omega_{\rm in})}^2.
\end{align}
Combining~\eqref{5.9}--\eqref{7.11}, we obtain
\begin{align}\label{7.12}
 I_\perp\geq(1-M_*)a_*^3
             \|\partial_{x_3}g\|_{L^2(\Omega_{\rm in})}^2
       -C_{N,\rho_*}(1-\kappa)^{-2\Lambda_N}\|g\|_2^2.
\end{align}

By~\eqref{7.24}, the remaining principal terms satisfy
\begin{equation}\label{cre}
 \int_\Omega 4(T^a)^3|\nabla_\shortparallel g|^2\,dx
 +\int_{\Omega_{\rm out}}4(T^a)^3|\partial_{x_3}g|^2\,dx
 \geq\frac{a_*^3}{2}\left(
       \|\nabla_\shortparallel g\|_2^2
       +\|\partial_{x_3}g\|_{L^2(\Omega_{\rm out})}^2\right).
\end{equation}
Taking
\[
 \theta=\min\{a_*^3/4,(1-M_*)a_*^3/2\}
\]
in~\eqref{7.8} and combining~\eqref{7.29}, \eqref{7.12},
and~\eqref{cre}, we conclude
\[
 -\int_\Omega4(T^a)^3g\Delta g\,dx
 \geq\theta\|\nabla g\|_2^2
       -C_{N,\rho_*}(1-\kappa)^{-2\Lambda_N}\|g\|_2^2.
\]
Since $2\Lambda_N\leq K_N$ and $\rho_*$ is now fixed, this proves
\eqref{7.6} with
$c_{\rm coh}=\min\{a_*^3/4,(1-M_*)a_*^3/2\}$.
\end{proof}

We now consider the linearized system
\begin{align}\label{2.3}
 \begin{cases}
 \v\beta\cdot\nabla I_R^\v+\kappa(I_R^\v-\bar I_R^\v)
 =(1-\kappa)(4(T^a)^3T_R^\v-I_R^\v)+R_1,\\
 -\v^2\Delta T_R^\v
 =(1-\kappa)4\pi(\bar I_R^\v-4(T^a)^3T_R^\v)+R_2,
 \end{cases}
\end{align}
with
\[
 I_R^\v|_{\gamma_-}=0,\qquad T_R^\v|_{\partial\Omega}=0.
\]
Here $R_1,R_2$ are general source terms.

For $S\in L^\infty(\Omega\times\mathbb S^2)$, the characteristic
argument used in the proof of Lemma~\ref{lem4.2} gives a unique
bounded solution of
\begin{align}\label{2.4}
 \v\beta\cdot\nabla\phi+\phi-\kappa\bar\phi=S,
 \qquad \phi|_{\gamma_-}=0,
\end{align}
satisfying
\begin{align}\label{2.8}
 \max\{\|\phi\|_{L^\infty},|\phi|_{L^\infty(\gamma_+)}\}
 \leq\frac{1}{1-\kappa}\|S\|_{L^\infty}.
\end{align}

\begin{lemma}\label{lem2.1}
Let $N>5$ be an integer and assume the hypotheses and the smallness
condition of Lemma~\ref{lem7.2}, with $\varepsilon_*$ sufficiently small.
For $R_1\in L^\infty(\Omega\times\mathbb S^2)$ and $R_2\in L^2(\Omega)$,
system~\eqref{2.3} with $I_R^\v|_{\gamma_-}=0$ and
$T_R^\v|_{\partial\Omega}=0$ has a unique solution with
$T_R^\v\in H^2(\Omega)\cap H^1_0(\Omega)$, periodic in $x_\shortparallel$,
and $I_R^\v\in L^\infty(\Omega\times\mathbb S^2)$. Moreover,
\begin{align}
 \|T_R^\v\|_{L^\infty}+\|I_R^\v\|_{L^\infty}
 +|I_R^\v|_{L^\infty(\gamma_+)}
 \leq\frac{C_N}{\v^3}(1-\kappa)^{-K_N/2}
 \bigl(\|R_1\|_{L^\infty}+\|R_2\|_{L^2_x}\bigr),
\end{align}
where $C_N$ is independent of $\v$ and $\kappa$.
\end{lemma}
\begin{proof}
Lemma~\ref{lem7.2} gives $a_*^3/2\leq4(T^a)^3\leq C_N$.
For fixed $\v,\kappa$, the decoupled construction used in the proof
of Lemma~\ref{lem5.1} gives an affine temperature map on $L^2(\Omega)$
whose linear part is compact, by the elliptic estimate and
$H^2(\Omega)\Subset L^2(\Omega)$.
The Fredholm alternative reduces existence to showing that
$R_1=R_2=0$ implies $T_R^\v=I_R^\v=0$.

\medskip\noindent\emph{Step 1. Energy estimate.}
Multiplying~\eqref{2.3}$_1$ by $I_R^\v$ and integrating over
$\Omega\times\mathbb S^2$, we obtain
\begin{align}\label{2.30}
 &\frac{\v}{2}|I_R^\v|_{L^2(\gamma_+)}^2
 +\kappa\|I_R^\v-\bar I_R^\v\|_2^2\nonumber\\
 &\qquad=(1-\kappa)\iint(4(T^a)^3T_R^\v-I_R^\v)I_R^\v
 +\iint R_1I_R^\v.
\end{align}
Multiplying~\eqref{2.3}$_2$ by $4(T^a)^3T_R^\v$,
integrating over $\Omega$, and applying~\eqref{7.6}, we obtain
\begin{align}\label{2.31}
 &c_{\rm coh}\v^2\|\nabla T_R^\v\|_2^2
 +(1-\kappa)\iint4(T^a)^3T_R^\v(4(T^a)^3T_R^\v-I_R^\v)\nonumber\\
 &\qquad\leq C_N(1-\kappa)^{-K_N}\v^2\|T_R^\v\|_2^2
 +4\int_\Omega R_2(T^a)^3T_R^\v.
\end{align}
Adding~\eqref{2.30} and~\eqref{2.31}, we obtain
\begin{align}\label{2.32}
 &\frac{\v}{2}|I_R^\v|_{L^2(\gamma_+)}^2
 +c_{\rm coh}\v^2\|\nabla T_R^\v\|_2^2
 +\kappa\|I_R^\v-\bar I_R^\v\|_2^2+(1-\kappa)\|4(T^a)^3T_R^\v-I_R^\v\|_2^2\nonumber\\
 &\quad\leq C_N(1-\kappa)^{-K_N}\v^2\|T_R^\v\|_2^2
 +\iint R_1I_R^\v+4\int_\Omega R_2(T^a)^3T_R^\v.
\end{align}

\medskip\noindent\emph{Step 2. The $L^2$ estimate for $T_R^\v$.}
Let $l$ be the periodic solution of
\begin{equation}\label{lpt}
 -\Delta l=T_R^\v\quad\text{in }\Omega,\qquad
 l=0\quad\text{on }\partial\Omega.
\end{equation}
The elliptic and trace estimates give
\begin{align}\label{ltb}
 \|l\|_{H^2}+\|\nabla l\|_{L^2(\partial\Omega)}
 &\leq C\|T_R^\v\|_2,\nonumber\\
 \|l\|_{H^3}+\|D^2l\|_{L^2(\partial\Omega)}
 &\leq C\|T_R^\v\|_{H^1}.
\end{align}
Multiplying~\eqref{2.3}$_1$ by $l$, integrating over
$\Omega\times\mathbb S^2$, and using $l|_{\partial\Omega}=0$, we obtain
\begin{align}\label{2.34}
 (1-\kappa)\int_\Omega\langle I_R^\v-4(T^a)^3T_R^\v\rangle l
 =\v\iint I_R^\v\,\beta\cdot\nabla l+\iint R_1l.
\end{align}
Testing~\eqref{2.3}$_1$ by $\beta\cdot\nabla l$ and integrating
by parts in $x$, we obtain, since $\langle\beta\rangle=0$,
\begin{align}\label{2.35}
 \iint I_R^\v\,\beta\cdot\nabla l
 &=\v\iint I_R^\v(\beta\cdot\nabla)^2l
 -\v\int_{\gamma_+}(\beta\cdot n)I_R^\v\,\beta\cdot\nabla l\nonumber\\
 &\quad+\iint R_1\beta\cdot\nabla l.
\end{align}
Testing~\eqref{2.3}$_1$ by $(\beta\cdot\nabla)^2l$,
integrating by parts in $x$, and substituting the resulting identity
and~\eqref{2.35} into~\eqref{2.34}, we obtain
\begin{align}\label{2.36}
 &(1-\kappa)\int_\Omega\langle I_R^\v-4(T^a)^3T_R^\v\rangle l\nonumber\\
 &\quad=\v^3\iint(I_R^\v-\bar I_R^\v)(\beta\cdot\nabla)^3l
 -\v^3\int_{\gamma_+}(\beta\cdot n)I_R^\v(\beta\cdot\nabla)^2l\nonumber\\
 &\qquad+\frac{4\pi}{3}\v^2
 \int_\Omega\bigl(\kappa\bar I_R^\v+(1-\kappa)4(T^a)^3T_R^\v\bigr)\Delta l\nonumber\\
 &\qquad-\v^2\int_{\gamma_+}(\beta\cdot n)I_R^\v\,\beta\cdot\nabla l
 +\v^2\iint R_1(\beta\cdot\nabla)^2l\nonumber\\
 &\qquad+\v\iint R_1\beta\cdot\nabla l+\iint R_1l.
\end{align}
Here we used $I_R^\v|_{\gamma_-}=0$ and
\[
 \int_{\mathbb S^2}\beta_i\beta_j\,\mathrm d\beta
 =\frac{4\pi}{3}\delta_{ij},\qquad
 \int_{\mathbb S^2}\beta_i\beta_j\beta_k\,\mathrm d\beta=0.
\]
Testing~\eqref{2.3}$_2$ by $l$ and using~\eqref{lpt},
we obtain
\begin{align}\label{2.38}
 \v^2\|T_R^\v\|_2^2
 +(1-\kappa)\int_\Omega\langle4(T^a)^3T_R^\v-I_R^\v\rangle l
 =\int_\Omega R_2l.
\end{align}
Testing~\eqref{2.3}$_2$ by $T_R^\v$, we also obtain
\begin{align}\label{2.39-0}
 &\v^2\|\nabla T_R^\v\|_2^2
 +16\pi(1-\kappa)\int_\Omega(T^a)^3(T_R^\v)^2\nonumber\\
 &\qquad=4\pi(1-\kappa)\int_\Omega\bar I_R^\v T_R^\v
 +\int_\Omega R_2T_R^\v.
\end{align}

Substituting $\Delta l=-T_R^\v$ into~\eqref{2.36} and
using~\eqref{2.39-0} to eliminate
$\int_\Omega\bar I_R^\v T_R^\v$, we obtain
\begin{align}\label{2.39}
 &(1-\kappa)\int_\Omega\langle4(T^a)^3T_R^\v-I_R^\v\rangle l\nonumber\\
 &\quad=\frac{16\pi}{3}\v^2\int_\Omega(T^a)^3(T_R^\v)^2
 +\frac{\kappa\v^4}{3(1-\kappa)}\|\nabla T_R^\v\|_2^2\nonumber\\
 &\qquad-\frac{\kappa\v^2}{3(1-\kappa)}\int_\Omega R_2T_R^\v
 -\v^3\iint(I_R^\v-\bar I_R^\v)(\beta\cdot\nabla)^3l\nonumber\\
 &\qquad+\v^3\int_{\gamma_+}(\beta\cdot n)I_R^\v(\beta\cdot\nabla)^2l
 +\v^2\int_{\gamma_+}(\beta\cdot n)I_R^\v\,\beta\cdot\nabla l\nonumber\\
 &\qquad-\v^2\iint R_1(\beta\cdot\nabla)^2l
 -\v\iint R_1\beta\cdot\nabla l-\iint R_1l.
\end{align}
Substituting~\eqref{2.39} into~\eqref{2.38} and
applying~\eqref{ltb}, we obtain
\begin{align}\label{2.40}
 &\v^2\|T_R^\v\|_2^2
 +\frac{16\pi}{3}\v^2\int_\Omega(T^a)^3(T_R^\v)^2
 +\frac{\kappa\v^4}{3(1-\kappa)}\|\nabla T_R^\v\|_2^2\nonumber\\
 &\quad\leq C\v^3\|I_R^\v-\bar I_R^\v\|_2\|T_R^\v\|_{H^1}
 +C\v^3|I_R^\v|_{L^2(\gamma_+)}\|T_R^\v\|_{H^1}\nonumber\\
 &\qquad+C\v^2|I_R^\v|_{L^2(\gamma_+)}\|T_R^\v\|_2
 +C\bigl(\|R_1\|_2+\|R_2\|_2\bigr)\|T_R^\v\|_2\nonumber\\
 &\qquad+\frac{C\v^2}{1-\kappa}\|R_2\|_2\|T_R^\v\|_2.
\end{align}
The smallness condition in Lemma~\ref{lem7.2}, with
$\varepsilon_*\leq1$ and $K_N\geq7$, gives
\begin{equation}\label{lsm}
 \v^2\leq\varepsilon_*^2(1-\kappa)^{2K_N}\leq1-\kappa.
\end{equation}
By Young's inequality and~\eqref{lsm},
\[
 \frac{C\v^2}{1-\kappa}\|R_2\|_2\|T_R^\v\|_2
 \leq\frac{\v^2}{8}\|T_R^\v\|_2^2
 +\frac{C\v^2}{(1-\kappa)^2}\|R_2\|_2^2
 \leq\frac{\v^2}{8}\|T_R^\v\|_2^2
 +\frac{C}{\v^2}\|R_2\|_2^2.
\]
Applying Young's inequality to the remaining terms
in~\eqref{2.40} and absorbing the resulting multiples of
$\v^2\|T_R^\v\|_2^2$, we obtain
\begin{align}\label{2.43}
 \v^2\|T_R^\v\|_2^2
 &\leq C\v^3\|I_R^\v-\bar I_R^\v\|_2^2
 +C\v^3\|\nabla T_R^\v\|_2^2
 +C\v^2|I_R^\v|_{L^2(\gamma_+)}^2\nonumber\\
 &\quad+\frac{C}{\v^2}
 \bigl(\|R_1\|_2^2+\|R_2\|_2^2\bigr).
\end{align}
Since $\langle I_R^\v-\bar I_R^\v\rangle=0$, we have
\begin{align}\label{las}
 \|4(T^a)^3T_R^\v-I_R^\v\|_2^2
 &=4\pi\|4(T^a)^3T_R^\v-\bar I_R^\v\|_{L^2_x}^2\nonumber\\
 &\quad+\|I_R^\v-\bar I_R^\v\|_2^2.
\end{align}
Writing $I_R^\v=(I_R^\v-4(T^a)^3T_R^\v)+4(T^a)^3T_R^\v$
and applying Young's inequality with $4(T^a)^3\leq C_N$, we obtain
\begin{align}\label{2.41}
 &\left|\iint R_1I_R^\v+4\int_\Omega R_2(T^a)^3T_R^\v\right|\nonumber\\
 &\quad\leq\frac{\v^2}{2}\|I_R^\v-4(T^a)^3T_R^\v\|_2^2
 +\frac{\v^2}{2}\|T_R^\v\|_2^2
 +\frac{C_N}{\v^2}\bigl(\|R_1\|_2^2+\|R_2\|_2^2\bigr).
\end{align}
Substituting~\eqref{2.41} into~\eqref{2.32}, we use
$\v^2\leq1-\kappa$ from~\eqref{lsm} to absorb the first term on the right of~\eqref{2.41} into half of
$(1-\kappa)\|4(T^a)^3T_R^\v-I_R^\v\|_2^2$.
The two temperature terms combine into
$C_N(1-\kappa)^{-K_N}\v^2\|T_R^\v\|_2^2$, since
$C_N(1-\kappa)^{-K_N}+\tfrac12
\leq(C_N+\tfrac12)(1-\kappa)^{-K_N}$.
Multiplying~\eqref{2.43} by $(1-\kappa)^{-K_N}$ and
using~\eqref{las} with
$\kappa+(1-\kappa)/2\geq1/2$, we obtain
\begin{align}\label{lte}
 &(1-\kappa)^{-K_N}\v^2\|T_R^\v\|_2^2\nonumber\\
 &\quad\leq C_N\v(1-\kappa)^{-K_N}
 \Bigl[\frac{\v}{2}|I_R^\v|_{L^2(\gamma_+)}^2
 +c_{\rm coh}\v^2\|\nabla T_R^\v\|_2^2\nonumber\\
 &\qquad\qquad+\kappa\|I_R^\v-\bar I_R^\v\|_2^2
 +\frac{1-\kappa}{2}\|4(T^a)^3T_R^\v-I_R^\v\|_2^2\Bigr]\nonumber\\
 &\qquad+\frac{C_N}{\v^2}(1-\kappa)^{-K_N}
 \bigl(\|R_1\|_2^2+\|R_2\|_2^2\bigr).
\end{align}
Substituting~\eqref{lte} into~\eqref{2.32}
after the preceding use of~\eqref{2.41}, we absorb the terms in
brackets by choosing $\varepsilon_*$ so that
\[
 C_N\v(1-\kappa)^{-K_N}\leq C_N\varepsilon_*\leq\tfrac12.
\]
Using~\eqref{2.43} once more, we obtain
\begin{align}\label{2.47}
 &\v|I_R^\v|_{L^2(\gamma_+)}^2
 +\v^2\|\nabla T_R^\v\|_2^2+\v^2\|T_R^\v\|_2^2
 +\kappa\|I_R^\v-\bar I_R^\v\|_2^2\nonumber\\
 &\qquad+(1-\kappa)\|4(T^a)^3T_R^\v-I_R^\v\|_2^2\nonumber\\
 &\quad\leq\frac{C_N}{\v^2}(1-\kappa)^{-K_N}
 \bigl(\|R_1\|_2^2+\|R_2\|_2^2\bigr).
\end{align}

\medskip\noindent\emph{Step 3. The $L^\infty$ estimate.}
By~\eqref{2.47} and~\eqref{las},
\begin{equation}\label{lmd}
 (1-\kappa)\|4(T^a)^3T_R^\v-\bar I_R^\v\|_{L^2_x}
 \leq\frac{C_N}{\v}(1-\kappa)^{(1-K_N)/2}
 \bigl(\|R_1\|_2+\|R_2\|_2\bigr).
\end{equation}
Applying the Dirichlet elliptic estimate to~\eqref{2.3}$_2$
and using $H^2(\Omega)\hookrightarrow L^\infty(\Omega)$
and~\eqref{lmd}, we obtain
\begin{align}\label{lti}
 \|T_R^\v\|_\infty
 &\leq\frac{C}{\v^2}\left[
 (1-\kappa)\|4(T^a)^3T_R^\v-\bar I_R^\v\|_{L^2_x}
 +\|R_2\|_2\right]\nonumber\\
 &\leq\frac{C_N}{\v^3}(1-\kappa)^{-K_N/2}
 \bigl(\|R_1\|_2+\|R_2\|_2\bigr).
\end{align}
Applying~\eqref{2.8} to~\eqref{2.3}$_1$ with
$S=(1-\kappa)4(T^a)^3T_R^\v+R_1$
and using~\eqref{lti}, we obtain
\begin{align*}
 \|I_R^\v\|_\infty+|I_R^\v|_{L^\infty(\gamma_+)}
 &\leq C_N\|T_R^\v\|_\infty
 +\frac{2}{1-\kappa}\|R_1\|_\infty\\
 &\leq\frac{C_N}{\v^3}(1-\kappa)^{-K_N/2}
 \bigl(\|R_1\|_\infty+\|R_2\|_2\bigr),
\end{align*}
where we used $\v\leq1$ and $K_N\geq2$.

Taking $R_1=R_2=0$ in~\eqref{2.47} gives $T_R^\v=I_R^\v=0$.
By linearity, this also proves uniqueness.
\end{proof}

\subsection{Nonlinear problem}
Take $\rho=\rho_*(1-\kappa)^5$ as in Lemma~\ref{lem7.2}.
Applying $e^{-z}\leq C_mz^{-m}$ with $m=N+1,N,N-1$ to the three
terms in the exponential prefactor of~\eqref{7.0}, respectively, we obtain
\begin{align}\label{eta}
 \left(1+\frac{\v}{\rho}+\frac{\v^2}{\rho^2}\right)
 e^{-c_N(1-\kappa)\rho/\v}
 \leq C_{N,\rho_*}\v^{N+1}(1-\kappa)^{-6(N+1)}.
\end{align}
Combining~\eqref{7.0}, \eqref{rpb}, and
\eqref{eta}, and using
$\|R_2\|_2\leq|\Omega|^{1/2}\|R_2\|_\infty$, we obtain
\begin{align}\label{rab}
 \|R_1\|_\infty+\|R_2\|_2
 \leq C_{N,\rho_*}\v^{N+1}
       (1-\kappa)^{-(20\cdot4^N+14N+34)}.
\end{align}
Taking $N=6$ in~\eqref{rab}, we obtain
\begin{align*}
 \|R_1\|_\infty+\|R_2\|_{L^2_x}
 \leq C(1-\kappa)^{-(20\cdot4^6+118)}\v^7.
\end{align*}
Set $I_{R,0}^\v=T_{R,0}^\v=0$. For $j\geq0$, define
$(I_{R,j+1}^\v,T_{R,j+1}^\v)$ by
\begin{align}\label{nit}
 &\v\beta\cdot\nabla I_{R,j+1}^\v+I_{R,j+1}^\v
 -\kappa\bar I_{R,j+1}^\v
 -(1-\kappa)4(T^a)^3T_{R,j+1}^\v\nonumber\\
 &\qquad=-R_1+(1-\kappa)
 \bigl[6(T^a)^2(T_{R,j}^\v)^2
       +4T^a(T_{R,j}^\v)^3+(T_{R,j}^\v)^4\bigr],\nonumber\\
 &-\v^2\Delta T_{R,j+1}^\v
 -(1-\kappa)\langle I_{R,j+1}^\v-4(T^a)^3T_{R,j+1}^\v\rangle\nonumber\\
 &\qquad=-R_2-4\pi(1-\kappa)
 \bigl[6(T^a)^2(T_{R,j}^\v)^2
       +4T^a(T_{R,j}^\v)^3+(T_{R,j}^\v)^4\bigr],
\end{align}
with
\[
 I_{R,j+1}^\v|_{\gamma_-}=0,\qquad
 T_{R,j+1}^\v|_{\partial\Omega}=0.
\]
Under the smallness condition of Lemma~\ref{lem7.2},
Lemma~\ref{lem2.1} gives a unique solution at each step.
Taking $j=0$ in~\eqref{nit} and applying
Lemma~\ref{lem2.1} with~\eqref{rab} at $N=6$, we obtain
\begin{align}\label{nif}
 \|T_{R,1}^\v\|_\infty+\|I_{R,1}^\v\|_\infty
 &\leq C_\kappa\v^4,\qquad
 C_\kappa\leq C(1-\kappa)^{-(K_6/2+20\cdot4^6+118)}.
\end{align}
Define
\begin{equation}\label{Khat}
 \widehat K_6:=K_6+20\cdot4^6+118.
\end{equation}
Choose a fixed $\widehat\varepsilon_*>0$ sufficiently small that
\begin{align}\label{ns}
 \v(1-\kappa)^{-\widehat K_6}\leq\widehat\varepsilon_*
 \quad\Longrightarrow\quad
 C\v(1-\kappa)^{-\widehat K_6}\leq\tfrac12,
 \qquad 2C_\kappa\v^4\leq1.
\end{align}
Indeed, by~\eqref{nif},
\[
 2C_\kappa\v^4
 \leq C\widehat\varepsilon_*^4
       (1-\kappa)^{3\widehat K_6+K_6/2}
 \leq C\widehat\varepsilon_*^4\leq1.
\]
Since $\widehat K_6\geq K_6$, choosing
$\widehat\varepsilon_*\leq\varepsilon_*$ also ensures the
smallness condition of Lemma~\ref{lem7.2}.

We prove by induction that
\begin{equation}\label{niu}
 \|T_{R,j}^\v\|_\infty\leq2C_\kappa\v^4,
 \qquad j\geq0.
\end{equation}
The cases $j=0,1$ follow from the initial choice and
\eqref{nif}. If~\eqref{niu} holds at
order $j$, then $\|T_{R,j}^\v\|_\infty\leq1$. Applying
Lemma~\ref{lem2.1} to~\eqref{nit} and using
$T^a\leq C$, we obtain
\begin{align*}
 \|T_{R,j+1}^\v\|_\infty
 &\leq C_\kappa\v^4
 +C\v^{-3}(1-\kappa)^{1-K_6/2}\|T_{R,j}^\v\|_\infty^2\\
 &\leq C_\kappa\v^4
 +C\v(1-\kappa)^{-\widehat K_6}\|T_{R,j}^\v\|_\infty\\
 &\leq C_\kappa\v^4+\tfrac12\|T_{R,j}^\v\|_\infty
 \leq2C_\kappa\v^4,
\end{align*}
where the second and third inequalities use
\eqref{nif}, \eqref{niu}, and~\eqref{ns}.

Subtracting~\eqref{nit} at indices $j$ and $j-1$
and applying Lemma~\ref{lem2.1}, we obtain, for $j\geq1$,
\begin{align}\label{nid}
 &\|T_{R,j+1}^\v-T_{R,j}^\v\|_\infty
 +\|I_{R,j+1}^\v-I_{R,j}^\v\|_\infty\nonumber\\
 &\quad\leq C\v^{-3}(1-\kappa)^{1-K_6/2}
 \bigl(\|T_{R,j}^\v\|_\infty+\|T_{R,j-1}^\v\|_\infty\bigr)
 \|T_{R,j}^\v-T_{R,j-1}^\v\|_\infty\nonumber\\
 &\quad\leq\tfrac12\|T_{R,j}^\v-T_{R,j-1}^\v\|_\infty
 \leq2^{-j}C_\kappa\v^4.
\end{align}
Here the quartic difference is bounded using~\eqref{niu},
and the last two inequalities follow from~\eqref{ns}
and~\eqref{nif}.
By~\eqref{nid}, $(I_{R,j}^\v,T_{R,j}^\v)$ converges in $L^\infty$ to
$(I_R^\v,T_R^\v)$. The Dirichlet estimate applied to
\eqref{nit}$_2$ also gives convergence of
$T_{R,j}^\v$ in $H^2(\Omega)$, and the characteristic formula
preserves the zero incoming data in the limit.
Passing to the limit in~\eqref{nit}, we obtain
\eqref{re}. Applying the same difference estimate to two solutions
with $\|T_R^\v\|_\infty\leq2C_\kappa\v^4$ proves uniqueness in this class.
By~\eqref{niu},
\[
 \|T_R^\v\|_\infty
 \leq2C_\kappa\v^4
 \leq C(1-\kappa)^{-(\widehat K_6-K_6/2)}\v^4.
\]
Applying~\eqref{2.8} to the radiation equation in
\eqref{re}, and using $\|T_R^\v\|_\infty\leq1$,
we obtain
\begin{align*}
 \|T_R^\v\|_\infty+\|I_R^\v\|_\infty
 &\leq C\left[
 (1-\kappa)^{-(\widehat K_6-K_6/2)}\v^4
 +(1-\kappa)^{-(\widehat K_6-K_6+1)}\v^7\right]\\
 &\leq C(1-\kappa)^{-\widehat K_6}\v^4,
\end{align*}
where the second term comes from $(1-\kappa)^{-1}\|R_1\|_\infty$.
The profile bounds in Proposition~\ref{prop} also give
\begin{align*}
 &\sum_{k=1}^6\v^k\left(
 \|(I_k,T_k)\|_\infty
 +\|(\mathscr I_k,\mathscr T_k)\|_\infty
 +\|(\mathscr I_k^+,\mathscr T_k^+)\|_\infty\right)\\
 &\quad\leq C\v(1-\kappa)^{-(6+5\cdot4^6)}
 \leq C\v(1-\kappa)^{-K_6}.
\end{align*}
Taking $\widehat\varepsilon_*\leq1$, condition
\eqref{ns} implies
\[
 (1-\kappa)^{-\widehat K_6}\v^4
 \leq\widehat\varepsilon_*\v^3
 \leq\v(1-\kappa)^{-K_6}.
\]
Adding the remainder and the higher profiles proves~\eqref{1.49}.
For $0\leq\kappa\leq\kappa_0<1$, the error coefficient and the
smallness threshold are uniform when the boundary norms, positive
lower bounds, and spectral constants are uniformly controlled.
This completes the proof of Theorem~\ref{thm1}.

\section{Vanishing-scattering limit}\label{vsl}
\begin{proof}[Proof of Theorem \ref{thm1.5}]
The estimate in Theorem \ref{thm1} is uniform under the hypotheses
of Theorem \ref{thm1.5}. It remains to compare the leading profiles.
All half-line estimates below are taken at a fixed $x_\shortparallel$;
their constants are uniform in $x_\shortparallel$ and
$0\leq\kappa\leq\kappa_0$.

\underline{Step 1. Stability of the boundary layer.}
Set $\phi=\mathbf I-\mathbf I_{\mathrm{ns}}$ and
$g=\mathbf T-\mathbf T_{\mathrm{ns}}$. By~\eqref{sec},
\[
 \mathbf T^4-(\mathbf T_{\mathrm{ns}})^4=a_\kappa g.
\]
The bounds in Theorem \ref{thm2}, with $a_*$ defined in \eqref{4.00},
give $a_*^3\leq a_\kappa\leq4a_2^3$.
Subtracting \eqref{8.1} from \eqref{3.16}, we obtain
\begin{align}\label{8.3}
 \begin{cases}
 \beta_3\partial_\eta\phi+\phi-a_\kappa g=S,\\
 \partial_\eta^2g+\langle\phi-a_\kappa g\rangle=\langle S\rangle,
 \end{cases}
 \qquad S:=\kappa(\bar{\mathbf I}-\mathbf T^4),
\end{align}
where $S$ is independent of $\beta$, so $\langle S\rangle=4\pi S$.
The common boundary data imply
\begin{align}\label{8.4}
 \phi(0,\beta)=0\quad(\beta_3>0),\qquad g(0)=0.
\end{align}
Choose
\[
 \alpha=\frac{\vartheta_*+1-\kappa_0}{2},
 \qquad
 \vartheta_*<\alpha<1-\kappa_0\leq1-\kappa.
\]
By~\eqref{4.0} and~\eqref{4.00}, we have
$|S|\leq C\kappa e^{-\alpha\eta}$, with $C$ independent of $\kappa$,
since $1-\kappa-\alpha\geq(1-\kappa_0-\vartheta_*)/2>0$.
Hence
\begin{align}\label{sb}
 \|e^{\vartheta_*\eta}S\|_{L^\infty}
 &\leq C\kappa,\qquad
 \|e^{\vartheta_*\eta}S\|_{L^2}
 \leq\frac{C\kappa}{\sqrt{\alpha-\vartheta_*}}
 \leq C\kappa.
\end{align}
The $L^2$ norm here is over $\mathbb R_+\times\mathbb S^2$.

Moreover,~\eqref{4.60} gives $e^{\vartheta_*\eta}\partial_\eta g\in L^2(\mathbb R_+)$. 
By the positive bounds on $a_\kappa$ and the spectral assumption on
$a_\kappa^{1/3}$, the $L^\infty$ estimates in the proof of
Theorem~\ref{thm4} apply to~\eqref{8.3}--\eqref{8.4}, with zero
scattering coefficient and $4\mathbf T^3$ replaced by $a_\kappa$.
Together with~\eqref{sb}, they give
\[
 \|g\|_{L^\infty}+\|\phi\|_{L^\infty}
 \leq C\left(
 \|e^{\vartheta_*\eta}S\|_{L^2}
 +\|e^{\vartheta_*\eta}S\|_{L^\infty}\right)
 \leq C\kappa.
\]
Since $\lim_{\eta\to\infty}g=T_{0,\infty}-T_{0,\infty,\mathrm{ns}}$, we obtain
\begin{align}\label{8.5}
 &\|T_{0,\infty}-T_{0,\infty,\mathrm{ns}}\|_{L^\infty(\mathbb T^2)}
 +\|\mathbf I-\mathbf I_{\mathrm{ns}}\|_{L^\infty}
 +\|\mathbf T-\mathbf T_{\mathrm{ns}}\|_{L^\infty}
 \leq C\kappa.
\end{align}

For each $\kappa$, the two truncated layers in \eqref{Bd} and
\eqref{8.2-1} use the same cutoff $\chi$, including the same width.
Their exact differences are
\[
 \begin{aligned}
 \mathscr T_0-\mathscr T_{0,\mathrm{ns}}
 &=\chi(\v\eta)
      \bigl[g-(T_{0,\infty}-T_{0,\infty,\mathrm{ns}})\bigr],\\
 \mathscr I_0-\mathscr I_{0,\mathrm{ns}}
 &=\chi(\v\eta)
      \bigl[\phi-(T_{0,\infty}^4-(T_{0,\infty,\mathrm{ns}})^4)\bigr].
 \end{aligned}
\]
Since $0\leq\chi\leq1$ and
$|T_{0,\infty}^4-(T_{0,\infty,\mathrm{ns}})^4|
 \leq4a_2^3|T_{0,\infty}-T_{0,\infty,\mathrm{ns}}|$,
estimate \eqref{8.5} yields
\begin{align}\label{8.6}
 &\|\mathscr I_0-\mathscr I_{0,\mathrm{ns}}\|_{L^\infty(\Omega\times\mathbb S^2)}
 +\|\mathscr T_0-\mathscr T_{0,\mathrm{ns}}\|_{L^\infty(\Omega)}
 \leq C\kappa.
\end{align}

\underline{Step 2. Stability of the interior profile.}
Recall $\mathcal H(z)=z+\frac{4\pi}{3}z^4$.
Equations~\eqref{3.5}, \eqref{Tb}, and \eqref{8.2-2} give
\[
 \begin{cases}
 \Delta\bigl(\mathcal H(T_0)-\mathcal H(T_{0,\mathrm{ns}})\bigr)=0
   &\text{in }\Omega,\\
 \mathcal H(T_0)-\mathcal H(T_{0,\mathrm{ns}})
 =\mathcal H(T_{0,\infty})-\mathcal H(T_{0,\infty,\mathrm{ns}})
   &\text{at }x_3=0,\\
 \mathcal H(T_0)-\mathcal H(T_{0,\mathrm{ns}})=0
   &\text{at }x_3=1.
 \end{cases}
\]
The maximum principle and~\eqref{8.5} imply
\[
 \|\mathcal H(T_0)-\mathcal H(T_{0,\mathrm{ns}})\|_{L^\infty(\Omega)}
 \leq\|\mathcal H(T_{0,\infty})
             -\mathcal H(T_{0,\infty,\mathrm{ns}})\|_{L^\infty(\mathbb T^2)}
 \leq C\kappa.
\]
Both interior temperatures lie in $[0,a_2]$. Since $\mathcal H'\geq1$
on this interval, its inverse is Lipschitz with constant one.
The quartic map has Lipschitz constant $4a_2^3$ there. Consequently,
\begin{align}\label{8.7}
 \|T_0-T_{0,\mathrm{ns}}\|_{L^\infty(\Omega)}
 +\|T_0^4-(T_{0,\mathrm{ns}})^4\|_{L^\infty(\Omega)}
 \leq C\kappa.
\end{align}

\underline{Step 3. Completion of the proof.}
The uniform estimate in Theorem \ref{thm1} gives
\[
 \|I^\v-T_0^4-\mathscr I_0\|_{L^\infty(\Omega\times\mathbb S^2)}
 +\|T^\v-T_0-\mathscr T_0\|_{L^\infty(\Omega)}
 \leq C\v.
\]
Combining this estimate with \eqref{8.6}--\eqref{8.7} and the
triangle inequality, we obtain
\begin{align*}
 &\|I^\v-(T_{0,\mathrm{ns}})^4-\mathscr I_{0,\mathrm{ns}}
                    \|_{L^\infty(\Omega\times\mathbb S^2)}\\
 &\qquad+\|T^\v-T_{0,\mathrm{ns}}-\mathscr T_{0,\mathrm{ns}}
                    \|_{L^\infty(\Omega)}
 \leq C(\v+\kappa).
\end{align*}
The constant and the smallness threshold for $\v$ are independent
of $\kappa\in[0,\kappa_0]$, as asserted.
\end{proof}

\appendix
\section{Proof of Lemma~\ref{lem5.1}}\label{app:lem5.1}
\begin{proof}
We write $\mathcal R\beta=(\beta_{\shortparallel},-\beta_3)$ for the reflected direction. For a given $g^B$, the solution to \eqref{5.2} can be expressed along backward characteristics as
\begin{align}\label{5.5}
	\phi^B(\eta,\beta)=
		p_b(\beta)e^{-\f{\eta}{\beta_3}}+\int_0^\eta \f{1}{\beta_3}e^{-\f{\eta-z}{\beta_3}}\left[\kappa\bar{\phi}^B(z)+(1-\kappa)4\mathbf{T}^3(z)g^B(z)+S(z,\beta)\right]\,\mathrm{dz},\quad \text{for }\beta_3>0,
	\end{align}
and
\begin{align}\label{5.5-1}
	\phi^B(\eta,\beta)&=	p_b(\mathcal{R}\beta)e^{-\f{2B-\eta}{|\beta_3|}}+\int_{\eta}^B\f{1}{|\beta_3|}e^{-\f{z-\eta}{|\beta_3|}}\left[\kappa\bar{\phi}^B(z)+(1-\kappa)4\mathbf{T}^3(z)g^B(z)+S(z,\beta)\right]\,\mathrm{dz}\nonumber\\
	&\quad +\int_0^B\f{1}{|\beta_3|}e^{-\f{2B-\eta-z}{|\beta_3|}}\left[\kappa\bar{\phi}^B(z)+(1-\kappa)4\mathbf{T}^3(z)g^B(z)+S(z,\mathcal R\beta)\right]\,\mathrm{dz}, \quad \text{for }\beta_3<0.
\end{align}
The last integral in \eqref{5.5-1} follows the characteristic before
reflection at $\eta=B$, where its direction is $\mathcal R\beta$.
Integrating \eqref{5.5} and \eqref{5.5-1} over the two hemispheres, we obtain
\begin{align}\label{5.6}
	\langle \phi^B(\eta,\cdot)\rangle =&\int_{\beta_3>0}\phi^B(\eta,\cdot)\, \mathrm{d}\beta +\int_{\beta_3<0}\phi^B(\eta,\cdot)\, \mathrm{d}\beta\nonumber\\
	=&\int_{\beta_3>0}	p_b(\beta)e^{-\f{\eta}{\beta_3}}\,\mathrm{d}\beta+\int_{\beta_3<0}	p_b(\mathcal{R}\beta)e^{-\f{2B-\eta}{|\beta_3|}}\,\mathrm{d}\beta\nonumber\\
	&+\int_0^B\int_{\beta_3>0}\f{1}{\beta_3}e^{-\f{|\eta-z|}{\beta_3}}\left[\kappa\bar{\phi}^B(z)+(1-\kappa )4\mathbf{T}^3(z)g^B(z)\right]\,\mathrm{d}\beta\,\mathrm{d}z\nonumber\\
	&+\int_0^B\int_{\beta_3>0}\f{1}{\beta_3}e^{-\f{2B-\eta-z}{\beta_3}}\left[\kappa\bar{\phi}^B(z)+(1-\kappa)4\mathbf{T}^3(z)g^B(z)+S(z,\beta)\right]\,\mathrm{d}\beta\,\mathrm{d}z\nonumber\\
	&+\int_0^\eta\int_{\beta_3>0}\f{1}{\beta_3}e^{-\f{\eta-z}{\beta_3}}S(z,\beta)\,\mathrm{d}\beta\,\mathrm{d}z+\int_\eta^B\int_{\beta_3>0}\f{1}{\beta_3}e^{-\f{z-\eta}{\beta_3}}S(z,\mathcal R\beta)\,\mathrm{d}\beta\,\mathrm{d}z.
\end{align}

\textbf{Step 1. Approximate solutions.} To prove the existence for \eqref{5.2}, we construct a sequence $\{\phi_k,g_k\}_{k=1}^\infty$ as follows:
\begin{align}\label{5.7}
	\begin{cases}
	\beta_3\partial_\eta\phi_k+\kappa(\phi_k-\bar{\phi}_k)=(1-\kappa)(4\mathbf{T}^3g_k-\phi_k)+S,\\
	\partial_\eta^2g_k+(1-\kappa)4\pi(\bar{\phi}_{k-1}-4\mathbf{T}^3g_k)=\langle S\rangle,
	\end{cases}
\end{align} 
with boundary conditions 
\begin{align*}
	&\phi_k(0,\beta)|_{\gamma_-}=p_b(\beta),\quad \phi_k(B,\beta)|_{\gamma_-}=\phi_k(B,\mathcal{R}\beta),\\
	&g_k(0)=0,\quad\partial_\eta g_k(B)=0,
\end{align*}
where $k\geq 1$ and we define $g_0(\eta,\beta)\equiv0$. Once $g_{k-1}$ is obtained, it follows from Lemma \ref{lem4.2} that there is a unique solution $\phi_{k-1}$. Then \eqref{5.7}$_2$ admits a solution  $g_k$. We denote $g_k=\mathcal{G}g_{k-1}$ for $k\geq 1$.
We now show that $\mathcal{G}$ is a contraction on the weighted space $L^2_{16\mathbf{T}^6}([0,B])$. We first verify that $\mathcal{G}$ maps this space into itself.

Multiplying \eqref{5.7}$_1$ by $\phi_k$ and integrating over $[0,B]\times\mathbb S^2$, we obtain
\begin{align}\label{5.13}
	&\f{1}{2}|\phi_k|_{L^2(\gamma_+)}^2+\|\phi_k-\bar{\phi}_k\|_{L^2}^2+(1-\kappa)4\pi\|\bar{\phi}_k\|_{L^2}^2\nonumber\\
	&=\f12|p_b|_{L^2(\gamma_-)}^2+(1-\kappa)4\pi \int_0^B	4\mathbf{T}^3g_k\bar{\phi}_k\,\mathrm{d}\eta+\int_0^B\int_{\mathbb{S}^2}\phi_k\cdot S\, \mathrm{d}\beta\,\mathrm{d}\eta,
\end{align}
Here the reflecting boundary condition cancels the contribution at $B$, and we used
$\|\phi\|_{L^2}^2=\|\phi-\bar\phi\|_{L^2}^2+4\pi\|\bar\phi\|_{L^2}^2$.

 Multiplying \eqref{5.7}$_2$ by $4\mathbf{T}^3g_k$, and integrating over $[0,B]$, we obtain
\begin{align}\label{ite}
	&-\int_0^B 4\mathbf{T}^3g_k\partial_\eta^2 g_k\,\mathrm{d}\eta+4\pi(1-\kappa)\int_0^B16\mathbf{T}^6g_k^2\,\mathrm{d}\eta\nonumber\\
	&=4\pi(1-\kappa)\int_0^B4\mathbf{T}^3g_k\bar{\phi}_{k-1}\, \mathrm{d}\eta-\int_0^B4\mathbf{T}^3g_k\cdot\langle S\rangle \, \mathrm{d}\eta.
\end{align}
To apply the spectral inequalities on $[0,B]$, extend a function
$f\in H^1(0,B)$ with $f(0)=0$ constantly by $f(B)$ for $\eta>B$.
The extension has zero derivative outside $[0,B]$. Applying
\eqref{asA} and \eqref{5.0-1}, with approximation by smooth functions
for the latter, and discarding the nonnegative exterior term, we obtain,
for $\sigma=\vartheta$ and $\sigma=0$, respectively,
\[
\int_0^B e^{2\sigma\eta}36\mathbf T|\partial_\eta\mathbf T|^2f^2\,\mathrm d\eta
\leq M\int_0^B e^{2\sigma\eta}4\mathbf T^3|\partial_\eta f|^2\,\mathrm d\eta.
\]
Cauchy--Schwarz inequality therefore gives
\[
\left|\int_0^B e^{2\sigma\eta}\partial_\eta(4\mathbf T^3)
 f\partial_\eta f\,\mathrm d\eta\right|
\leq\sqrt M\int_0^B e^{2\sigma\eta}4\mathbf T^3
 |\partial_\eta f|^2\,\mathrm d\eta.
\]
Since $M<1$, this proves the finite-interval versions of
\eqref{5.0-2} and \eqref{5.0-3}. Applying \eqref{5.0-2} to $f=g_k$ and
using $g_k(0)=0$ and $\partial_\eta g_k(B)=0$, we obtain 
\begin{align*}
	-\int_0^B4\mathbf{T}^3g_k\partial_\eta^2 g_k\,\mathrm{d}\eta=\int_0^B4\mathbf{T}^3|\partial_\eta g_k|^2\,\mathrm{d}\eta+\int_0^B(\partial_\eta4\mathbf{T}^3) g_k\partial_\eta g_k\,\mathrm{d}\eta\geq 0.
\end{align*}

Dividing~\eqref{ite} by $4\pi(1-\kappa)$, using the
preceding nonnegativity, and applying Young's inequality
to the two terms on the right, we obtain
\begin{align*}
	\int_0^B16\mathbf T^6g_k^2\,\mathrm d\eta
	&\leq \int_0^B4\mathbf T^3g_k\bar\phi_{k-1}\,\mathrm d\eta
	+\frac{1}{4\pi(1-\kappa)}
	\left|\int_0^B4\mathbf T^3g_k\langle S\rangle\,
	\mathrm d\eta\right|\\
	&\leq \frac34\int_0^B16\mathbf T^6g_k^2\,\mathrm d\eta
	+\frac12\int_0^B|\bar\phi_{k-1}|^2\,\mathrm d\eta
	+\frac{C}{(1-\kappa)^2}
	\int_0^B|\langle S\rangle|^2\,\mathrm d\eta.
\end{align*}
Absorbing the first term on the right gives
\begin{align}\label{5.11-2}
	\int_0^B16\mathbf T^6g_k^2\,\mathrm d\eta
	&\leq 2\int_0^B|\bar\phi_{k-1}|^2\,\mathrm d\eta
	+\frac{C}{(1-\kappa)^2}
	\int_0^B|\langle S\rangle|^2\,\mathrm d\eta.
\end{align}

Applying Young's inequality to \eqref{5.13}, we obtain
\begin{align}\label{5.14}
	\|\bar{\phi}_k\|_{L^2([0,B])}^2\leq \f{C}{1-\kappa}|p_b|_{L^2(\gamma_-)}^2+\f{C}{(1-\kappa)^2}\|S\|_{L^2}^2+C\|4\mathbf{T}^3g_{k}\|_{L^2([0,B])}^2.
\end{align}
Substituting \eqref{5.14}, with index $k-1$, into \eqref{5.11-2}, we deduce that 
\begin{align*}
	\int_0^B16\mathbf{T}^6g_k^2\,\mathrm{d}\eta\leq \f{C}{1-\kappa}|p_b|_{L^2(\gamma_-)}^2+\f{C}{(1-\kappa)^2}\|S\|_{L^2}^2+ C\int_0^B16\mathbf{T}^6g_{k-1}^2\,\mathrm{d}\eta.
\end{align*}
Hence $\mathcal{G}$ maps from $L^2_{(16\mathbf{T}^6)}([0,B])$ into itself.

{\textbf{Step 2. Convergence.}} We now show that $\mathcal{G}$ is a contraction. Let $h_k=g_k-g_{k-1}$ and
$\Phi_k=\phi_k-\phi_{k-1}$. Then the difference satisfies
\begin{align}\label{5.16}
	\begin{cases}
		\beta_3\partial_\eta\Phi_k+\kappa(\Phi_k-\bar{\Phi}_k)=(1-\kappa)(4\mathbf{T}^3h_k-\Phi_k),\\
		-\partial_\eta^2h_k=(1-\kappa)4\pi(\bar{\Phi}_{k-1}-4\mathbf{T}^3h_k)
	\end{cases}
\end{align}
with boundary conditions
\begin{align*}
	&\Phi_k(0,\beta)|_{\gamma_-}=0,\quad \Phi_k(B,\beta)|_{\gamma_-}=\Phi_k(B,\mathcal{R}\beta),\\
	&h_k(0)=0,\quad \partial_{\eta}h_k(B)=0,
\end{align*}
where $k\geq 2$. Multiplying \eqref{5.16}$_2$ by $4\mathbf T^3h_k$ and using the finite-interval version of \eqref{5.0-2} gives 
\begin{align}
	&\int_0^B16\mathbf{T}^6h_k^2\,\mathrm{d}\eta\leq \int_0^B4\mathbf{T}^3h_k\bar{\Phi}_{k-1}\, \mathrm{d}\eta,\label{5.17}\\
	&\f12|\Phi_k(0)|_{L^2(\gamma_+)}^2+\|\Phi_k-\bar{\Phi}_k\|_{L^2}^2+(1-\kappa)4\pi\|\bar{\Phi}_k\|_{L^2}^2=(1-\kappa)4\pi\int_0^B4\mathbf{T}^3h_k\bar{\Phi}_k\,\mathrm{d}\eta.\label{5.18}
\end{align}
We have from \eqref{5.18} that 
\begin{align}\label{5.19}
	\|\bar{\Phi}_k\|_{L^2}^2\leq \int_0^B4\mathbf{T}^3h_k\bar{\Phi}_k\, \mathrm{d}\eta\leq \int_0^B16\mathbf{T}^6h_k^2(\eta)\, \mathrm{d}\eta,
\end{align}
where we used the {\it Holder} inequality.

For $k\geq3$, applying~\eqref{5.6} with zero source and
incoming data to~\eqref{5.16}$_1$ at index $k-1$,
and substituting into~\eqref{5.17}, we obtain
\begin{align*}
	&4\pi\int_0^B16\mathbf T^6h_k^2\,\mathrm d\eta\\
	&\quad\leq
	\int_0^B\int_0^B\int_{\beta_3>0}
	\frac{e^{-|\eta-z|/\beta_3}
		+e^{-(2B-\eta-z)/\beta_3}}{\beta_3}
	4\mathbf T^3(\eta)h_k(\eta)\\
	&\qquad\qquad\times
	\bigl[\kappa\bar\Phi_{k-1}(z)
	+(1-\kappa)4\mathbf T^3(z)h_{k-1}(z)\bigr]\,
	\mathrm d\beta\,\mathrm dz\,\mathrm d\eta\\
	&\quad\leq
	\frac12\int_0^B\int_0^B\int_{\beta_3>0}
	\frac{e^{-|\eta-z|/\beta_3}
		+e^{-(2B-\eta-z)/\beta_3}}{\beta_3}\\
	&\qquad\qquad\times
	\bigl[16\mathbf T^6(\eta)h_k^2(\eta)
	+\kappa\bar\Phi_{k-1}^2(z)
	+(1-\kappa)16\mathbf T^6(z)h_{k-1}^2(z)\bigr]\,
	\mathrm d\beta\,\mathrm dz\,\mathrm d\eta,
\end{align*}
where the last inequality applies $2ab\leq a^2+b^2$
separately to the terms with coefficients $\kappa$ and
$1-\kappa$. For $\beta_3>0$, direct integration gives
\[
\int_0^B
\frac{e^{-|\eta-z|/\beta_3}
	+e^{-(2B-\eta-z)/\beta_3}}{\beta_3}\,\mathrm dz
=2-e^{-\eta/\beta_3}-e^{-(2B-\eta)/\beta_3}.
\]
Using this identity and the symmetry of the kernel in
$\eta,z$, and moving the term containing $h_k^2$
to the left, we obtain
\begin{align}\label{5.20}
	&\f12\int_0^B\int_{\beta_3>0}(2+e^{-\f{\eta}{\beta_3}}+e^{-\f{2B-\eta}{\beta_3}})16\mathbf{T}^6h_k^2\,\mathrm{d}\beta\,\mathrm{d}\eta\nonumber\\
	&\leq \f12\int_0^B\int_{\beta_3>0}(2-e^{-\f{\eta}{\beta_3}}-e^{-\f{2B-\eta}{\beta_3}})\left[\kappa\bar{\Phi}_{k-1}^2(\eta)+(1-\kappa)16\mathbf{T}^6h_{k-1}^2(\eta)\right]\, \mathrm{d}\beta\,\mathrm{d}\eta.
\end{align}

Since for $\beta_3>0$,
\begin{align*}
	\min_{\eta\in[0,B]}(e^{-\f{\eta}{\beta_3}}+e^{-\f{2B-\eta}{\beta_3}})=2e^{-\f{B}{\beta_3}},
\end{align*}
\eqref{5.19} and \eqref{5.20} gives
\begin{align}\label{5.21}
	\int_0^B16\mathbf{T}^6(\eta)h_k^2(\eta)\,\mathrm{d}\eta\leq \d \int_0^B16\mathbf{T}^6(\eta)h_{k-1}^2(\eta)\,\mathrm{d}\eta,\qquad k\geq3,
\end{align}
with
\begin{align*}
	\d=\f{\int_{\beta_3>0}(1-e^{-\f{B}{\beta_3}})\,\mathrm{d}\beta}{\int_{\beta_3>0}(1+e^{-\f{B}{\beta_3}})\,\mathrm{d}\beta}<1.
\end{align*}
Then \eqref{5.21} can be written as 
\begin{align*}
	\|h_k\|_{L^2(16\mathbf{T}^6)([0,B])}^2\leq \d \|h_{k-1}\|_{L^2(16\mathbf{T}^6)([0,B])}^2,
\end{align*}
which is 
\begin{align*}
		\|g_k-g_{k-1}\|_{L^2(16\mathbf{T}^6)([0,B])}^2\leq \d^{k-2}\|g_2-g_1\|_{L^2(16\mathbf{T}^6)([0,B])}^2,\qquad k\geq2.
\end{align*}
The same calculation applies to any two inputs of $\mathcal G$, since their associated transport solutions have the same source and incoming data. Thus $\mathcal{G}:L^2_{16\mathbf{T}^6}([0,B])\to L^2_{16\mathbf{T}^6}([0,B])$ is a contraction with Lipschitz constant $\sqrt{\d}<1$. The contraction factor $\d$ depends on $B$, so this argument constructs the solution for each fixed finite interval. 

Additionally, it follows from \eqref{5.18} and \eqref{5.19} that 
\begin{align*}
	\|\Phi_k-\bar{\Phi}_k\|_{L^2}^2+(1-\kappa)4\pi\|\bar{\Phi}_k\|_{L^2}^2\leq (1-\kappa)4\pi\d^{k-2}\|g_2-g_1\|_{L^2(16\mathbf{T}^6)([0,B])}^2,\quad k\geq2,
\end{align*}
Since the left-hand side is at least $(1-\kappa)\|\Phi_k\|_{L^2}^2$, we obtain
\begin{align*}
	\|\phi_k-\phi_{k-1}\|_{L^2([0,B]\times{\mathbb{S}^2})}^2\leq4\pi\d^{k-2}\|g_2-g_1\|_{L^2(16\mathbf{T}^6)([0,B])}^2,\qquad k\geq2.
\end{align*}
The geometric bound makes $\{\phi_k\}$ a Cauchy sequence in $L^2$, by summing the norms of successive differences. By the Banach fixed point theorem for $\mathcal G$ and passage to the limit in \eqref{5.7}, there exists a unique solution $(\phi^B,g^B)$ to the equation \eqref{5.2}. By \eqref{5.2}$_2$, $\partial_\eta^2g^B\in L^2([0,B])$ and thus $g^B\in C^1([0,B])$. Since $g^B\in L^2_{(16\mathbf{T}^6)}([0,B])$ and $\phi^B\in L^2([0,B])$, by \eqref{5.5} and \eqref{5.5-1}, $\phi^B\in C([0,B]\times \mathbb{S}^2\backslash\gamma_0)$.

\textbf{Step 3. Weighted estimate.} The weighted estimate is derived in a similar manner as for the nonlinear case.  First, we have that
\begin{align}\label{5.21-1}
\partial_\eta g^B(\eta) = \langle \beta_3 \phi^B(\eta, \cdot) \rangle, \quad \text{for any } \eta \in [0, B]. 
\end{align}
This can be seen from \eqref{5.2} and the boundary conditions \eqref{5.3}--\eqref{5.4}. To obtain the weighted estimate, we multiply \eqref{5.2}$_2$ by $-e^{2\vartheta\eta} 4\mathbf{T}^3 g^B$, and also multiply \eqref{5.2}$_1$ by $e^{2\vartheta\eta} \phi^B$ and integrate over $[0, B]$ and $\mathbb{S}^2$, then we  obtain
\begin{align}\label{5.22}
	&-\int_0^B e^{2\vartheta\eta} 4\mathbf{T}^3 g^B \partial_\eta^2 g^B \mathrm{d}\eta + \int_0^B \int_{\mathbb{S}^2} \beta_3 e^{2\vartheta\eta} \partial_\eta \frac{(\phi^B)^2}{2} \mathrm{d}\beta\,\mathrm{d}\eta \nonumber\\
	&+\kappa\int_0^B \int_{\mathbb{S}^2} e^{2\vartheta\eta} (\phi^B-\bar{\phi}^B)^2 \mathrm{d}\beta\,\mathrm{d}\eta +(1-\kappa) \int_0^B \int_{\mathbb{S}^2} e^{2\vartheta\eta} (\phi^B - 4\mathbf{T}^3 g^B)^2 \mathrm{d}\beta\,\mathrm{d}\eta \nonumber\\
	&= \int_0^B \int_{\mathbb{S}^2} e^{2\vartheta\eta} (\phi^B - 4\mathbf{T}^3 g^B) S \mathrm{d}\beta\,\mathrm{d}\eta. 
\end{align}

Integrating the first term on the left of~\eqref{5.22}
by parts, using $g^B(0)=0$ and
$\partial_\eta g^B(B)=0$, and applying the finite-interval
version of~\eqref{5.0-3} established in Step 1, we obtain
\begin{align}\label{5.24-0}
	&-\int_0^B e^{2\vartheta\eta}
	4\mathbf T^3g^B\partial_\eta^2g^B\,\mathrm d\eta
	\nonumber\\
	&\quad=
	\int_0^B e^{2\vartheta\eta}
	4\mathbf T^3|\partial_\eta g^B|^2\,\mathrm d\eta
	+\int_0^B e^{2\vartheta\eta}
	\partial_\eta(4\mathbf T^3)g^B
	\partial_\eta g^B\,\mathrm d\eta
+2\vartheta\int_0^B e^{2\vartheta\eta}
	4\mathbf T^3g^B\partial_\eta g^B\,\mathrm d\eta
	\nonumber\\
	&\quad\geq
	2\vartheta\int_0^B e^{2\vartheta\eta}
	4\mathbf T^3g^B\partial_\eta g^B\,\mathrm d\eta.
\end{align}

The reflection condition makes the flux at $\eta=B$ vanish, and hence
\begin{align}\label{5.23}
	&\int_0^B \int_{\mathbb{S}^2} \beta_3 e^{2\vartheta\eta} \partial_\eta \frac{(\phi^B)^2}{2} \mathrm{d}\beta\,\mathrm{d}\eta \nonumber\\
	&= -\vartheta \int_0^B \int_{\mathbb{S}^2} e^{2\vartheta\eta} \beta_3 (\phi^B)^2 \mathrm{d}\beta\,\mathrm{d}\eta+ \frac{1}{2} \int_{\mathbb{S}^2} e^{2\vartheta\eta} \beta_3 (\phi^B)^2 \mathrm{d}\beta \bigg|_0^B \nonumber\\
	&= -\vartheta \int_0^B \int_{\mathbb{S}^2} e^{2\vartheta\eta} \beta_3 (\phi^B)^2 \mathrm{d}\beta\,\mathrm{d}\eta  - \frac{1}{2} \int_{\mathbb{S}^2} \beta_3 (\phi^B)^2(0,\cdot) \mathrm{d}\beta. 
\end{align}
The two terms produced by differentiating the exponential weight can be combined into a square involving $\phi^B-4\mathbf{T}^3g^B$. Indeed, using the moment identity \eqref{5.21-1} and $\langle\beta_3\rangle=0$, we obtain
\begin{align}\label{5.24}
2\vartheta \int_0^B e^{2\vartheta\eta} 4\mathbf{T}^3 g^B \partial_\eta g^B \mathrm{d}\eta - \vartheta \int_0^B \int_{\mathbb{S}^2} e^{2\vartheta\eta} \beta_3 (\phi^B)^2 \mathrm{d}\beta\,\mathrm{d}\eta
= -\vartheta \int_0^B \int_{\mathbb{S}^2}e^{2\vartheta\eta}\beta_3 (\phi^B - 4\mathbf{T}^3 g^B)^2 \mathrm{d}\beta\,\mathrm{d}\eta. 
\end{align}

Moreover, for the right hand side of \eqref{5.22}, we have
\begin{align*}
\int_0^B \int_{\mathbb{S}^2} e^{2\vartheta\eta} (\phi^B - 4\mathbf{T}^3 g^B) S \mathrm{d}\beta\,\mathrm{d}\eta
\leq \xi \int_0^B \int_{\mathbb{S}^2} e^{2\vartheta\eta} (\phi^B - 4\mathbf{T}^3 g^B)^2  \mathrm{d}\beta\,\mathrm{d}\eta + \frac{1}{2\xi} \int_0^B\int_{\mathbb{S}^2} e^{2\vartheta\eta} |S|^2  \mathrm{d}\beta\,\mathrm{d}\eta.
\end{align*}
Combining the above inequality with \eqref{5.24-0}, \eqref{5.23} and \eqref{5.24}, the equality \eqref{5.22} implies \eqref{5.25}. Thus we have finished the proof of Lemma \ref{lem5.1}.
\end{proof}

\medskip

\noindent{\bf Acknowledgments.}
Qiangchang Ju's research is supported by National Natural Science Foundation of China  No. 12131007.

\medskip

\noindent{\bf Conflict of Interest:} The authors declare that they have no conflict of interest.

\end{document}